\documentclass[twoside,10pt]{article}
\usepackage{graphicx}
\usepackage{subfigure}
\usepackage{ifpdf}

\pdfoutput=1

\ifpdf
      \DeclareGraphicsExtensions{.pdf,.jpg,.png}
\else
     \DeclareGraphicsExtensions{.eps}
\fi

\usepackage{amsfonts}
\usepackage{amscd}

\newcommand\Z{\ensuremath{\mathbb Z}}

\newcommand\R{\ensuremath{\mathbb R}}
\newcommand\T{\ensuremath{\mathbb T}}
\newcommand\Sc{\ensuremath{\mathbb S}}
\newcommand\Pj{\ensuremath{\mathbb P}}
\newcommand\C{\ensuremath{\mathbb C}}

\newcommand\ee{\ensuremath\mathrm{e}}

\newcommand\eps\varepsilon
\newcommand\A{\ensuremath{\mathcal A}}
\newcommand\B{\ensuremath{\mathcal B}}
\newcommand\Cc{\ensuremath{\mathcal C}}

\newcommand\E{\ensuremath{\mathcal E}}

\newcommand\Hc{\ensuremath{\mathcal H}}

\newcommand\Lc{\ensuremath{\mathcal L}}
\newcommand\M{\ensuremath{\mathcal M}}
\newcommand\Nc{\ensuremath{\mathcal N}}
\newcommand\Pc{\ensuremath{\mathcal P}}

\newcommand\Rc{\ensuremath{\mathcal R}}

\newcommand\U{\ensuremath{\mathcal U}}

\newcommand\W{\ensuremath{\mathcal W}}

\newcommand\const{\ensuremath{\rm const}}

\newcommand\ds{\displaystyle}
\newcommand{\dfrac}[2]{\ensuremath\displaystyle\frac{#1}{#2}}
\newcommand{\tfrac}[2]{\ensuremath\textstyle\frac{#1}{#2}}

\newcommand\p[1]{\left(#1\right)}
\newcommand\pq[1]{\left[#1\right]}

\newcommand\po[1]{\left(#1\right)}  
\newcommand\poc[1]{\left(#1\right]} 
\newcommand\pco[1]{\left[#1\right)} 
\newcommand\scprod[2]{\left\langle#1,#2\right\rangle}
\newcommand\abs[1]{\left|#1\right|}

\newcommand\df{\ensuremath\mathrm{d}}
\newcommand\Df{{\rm D}}
\newcommand\deriv[2]{\frac{\df#1}{\df#2}}
\newcommand\dderiv[2]{\frac{\df^2#1}{\df#2^2}}
\newcommand\pd[2]{\frac{\partial#1}{\partial#2}} 
\newcommand\pdd[2]{\frac{\partial^2#1}{\partial#2^{\,2}}}
\newcommand\pddm[3]{\frac{\partial^2#1}{\partial#2\,\partial#3}}

\newcommand\tl{\tilde}
\newcommand{\wt}[1]{\ensuremath\widetilde{#1}}
\newcommand{\wh}[1]{\ensuremath\widehat{#1}}

\newcommand\Ord{{\mathcal O}}

\newcommand{\diag}{\ensuremath\mathrm{diag}}
\newcommand\ut{{\rm u}}  
\newcommand\st{{\rm s}}  
\newcommand\Id{{\rm Id}}

\newcommand\vect[2]{\p{\begin{array}{c}#1\\[3pt]#2\end{array}}}
\newcommand\mmatrix[4]{\p{\begin{array}{cc}#1&#2\\[3pt]#3&#4\end{array}}}
\newcommand\symmatrix[3]{\mmatrix{#1}{#2}{#2}{#3}}

\newtheorem{theorem}{Theorem}
\newtheorem{proposition}[theorem]{Proposition}
\newtheorem{lemma}[theorem]{Lemma}

\newcommand\bremark{\noindent{\bf Remark.}\ \ }
\newcommand\eremark{\bigskip}
\newcommand\proof{\noindent\emph{Proof.}\ \ }
\newcommand\proofof[1]{\noindent\emph{Proof of #1.}\ \ }
\newcommand\qed{\ \ \null\nolinebreak\hfill$\frame{\large\phantom a}$}
\newcommand\hyp[1]{{\bf(#1)}}
\newcommand\paragr[1]{\noindent\textbf{\emph{#1.}}\quad}

\newcommand\beq{\begin{equation}}
\newcommand\eeq{\end{equation}}
\newcommand\bea{\begin{eqnarray}}
\newcommand\eea{\end{eqnarray}}
\newcommand\bean{\begin{eqnarray*}}
\newcommand\eean{\end{eqnarray*}}
\newcommand\btm{\vspace{-\baselineskip}\begin{itemize}}
\newcommand\etm{\end{itemize}\vspace{-\baselineskip}}

\def\cprime{\/$'$}
\def\bprime{\boldmath$'$}

\begin{document}

\title{Transversality of homoclinic orbits to hyberbolic equilibria in a
  Hamiltonian system, via the Hamilton--Jacobi equation}
\author{\sc
    Amadeu Delshams,
  \ Pere Guti\'errez,
  \ Juan R. Pacha
  \footnote{The authors were supported in part by
    the Spanish MICINN-FEDER grant MTM2009-06973, and
    the Catalan CUR-DIUE grant 2009SGR859.
    Besides, the author JRP was also supported by \mbox{MICINN-PN}
    (I+D+I) 2008--2011 (grant JC2009-00306).}
\\[12pt]
\parbox{11.5cm}{
  \small
  \begin{itemize}
  \item[]
    Dep. de Matem\`atica Aplicada I\\
    Universitat Polit\`ecnica de Catalunya\\
    Av. Diagonal 647, 08028 Barcelona\\
    {\footnotesize
      \texttt{amadeu.delshams@upc.edu},
      \texttt{pere.gutierrez@upc.edu},
      \texttt{juan.ramon.pacha@upc.edu}}
  \end{itemize}
}}
\date{\normalsize November 21, 2011}
\maketitle
\begin{abstract}
We consider a Hamiltonian system with 2 degrees of freedom, with a hyperbolic
equilibrium point having a loop or homoclinic orbit (or, alternatively, two
hyperbolic equilibrium points connected by a heteroclinic orbit), as a step
towards understanding the behavior of nearly-integrable Hamiltonians near
double resonances. We provide a constructive approach to study whether the
unstable and stable invariant manifolds of the hyperbolic point intersect
transversely along the loop, inside their common energy level. For the system
considered, we establish a necessary and sufficient condition for the
transversality, in terms of a Riccati equation whose solutions give the slope
of the invariant manifolds in a direction transverse to the loop. The key point
of our approach is to write the invariant manifolds in terms of generating
functions, which are solutions of the Hamilton--Jacobi equation. In some
examples, we show that it is enough to analyse the phase portrait of the
Riccati equation without solving it explicitly. Finally, we consider an
analogous problem in a perturbative situation. If the invariant manifolds of
the unperturbed loop coincide, we have a problem of splitting of separatrices.
In this case, the Riccati equation is replaced by a Mel$'$nikov potential
defined as an integral, providing a condition for the existence of a perturbed
loop and its transversality. This is also illustrated with a concrete example.
\par\vspace{12pt}
\end{abstract}

\section{Introduction}

\subsection{Setup and main results}\label{secsetup}

The study of the behavior of a Hamiltonian system near a double resonance
is one of the main difficulties related with Arnol{\cprime}d diffusion,
a phenomenon of instability in perturbations of integrable Hamiltonian systems
with more than 2~degrees of freedom. Such a behavior is usually studied
with the help of resonant normal forms. Neglecting the remainder, the normal
form can be reduced to a Hamiltonian with 2~degrees of freedom, that in general
is not integrable. As a first step towards studying the complete system
near the resonance, a good understanding of this reduced Hamiltonian
is very important, and particularly the intersections between the invariant
manifolds of equilibrium points, along a homoclinic orbit.

As a model for the reduced system, we consider a
\emph{classical Hamiltonian} with 2~degrees of freedom,
of the type kinetic energy plus potential energy,
with the 2-dimensional torus $\T^2$ as the configuration manifold
(see Section~\ref{secmotiv} for more details).
In fact, our approach will also be valid in a more general manifold.

Let $H$ be a Hamiltonian with 2~degrees of freedom,
defined on a phase space $T^*Q$,
where $Q$ is a 2-dimensional configuration manifold.
Let $q=(q_1,q_2)\in U\subset\R^2$ be some local coordinates for $Q$.
Then, we have canonical coordinates $x=(q,p)\in U\times\R^2$ for $T^*Q$,
with the standard symplectic form $\Omega=\df q\wedge\df p$,
whose associated matrix is
$J=\p{\begin{array}{cc}0&\Id\\-\Id&0\end{array}}$.
In these coordinates, our Hamiltonian takes the form
\beq\label{ham}
  H(q,p)=\tfrac12\scprod{B(q)p}p+V(q),
\eeq
with a positive definite (symmetric) matrix function $B(q)$,
and a scalar function $V(q)$,
providing the \emph{kinetic energy} and the \emph{potential} respectively.
The Hamiltonian equations are
$\dot x=X_H(x)=J\,\nabla H(x)=\p{\pd Hp,-\pd Hq}$,
namely
\beq\label{hameq}
  \dot q=B(q)p,
  \qquad
  \dot p=-\tfrac12\,\pd{}q[\scprod{B(q)p}p]-\nabla V(q),
\eeq
We assume that $B(q)$ and $V(q)$ are \emph{smooth} functions on $U$,
i.e.~they are $\Cc^r$ with $r\ge2$, or analytic.
Then, the Hamiltonian equations~(\ref{hameq}) are $\Cc^{r-1}$ or analytic.

For a given \emph{homoclinic orbit} or \emph{loop},
biasymptotic to a \emph{hyperbolic equilibrium point},
our goal is to provide a constructive approach to study
the \emph{transversality} of the \emph{invariant manifolds} along that orbit,
inside the energy level where they are contained.
In fact, we develop our aproach for the case of a \emph{heteroclinic orbit},
which makes no difference with respect to the homoclinic case.
We denote $O$, $\wt O$ two (possibly equal) hyperbolic equilibrium points,
and $\W^{\ut,\st}$, $\wt\W^{\ut,\st}$ their respective unstable and stable
invariant manifolds. Let $\gamma$ be a heteroclinic (or homoclinic) orbit,
that we assume known, connecting the two points,
i.e.~$\gamma\subset\W^\ut\cap{\wt\W}^\st$,
and we we have to study whether such intersection is transverse.

We consider an open neighborhood $\U$ of the first hyperbolic point $O$,
with coordinates $(q,p)\in U\times\R^2$ as in~(\ref{ham}).
Of course, this neighborhood $\U$ may not contain the whole orbit $\gamma$,
nor the second point $\wt O$.
For the point $\wt O$, we consider a neighborhood $\wt\U$
with coordinates $(\tl q,\tl p)\in\wt U\times\R^2$.
We assume that the neighborhoods $\U$ and $\wt\U$ have intersection,
in which the symplectic change
between the coordinates $(q,p)$ and $(\tl q,\tl p)$
is induced by a change in the configuration manifold $Q$:
\beq\label{qchange}
  \tl q=\chi(q),
  \qquad
  \tl p=\Df\chi(q)^{-\top}p
\eeq
(where $\Df\chi(q)$ stands for the Jacobian matrix of the change,
and we use the notation $A^{-\top}$ for the inverse of the transpose of
a matrix $A$).
In the intersection $\U\cap\wt\U$, we will study the transversality
between the unstable manifold $\W^\ut$ of $O$ and the stable manifold
$\wt\W^\st$ of $\wt O$ along the orbit $\gamma$.

The transversality between the invariant manifolds will be studied
in the coordinates $(q,p)$ of $\U$.
When restricted to the neighborhood $\U$, we may refer to
the \emph{`outgoing parts'} of $\gamma$ and $\W^\ut$ as
the local outgoing orbit and the local unstable manifold respectively,
before leaving (forward in time) the neighborhood $\U$.
On the other hand, since the global manifold $\wt\W^\st$
contains the whole orbit $\gamma$, it also enters (backward in time)
in the neighborhood $\U$, and will be compared with $\W^\ut$.

Let us describe our \emph{hypotheses} of the Hamiltonian~(\ref{ham}),
expressed in the coordinates $(q,p)$ of the neighborbood $\U$.
First, we assume:
\btm
\item[\hyp{H1}] the potential $V(q)$ has a nondegenerate maximum at $q=(0,0)$,
  with $V(0,0)=0$.
\etm
Thus, we are assuming that the point $O$ is the origin $(q,p)=(0,0,0,0)$,
and this hypothesis says that $O$
is a hyperbolic equilibrium point of the Hamiltonian $H$.
Notice that the orbit $\gamma$ is then contained
in the zero energy level of $H$.
For the sake of simplicity, we also assume:
\btm
\item[\hyp{H2}] the outgoing part of $\gamma$ satisfies $q_2=0$,
  with $q_1$ increasing along the orbit;
\item[\hyp{H3}] the expansion of $B(q_1,q_2)$ in $q_2$
  has no term of order one, i.e.~$\pd B{q_2}(q_1,0)=0$.
\etm
In fact, there always exist local coordinates $q=(q_1,q_2)$
such that~\hyp{H2} is satisfied, though it may be difficult to construct
them explicitly in a concrete case
(nevertheless, see the example of Section~\ref{secpend2}).
We point out that hypothesis~\hyp{H2} imposes strong restrictions
on the form of the functions $B(q)$ and $V(q)$ in~(\ref{ham}),
that define the Hamiltonian $H$ in the chosen coordinates $(q,p)$
(see Section~\ref{secrestrictions}).
Concerning~\hyp{H3}, this is not an essential hypothesis,
but we assume it just to ease the computations, for it is satisfied in
all the examples considered in this paper
(see technical remarks in Sections~\ref{secrestrictions} and~\ref{secriccati}).

Analogous hypotheses to~\hyp{H1--H3} could be formulated
in the local coordinates $\tl q=(\tl q_1,\tl q_2)$,
concerning the second point $\wt O$ and the ingoing part of $\gamma$.
Notice that imposing hypothesis~\hyp{H2} for both the coordinates
$q$ and $\tl q$ implies that the change~(\ref{qchange})
satisfies the relation $\chi(q_1,0)=(\tl q_1,0)$.
We will express this relation as $\chi_0(q_1)=\tl q_1$.
Nevertheless, if some reversibility properties can be applied,
it will be enough to impose hypotheses~\hyp{H1--H3}
only for the coordinates $q$, since they imply the analogous ones for $\tl q$.

To establish the transversality along $\gamma$,
we shall express the invariant manifolds
in terms of \emph{generating funcions}.
We show in Section~\ref{secgenerating} that,
in a neighborhood of the point $O$,
the unstable manifold $\W^\ut$ can be seen,
in the coordinates $(q,p)$,
as a \emph{graph} of the form
\beq\label{eq:generating}
  p=\nabla S^\ut(q),
  \qquad
  q\in U.
\eeq
This gradient form is closely related to the \emph{Lagrangian properties}
of the invariant manifolds.
The generating function $S^\ut(q)$ can be extended along a neighborhood of
the orbit $\gamma$ as far as the unstable manifold can be
expressed as a graph.
We will assume that the neighborhood $\U$ is such that
the form~(\ref{eq:generating}) is valid for the whole part of $\W^\ut$
inside $\U$.
The same can be done for the stable manifold $\wt\W^\st$
of the second point $\wt O$.
In the coordinates $(\tl q,\tl p)$ of the neighborhood~$\wt\U$,
this manifold also becomes a graph
\begin{equation}\label{eq:generating-wt}
  \tl p=\nabla\wt S^\st(\tl q),
  \qquad
  \tl q\in\wt U.
\end{equation}
We shall assume that the neighborhoods intersect,
and both generating functions can be extended up to this intersection:
\btm
\item[\hyp{H4}] there exists a point $P\in \gamma\cap\U\cap\wt\U$ with
  $q$-coordinates $(q_{1}^*,0)\in U$ and $(\tl{q}_{1}^*,0) =
  \chi(q_{1}^*,0)\in\wt{U}$, such that the generating functions
  in~(\ref{eq:generating}) and~(\ref{eq:generating-wt}) can be extended along
  the outgoing and ingoing parts of $\gamma$,
  up to $q_{1}^*$ and $\tl{q}_{1}^*$ respectively.
\etm
Under this hypothesis (to be discussed in the examples),
applying the symplectic change~(\ref{qchange}) to~(\ref{eq:generating-wt})
provides a generating function, that we denote $\wh S^\st(q)$,
for the global manifold $\wt\W^\st$ in the coordinates $(q,p)$:
\begin{equation}\label{eq:generating-wh}
  p=\nabla\wh S^\st(q),
  \qquad
  q\in U\cap\chi^{-1}(\wt U).
\end{equation}
It is easy to check, from the expression of the change~(\ref{qchange}),
that $\wh S^\st=\wt S^\st\circ\chi+\const$.
In this way, both manifolds $\W^\ut$ and $\wt\W^\st$ are expressed,
in their common neighborhood $\U\cap\wt\U$, in terms of
the same coordinates $(q,p)$. Hence, we can study the transversality of their
intersection along $\gamma$ by comparing the generating functions
$S^\ut(q)$ and $\wh S^\st(q)$.

More precisely, the transversality of the invariant manifolds will be studied
by comparing, inside the $3$-dimensional energy level containing them, the
\emph{slope} of the coordinate $p_{2}$ with respect $q_{2}$, which is a
transverse direction to $\gamma$. Since
in~(\ref{eq:generating}) and~(\ref{eq:generating-wh}) we have
$\ds p_{2}=\pd{S^{\ut}}{q_{2}}$ and $\ds p_{2}=\pd{\wh{S}^{\st}}{q_{2}}$
respectively, we define the functions
\beq\label{defT}
   T^{\ut}(q_{1}):= \pdd{S^{\ut}}{q_{2}}(q_{1},0), \qquad
   \wh{T}^{\st}(q_{1}):= \pdd{\wh{S}^{\st}}{q_{2}}(q_{1},0).
\eeq
Such functions are defined, respectively, in intervals
$U\cap\{q_2=0\}$ and $I:=U\cap\chi^{-1}(\wt U)\cap\{q_2=0\}$.
Actually, a relation between $\wh{T}^{\st}(q_{1})$ and the function
\[
  \wt{T}^{\st}(\tl q_{1}):= \pdd{\wt{S}^{\st}}{\tl q_{2}}(\tl q_{1},0)
\]
can be explicitly given.
In terms of such functions, our \emph{main result}
can be stated as follows.

\begin{theorem}\label{teoriccati00}
The functions $T^\ut(q_1)$ and $\wt T^\st(\tl q_1)$
are solutions of two Riccati equations whose coefficients are given
explicitly from the coefficients up to order~$2$
in the expansion with respect to~$q_2$ of $B(q)$ and $V(q)$.
Assuming that they can be extended up to the common
point referred to in hypothesis~\hyp{H4}, namely $q^*_1$ or
$\tl q^*_1=\chi_0(q^*_1)$ respectively, then
a necessary and sufficient condition for the transversality of the
invariant manifolds $\W^\ut$ and $\wt\W^\st$ along $\gamma$
is that the following inequality is fulfilled:
\beq\label{transverse}
  T^\ut(q^*_1)\ne\wh T^\st(q^*_1),
\eeq
where the function $\wh T^\st(q_1)$ can be expressed,
through the change~(\ref{qchange}), in terms of
$\wt T^\st(\tl q_1)$, $B(q_1,0)$ and $V(q_1,0)$.
The transversality is kept along the whole orbit $\gamma$.
\end{theorem}

This theorem is deduced in Section~\ref{secriccati},
where we provide explicitly the \emph{Riccati equation}
for the function $T^\ut(q_1)$ (see Theorem~\ref{teoriccati}),
and show that it has a unique solution under a suitable initial condition.
To obtain this Riccati equation, the key point
is to consider the generating function $S^\ut(q)$
as a solution of the \emph{Hamilton--Jacobi equation\/}:
\[
  H(q,\nabla S^\ut(q))=0,
\]
and to use the expansion in $q_2$ of this equation.
As far as $T^\ut(q_1)$ is bounded, the unstable manifold $\W^\ut$
admits a generating function $S^\ut(q)$ as in~(\ref{eq:generating}).
Of course, similar considerations can be formulated for the
function $\wt T^\st(\tl q_1)$, leading to an analogous Riccati equation.
Nevertheless, many examples satisfy a \emph{reversibility} condition
and the solution of the Riccati equation for the stable case
can be deduced directly from the unstable one
(see Section~\ref{secreversible}).
The transversality condition~(\ref{transverse})
can be checked from the Riccati equation,
even in some cases where it cannot be solved explicitly,
through a qualitative study of its phase portrait,
as we show in some examples.

\paragr{Organization of the paper}
The results of Section~\ref{secriccati0}
are valid for any \mbox{2-dimensional} configuration manifold $Q$, whereas in
Sections~\ref{sectorus} and~\ref{secperturbed} we restrict ourselves to the
case of a 2-torus: $Q=\T^2$.
On the other hand, we deal with a non-perturbative situation
in Sections~\ref{secriccati0} and~\ref{sectorus},
and with a perturbative situation in Section~\ref{secperturbed}.

We start in Section~\ref{secgenerating} by studying the invariant manifolds
$\W^{\ut,\st}$ in a neighborhood of the hyperbolic point $O$,
showing that they can be expressed in terms
of (local) generating functions $S^{\ut,\st}(q)$.
In Section~\ref{secgenerating-loop}, we consider the extension of the
generating function $S^\ut(q)$ along
the homoclinic/heteroclinic orbit $\gamma$,
and consider its Taylor expansion in a transverse direction.
We also discuss the transversality condition, from the comparison
of the generating functions $S^\ut(q)$ and $\wh S^\st(q)$.
In fact, we give the results for the unstable manifold of the point $O$,
whereas the results for the stable manifold of $\wt O$ would be analogous.
Expanding in $q_2$ the generating function $S^\ut(q)$,
we see that the coefficients of orders~0 and~1 are given by the orbit $\gamma$,
and the coefficient of order~2 is the function $T^\ut(q_1)$,
which allows us to give a necessary and sufficient condition
for the transversality of the invariant manifolds along $\gamma$.
The fact that the outgoing part of $\gamma$ is contained in $q_2=0$ imposes
several restrictions on the Hamiltonian and the generating functions,
which are studied in Section~\ref{secrestrictions}.
Next, we deduce in Section~\ref{secriccati} the Riccati equation
(Theorem~\ref{teoriccati}) by using that the generating functions are solutions
of the Hamilton--Jacobi equation, and we also establish the appropiate
initial condition for the Riccati equation.
As pointed out above, one should consider two Riccati equations,
for both the unstable and the stable manifolds.
Nevertheless, we show in Section~\ref{secreversible} that, if we consider
a reversible Hamiltonian, the solution of the Riccati equation for
the stable manifold can be deduced from the unstable one, and in this way
it is enough to consider one Riccati equation.
As an example, in Section~\ref{secdevaney} we revisit a Neumann's problem
considered in \cite{Devaney78},
consisting of an integrable system on $Q=\Sc^2$ with transversality of the
invariant manifolds.

In Section~\ref{secperiodic}, we reformulate the statements and results
of Section~\ref{secriccati0} for the case of the configuration
manifold $Q=\T^2$, where we can take advantage of the periodicity in $q_1$,
and $\gamma$ is a homoclinic orbit or loop.
This case is interesting in view of its relation to double resonances
in nearly-integrable Hamiltonians with more than 2~degrees of freedom
(see Section~\ref{secmotiv}).
As an application, we consider in Section~\ref{secpend1}
the example of two identical pendula connected by an interacting potential,
generalizing the results obtained in \cite{GelfreichS95} for the case of a
linear interaction. Using bounds of the solution of the Riccati equation for
this case, we provide a sufficient condition for the transversality.

The last part is devoted to the case of a perturbed Hamiltonian
$H_\eps=H+\eps H_*$. In Section~\ref{secmelnikov} we assume that,
for $\eps=0$, the unperturbed Hamiltonian $H$ has a loop $\gamma$
contained in $q_2=0$ as in the previous sections.
We consider two cases:
(a)~the unperturbed invariant manifolds intersect transversely along $\gamma$;
and (b)~the unperturbed invariant manifolds coincide, becoming a
2-dimensional separatrix containing $\gamma$ as an orbit belonging to
a 1-parameter family of loops $\bar\gamma_s$, $s\in\R$.
In Theorem~\ref{teoloopsAB}, we provide sufficient conditions ensuring,
for $\eps\ne0$ small enough,
the existence of a \emph{perturbed loop} $\gamma_\eps$
and the perturbed invariant manifolds intersect
transversely along it, although in the case~(b) we have to impose an
additional condition on the perturbation $H_*$
(we stress that the existence of a loop is assumed in the unperturbed
Hamiltonian $H$, but not in the perturbed one $H_\eps$).
Next, we show in Theorem~\ref{teomelnikov} that the additional condition
for the case~(b) can be expressed in terms of a Mel{\cprime}nikov potential,
defined in~(\ref{melnipot}--\ref{redumelnipot})
as the integral of $H_*$ along the unperturbed loops $\bar\gamma_s$,
which allows one to check the transversality in concrete cases.
Finally, in Section~\ref{secpend2} we apply these results to the an example
consisting of two pendula with different Lyapunov exponents,
plus a small interacting potential of order~$\eps$.

\subsection{Motivation}\label{secmotiv}

The main interest for the 2-d.o.f. model considered in this paper
lies in its close relation to \emph{resonant normal forms}.
For a given nearly-integrable Hamiltonian
$\Hc(\varphi,I)=h(I)+\eps f(\varphi,I)$,
with $N>2$ degrees of freedom, in action--angle variables, the
mechanism described in \cite{Arnold64} to detect instability
(Arnol{\cprime}d diffusion) is based on the connections between
invariant manifolds of $(N-1)$-dimensional hyperbolic invariant tori,
associated to simple resonances. Nevertheless, along the simple resonances
one also finds double resonances, which should be taken into account.

Let us give a brief description of resonant normal forms in this context
(for details, see for instance \cite{BenettinG86}, \cite[ch.~2]{LochakMS03}
and also \cite{DelshamsG01};
the ideas were initially developed in \cite{Nekhoroshev77}).
To study the behavior of the trajectories of $\Hc$
in the region close to a resonance of multiplicity $n$,
with $1\le n<N$ (associated to a module of resonances $\M\subset\Z^N$),
one carries out some steps of normalizing transformation,
in order to minimize the nonresonant terms of the Fourier expansion
in the angular variables $\varphi$. In this way, one obtains a symplectic
transformation $\Phi$ leading to new variables $(\psi,J)$,
in which the Hamiltonian becomes a
resonant normal form plus a small remainder: $\Hc\circ\Phi=\Gamma+R$,
where $\Gamma$ only depends on resonant combinations of angles.
By means of a linear change, we can assume that $\Gamma$ only depends
on $(\psi_1,\dots,\psi_n)=q$.
Writing $\psi=(q,\bar\psi)\in\T^n\times\T^{N-n}$ and
$J=(p,\bar J)\in\R^n\times\R^{N-n}$,
we have
$\Gamma=\Gamma(q,p,\bar J)$, $R=R(q,\bar\psi,p,\bar J)$,
and we can study the (truncated) normal form $\Gamma$
as a first approximation for the whole Hamiltonian~$\Hc$.

If we neglect the remainder, for the normal form $\Gamma(q,p,\bar J)$
we have $\dot{\bar J}=0$. Then, taking $\bar J$ as a parameter we can consider
an $n$-d.o.f. \emph{reduced normal form},
$\Gamma^0_{\bar J}(q,p)=\Gamma(q,p,\bar J)$,
and the behavior in the coordinates $(\bar\psi,\bar J)$
becomes just a set of $N-n$ rotors:
$\dot{\bar\psi}=\pd\Gamma{\bar J}$, $\dot{\bar J}=0$.
It is very important to understand the behavior of the reduced normal form
in the coordinates $(q,p)$ because this gives
a first approximation to the original $N$-d.o.f. Hamiltonian $\Hc$.

Under certain conditions, the reduced normal form $\Gamma^0_{\bar J}(q,p)$ has
equilibrium points, which provide a first approximation
for resonant $(N-n)$-dimensional invariant tori of the whole Hamiltonian $\Hc$.
In the same way, the $n$-dimensional invariant manifolds
of hyperbolic equilibrium points of the reduced normal form,
provide approximations for
the $N$-dimensional invariant manifolds associated to hyperbolic resonant tori,
also called \emph{hyperbolic KAM tori}
(whose splitting seems to be closely related to Arnol{\cprime}d diffusion).

We point out that the ``reduced'' invariant manifolds
can easily be studied in the case of a simple resonance ($n=1$),
because in this case the $1$-d.o.f. reduced normal form
is integrable and the invariant manifolds
become generically homoclinic connections.
For the whole Hamiltonian $\Hc$, the existence of homoclinic intersections
between the $N$-dimensional manifolds of hyperbolic KAM tori was shown
in \cite{Eliasson94}, \cite{DelshamsG00}.
Besides, the Poincar\'e--Mel{\cprime}nikov method can be used to measure
the splitting of the invariant manifolds in some restrictive models
(see \cite{LochakMS03}, \cite{DelshamsG04} for $N=3$).
As another related situation, the case of a loop asymptotic to a
center--center--saddle equilibrium was considered
in \cite{KoltsovaLDG05}, proving
the existence of homoclinic intersections and their
transversality for hyperbolic KAM tori, contained in the center manifold
and close to the equilibrium point.

But for a \emph{multiple resonance} ($n\ge2$),
in general the $n$-d.o.f. reduced normal form is \emph{non-integrable}
and hence the behavior of its invariant manifolds cannot be fully
understood, although it is possible to give some
partial results which can be useful for the whole Hamiltonian.
In this context, it is proved in \cite[\S1.10.2]{LochakMS03}
that, if there exists a homoclinic orbit $\gamma$ for the reduced normal form
such that the invariant manifolds $\W^{\ut,\st}$ \emph{intersect transversely}
along $\gamma$, then one can establish the existence of intersection
for the invariant manifolds of the whole Hamiltonian $\Hc$,
together with a lower bound ($N-n+1$) for the number of homoclinic orbits.
A similar result is proved in \cite{RudnevT06}, for a model
for the behavior near multiple resonances.
Concerning a double resonance ($n=2$), the case of a Hamiltonian with a loop
asymptotic to an equilibrium point with $2$~saddles and $N-2$~centers
has recently been considered in \cite{DelshamsGKP10},
showing under some restrictions the effective existence
of (transverse) homoclinic intersections associated to hyperbolic KAM tori
on the center manifold.
On the other hand, the dynamics near a double resonance has been studied in
\cite{Haller95} and \cite{Haller97}, but assuming that one of the resonances
is strong and the other one is weak
(in this situation, the reduced system is close to integrable).

Coming back to the classical model considered in Section~\ref{secsetup},
the Hamiltonian $H$ defined in~(\ref{ham})
can be seen as a particular model for the reduced normal form,
where our assumption that $B(q)$ is positive definite
can be related with a quasiconvexity condition on the integrable part $h(I)$,
and the potential $V(q)$ comes from the perturbation $f(\varphi,I)$
via the resonant normal form.
For such a model, we are assuming that
we are able to describe a homoclinic orbit $\gamma$,
and we are going to provide a condition which allows to check
the transversality of the invariant manifolds $\W^{\ut,\st}$
along~$\gamma$ in concrete examples.

We point out that the existence of homoclinic orbits for a Hamiltonian
of type~(\ref{ham}) can be established using variational methods.
It is shown in \cite{Bolotin78}
(see also \cite{BolotinR98a}, \cite{BolotinR98b})
that, if the configuration manifold $Q$ is compact (such as $Q=\T^2$)
and the potential $V(q)$ has a unique maximum point,
which is nondegenerate, then there exists a homoclinic orbit to this point.

\section{The generating functions as solutions of a Riccati equation}
\label{secriccati0}

\subsection{The local generating functions around a
  hyperbolic equilibrium point}
\label{secgenerating}

We show in this section that, in a neighbourhood $\U$ of a hyperbolic
equilibrium point $O$, the local invariant manifolds $\W^{\ut,\st}$ can be
written in terms of generating functions: $p=\nabla S^{\ut,\st}(q)$.
Although most of the results of this section are standard,
their proof is included here for the sake
of completeness and notational convenience.
At the end of the section, we also establish some local
properties of the functions $S^{\ut,\st}(q)$.

To carry out this local study, we use
the quadratic part of the Hamiltonian function~(\ref{ham}):
\begin{equation}\label{eq:expansion-H}
  H(q,p)=\tfrac12\scprod{\B p}p-\tfrac12\scprod{\A q}q+\Ord(\abs{(q,p)}^3),
\end{equation}
where $\A=-\Df^2V(0,0)$ is a symmetric matrix,
and $\B=B_0(0)=B(0,0)$ is a positive definite symmetric matrix.
First we show that, under hypothesis~\hyp{H1},
the point $O$ at the origin of the coordinates $(q,p)$ is
a hyperbolic equilibrium.

\begin{proposition}\label{prop:V}
If the potential $V(q)$ of the Hamiltonian~(\ref{ham}) has a
nondegenerate maximum at the origin $q=(0,0)$ of the configuration space,
then $(q,p)=(0,0,0,0)$ is an equilibrium point of hyperbolic type
of the Hamiltonian.
\end{proposition}
The following lemma will be used in the proof of this proposition.

\begin{lemma}\label{lemma:positive-d}
Let $G$ and $Q$ be two real symmetric matrices with $G$ positive definite.
Then, the product $GQ$ has all its eigenvalues real and positive
if and only if $Q$ is positive definite.
\end{lemma}
\proofof{Lemma~\ref{lemma:positive-d}}
Indeed, let us consider the Cholesky decomposition $G=LL^{\top}$. It
turns out that $GQ$ and $L^{\top}QL$ are similar matrices,
since $GQ=LL^{\top}Q\sim L^{-1}LL^\top QL=L^\top QL$.
Therefore, the lemma follows at once from the application of
Sylvester's law of inertia.
\qed

\proofof{Proposition~\ref{prop:V}}
If the potential $V(q)$ has a nondegenerate maximum at $q=(0,0)$,
then the matrix $\A = -\Df^{2}V(0,0)$ is positive definite. The Lyapunov
exponents of the equilibrium point $O$ are the eigenvalues of
the differential matrix of the Hamiltonian vector field:
\begin{displaymath}
\Df X_{H}(O) = J\Df^{2}H(O) = \biggl(\begin{array}{cc}  & \B \\
\A &  \end{array}\biggr).
\end{displaymath}
To show that these eigenvalues are all real, we consider the square,
\begin{displaymath}
\Df X_{H}(O)^{2} = \biggl(\begin{array}{cc} \B\A & \phantom{0} \\
\phantom{0} & \A\B \end{array}\biggr),
\end{displaymath}
where we have $(\B\A)^{\top}=\A\B$,
since $\A$ and $\B$ are both symmetric matrices. Then,
by Lemma~\ref{lemma:positive-d} all the eigenvalues of $\B\A$ are real and
positive, for $\A$ is positive definite.
We also see that that the product matrix $\B\A$ diagonalizes.
Indeed, if we consider the Cholesky decomposition $\B =LL^{\top}$,
we see that $L^{-1}\B\A L = L^{\top}\A L$, which is a symmetric matrix
and diagonalizes, and so does $\B\A$.

Let $\lambda_{1}^{\,2}, \lambda_{2}^{\,2}$
be the eigenvalues of $\B\A$, and $v_{1}, v_{2}$ the corresponding eigenvectors.
Since $\A\B$ is the transpose of $\B\A$, it also
diagonalizes with the same eigenvalues, whereas one sees immediately that
$w_{1}=\B^{-1} v_{1}$, $w_{2}=\B^{-1}v_{2}$ constitute a basis of eigenvectors:
\begin{displaymath}
\A\B w_{k} = \A v_{k} =
\B^{-1}\B\A v_{k} = \B^{-1}\cdot\lambda^{\,2}_{k}v_{k} =
\lambda_{k}^{\,2}w_{k},
\qquad k=1,2.
\end{displaymath}
Finally, the matrix $\Df X_{H}(O)$ has $\pm\lambda_{1},\pm\lambda_{2}$ as
eigenvalues with the four associated eigenvectors
\begin{equation}\label{eq:imvects}
\varpi_{1}^{\pm}= \biggl(\begin{array}{c} v_{1}\\
  \pm\lambda_{1}w_{1} \end{array}\biggr),\qquad
\varpi_{2}^{\pm}= \biggl(\begin{array}{c} v_{2}\\
  \pm\lambda_{2}w_{2} \end{array}\biggr)
\end{equation}
respectively, as can be easily checked:
\bean
  &&
  \Df X_{H}(O)\varpi_{k}^{\pm}=
  \biggl(\begin{array}{cc} & \B \\ \A &  \end{array}\biggl)
  \biggl(\begin{array}{c} v_{k}\\
    \pm\lambda_{k}w_{k} \end{array}\biggr)
  =\biggl(\begin{array}{c} \pm\lambda_{k}v_{k}\\
    \A v_{k} \end{array}\biggr)
\\
&&\qquad
  =\biggl(\begin{array}{c} \pm\lambda_{k}v_{k}\\
  \lambda_{k}^{2}w_{k}\end{array}\biggr) =
\pm\lambda_k\varpi_{k}^{\pm}.
\eean
Thus, the origin $(q,p)=(0,0)$ is an equilibrium point of hyperbolic type
with Lyapunov exponents $\pm\lambda_1$, $\pm\lambda_2$
(real and nonvanishing), and the proof of Proposition~\ref{prop:V} is complete.
\qed

\paragr{Existence of generating functions for the (local) invariant manifolds}
We have shown in the previous proof that the matrix $\B\A$ diagonalizes,
with real and positive eigenvalues.
Hence, we can write
\[
  \B\A=M\Lambda^2M^{-1},
  \qquad
  \A\B=N\Lambda^2N^{-1},
\]
where we define
\begin{equation}\label{eq:MN}
  \Lambda:=\symmatrix{\lambda_1}{}{\lambda_2},
  \qquad
  M :=\left(\!\!\begin{array}{ll} v_{1} & v_{2}
      \end{array}\!\!\right),
  \qquad
  N :=\left(\!\!\begin{array}{ll} w_{1} & w_{2}
      \end{array}\!\!\right) =
      \B^{-1} M.
\end{equation}
We can assume, without loss of generality,
that both $\lambda_{1}$ and $\lambda_{2}$ are positive.

The eigenvectors $\varpi_{1}^{\pm}, \varpi_{2}^{\pm}$ of $DX_{H}(O)$
introduced in~(\ref{eq:imvects}) give the
linear approximation of the (local) invariant manifolds $\W^{\ut,\st}$.
To be more precise, the eigenvectors $\varpi_1^+$, $\varpi_{2}^+$
(with positive eigenvalues) give the unstable manifold $\W^\ut$,
and the eigenvectors $\varpi_1^-$, $\varpi_2^-$
(with negative eigenvalues) give the stable manifold $\W^\ut$.
Hence, up to first order, each local invariant manifold
can be parameterized as
\begin{equation}\label{eq:param-wloc}
\renewcommand{\arraystretch}{1.5}
\W^{\ut}:
\left\{\begin{array}{l}
q= M \zeta +\Ord(|\zeta|^{2}),\\
p= N\Lambda \zeta + \Ord(|\zeta|^{2}),
\end{array}\right.\qquad\qquad
\W^{\st}:
\left\{\begin{array}{l}
q= M\zeta +\Ord(|\zeta|^{2}),\\
p= -N\Lambda \zeta + \Ord(|\zeta|^{2}),
\end{array}\right.
\end{equation}
with $q=(q_{1},q_{2})$, $p=(p_{1},p_{2})$,
and parameters $\zeta=(\zeta_1,\zeta_2)$.
On the other hand, since $\det M\ne 0$, then by the implicit function theorem:
\[
\zeta = M^{-1}q + \Ord(|q|^{2}),
\]
locally, in a neighborhood of $q= (0,0)$.
Therefore, substitution into the
second equations of~(\ref{eq:param-wloc}) yields,
\begin{equation}\label{eq:par-Wloc}
   p = g^{\ut}(q) = N\Lambda M^{-1} q + \Ord(|q|^{2}),
   \quad
   p = g^{\st}(q) = -N\Lambda M^{-1} q + \Ord(|q|^{2}).
\end{equation}
Thus, the manifolds $\W^{\ut,\st}$ can be expressed locally
as the graphs. Moreover, due to the fact that $\W^{\ut,\st}$
are \emph{Lagrangian manifolds}, the restriction of the standard $1$-form
$\theta = p_{1}\,\df q_{1}+p_{2}\,\df q_{2}$ on $\W^{\ut,\st}$,
\begin{displaymath}
  \left.\theta\right|_{\W^{\ut,\st}}= g_{1}^{\ut,\st}\,\df q_{1} +
g_{2}^{\ut,\st}\,\df q_{2},
\end{displaymath}
is a closed $1$-form, and hence locally exact.
Then, we have in~(\ref{eq:par-Wloc}) that the expressions
$g^{\ut,\st}=\p{g_{1}^{\ut,\st},g_{2}^{\ut,\st}}$ are the gradients
of \emph{generating functions} $S^{\ut,\st}$
defined in a neighborhood of $q=(0,0)$
(and uniquely determined up to constants).
Besides, it is not hard to see that the generating functions
are as smooth ($\Cc^r$ or analytic) as the initial Hamiltonian $H$.

\bremark
The result given above is local, but suitable for our purposes since,
as far as the invariant manifolds $\W^{\ut,\st}$
can be expressed as graphs $p=g^{\ut,\st}(q)$,
they admit generating functions of the form $S^{\ut,\st}(q)$
beyond a small neighborhood of the origin.
Actually, it is well-known that the manifolds $\W^{\ut,\st}$,
which are asymptotic to an equilibrium point,
are \emph{exact Lagrangian} manifolds, that is
$\left.\theta\right|_{\W^{\ut,\st}}=\df S^{\ut,\st}$
with generating functions defined globally in the whole manifolds,
$S^{\ut,\st}:\W^{\ut,\st}\longrightarrow\R$
(see for instance \cite{DelshamsR97}, \cite{LochakMS03}).
Nevertheless, the global manifolds are not necessarily graphs $p=p(q)$,
and in such a case the generating functions have to be expressed in some
other variables.
\eremark

Now, we define the following symmetric $2\times 2$~matrices
\begin{equation}\label{eq:def-E}
    \E^{\ut,\st} := \Df^{2}S^{\ut,\st}(0,0).
\end{equation}
The next lemma states some useful relations between the
matrices $\A$, $\B$ and~$\E^{\ut,\st}$.

\begin{lemma}\label{lemma:teohyppoint}
For the matrices $\E^{\ut,\st}$, we have:
\btm
\item[\rm(a)] $\E^\ut\B\E^\ut=\A$;
\item[\rm(b)] $\E^\ut=N\Lambda M^{-1}$, positive definite;
\item[\rm(c)] $\E^\st=-\E^\ut$, negative definite.
\etm
\end{lemma}

\proof
We show that part~(a) comes from the expansion
of the Hamilton--Jacobi equation in $q$,
taking the second-order terms. Indeed, expanding the gradient
$p=\nabla S^{\ut,\st}(q)$ at $q=0$ we have
\begin{eqnarray}
p &=& \nabla S^{\ut,\st}(q) = \nabla\left( S^{\ut,\st}(0) + \nabla
  S^{\ut,\st}(0) q + \frac{1}{2}\left\langle \Df^{2}S^{\ut,\st}(0)q,
q\right\rangle + \Ord(|q|^{3})\right)\nonumber\\
&&\qquad\qquad = \nabla S^{\ut,\st}(0) + \Df^{2}S^{\ut,\st}(0)q +
\Ord(|q|^{2}) = \E^{\ut,\st}q + \Ord(|q|^{2}) \label{eq:exp-DS}
\end{eqnarray}
(recall that $\nabla S^{\ut,\st}(0)=0$ for the origin of the
coordinates $(q,p)$ belongs
to both $\W^{\ut,\st}$, so no constant terms might appear in their
expansions, compare also~(\ref{eq:par-Wloc})). Next, substitution
of~(\ref{eq:exp-DS}) in~(\ref{eq:expansion-H})
leads to the following second-order expansion
of the Hamilton--Jacobi equation:
\begin{eqnarray*}
H\left(q,\nabla S^{\ut}(q)\right)&=& H\left(q,\E^{\ut}q
  +\Ord(|q|^{2})\right)\\
  &=& \frac{1}{2}\left\langle\B\left(\E^{\ut}q +\Ord(|q|^{2})\right),
     \E^{\ut}q +\Ord(|q|^{2})\right\rangle -
     \frac{1}{2}\left\langle\A q, q\right\rangle + \Ord(|q|^{3})\\
  &=& \frac{1}{2}
      \left\langle\B\E^{\ut}q,\E^{\ut}q\right\rangle -
      \frac{1}{2}\left\langle\A q, q\right\rangle + \Ord(|q|^{3}) = 0,
\end{eqnarray*}
from which the desired equality of part~(a) follows at once.

Let us show~(b). The expansions~(\ref{eq:par-Wloc}) and~(\ref{eq:exp-DS}) give
$\W^{\ut,\st}$, locally, as the graph of a function;
therefore we can identify
\begin{equation}\label{eq:D2S-2}
\E^{\ut}=N\Lambda M^{-1},\qquad
\E^{\st}=-N\Lambda M^{-1}.
\end{equation}
We see from the first of~(\ref{eq:D2S-2}) that $\B\E^{\ut}= \B N\Lambda
M^{-1}= M\Lambda M^{-1}$ (see the definition of matrix
$N$ in~(\ref{eq:MN})), hence $\B\E^{\ut}\sim \Lambda$,
which has positive eigenvalues, and $\B$ is positive definite.
By Lemma~\ref{lemma:positive-d}, the matrix $\E^{\ut}$ is
positive definite as well, proving part~(b). Finally, we deduce~(c) as an
immediate consequence of the second of~(\ref{eq:D2S-2}).
\qed

\bremark
Although the results of this section have been stated for the case of
2~degrees of freedom |so the matrices they refer to are $2\times 2$~matrices|,
the same results apply for $n$-degree-of-freedom
Hamiltonians and consequently, for $n\times n$~matrices.
\eremark

\subsection{The generating functions around a homoclinic or heteroclinic orbit}
\label{secgenerating-loop}

So far we have shown the local existence of the generating functions
$S^{\ut,\st}(q)$ in a small neighborhood of the hyperbolic point $O$,
as well as the tangent planes at $O$ of the
local unstable and stable manifolds $\W^{\ut,\st}$.
Now, our aim is to study whether such generating functions can be continued
along the homoclinic/heteroclinic orbit $\gamma$,
in order to study their transversality.
Recall from hypothesis~\hyp{H2} that, in the coordinates $(q,p)$
of the neighborhood $\U$, we are assuming that
this orbit is contained in $q_2=0$.
In fact, we may consider such expansions for any Lagrangian manifold
containing an orbit satisfying~\hyp{H2},
regardless of the fact that the orbit is asymptotic to hyperbolic points.

\paragr{Expansion of the generating functions in a transverse direction}
To start, we consider the expansion in~$q_{2}$ of the potential $V(q)$
and the matrix $B(q)$ of the Hamiltonian~(\ref{ham}):
\bea
  \label{expansV}
  V(q)
  &=
  &V_0(q_1)+V_1(q_1)\,q_2-\tfrac12Y(q_1)\,q_2^{\,2}+\Ord(q_2^{\,3}),
\\
  \label{expansB}
  B(q)
  &=
  &B_0(q_1)+\tfrac12B_2(q_1)\,q_2^{\,2}+\Ord(q_2^{\,3}),
\eea
with
\[
  B_j(q_1)=\symmatrix{b_{11j}(q_1)}{b_{12j}(q_1)}{b_{22j}(q_1)},
  \quad
  j=0,2
\]
(recall that we have $B_1(q_1)\equiv0$ according to hypothesis~\hyp{H3}).
Notice that the matrices introduced in~(\ref{eq:expansion-H})
can be expressed in terms of such expansions:
\[
  \A=\symmatrix{-V''_0(0)}{-V'_1(0)}{Y(0)},
  \qquad
  \B=\symmatrix{b_{110}(0)}{b_{120}(0)}{b_{220}(0)}.
\]

Here, we study the function $S^\ut(q)$
near the hyperbolic equilibrium point $O$,
as well as its continuation along $\gamma$,
for $q_2$ close to 0 and $(q_1,q_2)\in U$,
with the help of the Taylor expansion in the variable $q_2$
(notice that we are using $q_1$ to parameterize $\gamma$
inside the neighborhood $\U$,
and $q_2$ provides a transverse direction to it).
For both generating functions $S^{\ut,\st}(q)$,
associated to the invariant manifolds $\W^{\ut,\st}$ of the point $O$,
we consider the expansions
\beq\label{expans1}
  S^{\ut,\st}(q)
  =S_0^{\ut,\st}(q_1)+S_1^{\ut,\st}(q_1)\,q_2
   +\tfrac12T^{\ut,\st}(q_1)\,q_2^{\,2}+\Ord(q_2^{\,3}).
\eeq
The term of order~0 in these expansions is determined
up to an additive constant that, to fix ideas, will be chosen in such a way
that $S_0^{\ut,\st}(0)=0$.
Notice that the matrices defined in~(\ref{eq:def-E}) are
\beq\label{matrixE}
  \E^{\ut,\st}
  =\symmatrix{(S_0^{\ut,\st})''(0)}{(S_1^{\ut,\st})'(0)}{T^{\ut,\st}(0)}.
\eeq
Analogously, we consider the generating functions for
the invariant manifolds $\wt\W^{\ut,\st}$ of the point $\wt O$,
\beq\label{expans1b}
  \wt S^{\ut,\st}(\tl q)
  =\wt S_0^{\ut,\st}(\tl q_1)+\wt S_1^{\ut,\st}(\tl q_1)\,\tl q_2
   +\tfrac12\wt T^{\ut,\st}(\tl q_1)\,\tl q_2^{\,2}+\Ord(\tl q_2^{\,3}),
\eeq
also with $\wt S_0^{\ut,\st}(0)=0$
(we show in the subsequent sections that, if a suitable symmetry occurs,
the generating functions around $\wt O$
can be deduced from the ones obtained for $O$).

Under hypothesis~\hyp{H4}, there is a common neighborhood $\U\cap\wt\U$
where we can apply the change~(\ref{qchange}), induced by $\tl q=\chi(q)$.
With this change, we can write (a piece of) the stable manifold
$\wt\W^\st$ of $\wt O$ in the coordinates $(q,p)$,
in order to compare it with with the unstable manifold $\W^\ut$ of $O$.
Applying this change to equation~(\ref{eq:generating-wt}),
we obtain the equation
$p=\Df\chi(q)^\top\nabla\wt S^\st(\chi(q))=\nabla(\wt S^\st\circ\chi)(q)$,
i.e.~equation~(\ref{eq:generating-wh}),
also with a generating function
\beq\label{qchange2}
  \wh S^\st(q)=\wt S^\st\circ\chi(q)+\const
\eeq
(in this case, the additive constant will not be taken equal to zero;
see below).
Now, we can consider an analogous expansion of the function $\wh S^\st(q)$,
for $q_2$ close to 0 and $(q_1,q_2)\in U\cap\chi^{-1}(\wt U)$,
\beq\label{expans2}
  \wh S^\st(q)
  =\wh S_0^\st(q_1)+\wh S_1^\st(q_1)\,q_2+\tfrac12\wh T^\st(q_1)\,q_2^{\,2}
   +\Ord(q_2^{\,3}),
\eeq
where the coefficients can be determined from the ones in~(\ref{expans1b}),
applying the change $\chi$.
In particular, the function $\wh T^s(q_1)$ can be determined from
$\wt S_0^\st(\tl q_1)$, $\wt S_1^\st(\tl q_1)$ and $\wt T^\st(\tl q_1)$
(for an illustration, see the example of Section~\ref{secdevaney}).
Our aim is to compare the expansions of $S^\ut(q)$ and $\wh S^\st(q)$
in their common domain.

It is important to stress that the coefficients of orders~0 and~1
in the expansions~(\ref{expans1}) and~(\ref{expans2})
are determined by the orbit $\gamma$ itself.
Indeed, if we consider the function $S^\ut(q)$,
we see from~(\ref{eq:generating}) and hypothesis~\hyp{H2} that
the orbit $\gamma$ is given, in the neighborhood $\U$, by the equations
\beq\label{loopparam}
  \gamma:
  \qquad
  p_1=\pd{S^\ut}{q_1}(q_1,0)=(S_0^\ut)'(q_1),
  \qquad
  p_2=\pd{S^\ut}{q_2}(q_1,0)=S_1^\ut(q_1)
\eeq
(notice that at the hyperbolic point $O$
we have $(S_0^\ut)'(0)=S_1^\ut(0)=0$).
Since the same can be done with the function $\wh S^\st(q)$,
we deduce that the coefficients of orders~0 and~1 for both functions coincide.
We can introduce a common notation for them:
\begin{equation}\label{eq:defS0S1}
  S_0(q_1):=S_0^\ut(q_1)=\wh S_0^\st(q_1),
  \qquad
  S_1(q_1):=S_1^\ut(q_1)=\wh S_1^\st(q_1),
\end{equation}
for any $q_1\in I$, the common interval for both generating functions.
Notice that we can choose the additive constant in~(\ref{qchange2})
in such a way no additive constant appears in
the first equality of~(\ref{eq:defS0S1})
(this is not very important because the constant does not take part in the
gradient equations, but it is useful in order to fix ideas).

\paragr{The transversality condition}
Next, with the help of the expansions introduced above
we provide the condition for the transversality of
the \mbox{2-dimensional} invariant manifolds
$\W^\ut$ and $\wt\W^\st$ along $\gamma$.
This transversality must be considered inside
the 3-dimensional \emph{energy level} $\Nc$
containing both invariant manifolds.
Since $\Nc$ is given by the equation $H=0$,
and on $\gamma$ we have $\pd H{p_1}=\dot q_1\ne0$,
by the implicit function theorem we have that, near $\gamma$,
the energy level $\Nc$ can be parameterized as $p_1=g(q_1,q_2,p_2)$.
In the coordinates $(q_1,q_2,p_2)$ of $\Nc$,
from~(\ref{eq:generating}) we see that the unstable manifold $\W^\ut$ is given
by the equation
\[
  p_2=\pd{S^\ut}{q_2}(q_1,q_2)=S_1(q_1)+T^\ut(q_1)\,q_2+\Ord(q_2^{\,2}).
\]
This says that, in the coordinates $(q_1,q_2,p_2)$,
the coefficient $T^\ut(q_1)$ provides the slope of the manifold $\W^\ut$
in the direction of $q_2$, which is transverse to $\gamma$.
Analogous considerations can be formulated for
the stable manifold $\wt\W^\st$, whose slope is given by $\wh T^\st(q_1)$.
Then, a necessary and sufficient condition for the \emph{transversality}
of $\W^\ut$ and $\wt\W^\st$ along~$\gamma$ is that
the two slopes are different for some fixed $q^*_1\in I$,
as stated in~(\ref{transverse}).

\subsection{Restrictions on the functions defining the Hamiltonian}
\label{secrestrictions}

We assumed in hypothesis~\hyp{H2} that the orbit
$\gamma$ satisfies $q_2=0$ in the neighborhood $\U$ where
the coordinates $(q,p)$ are defined.
In this section, we show that this hypothesis implies several
equalities, that will be used later, involving
the coefficients $V_0(q_1)$, $V_1(q_1)$, $(S_0^\ut)'(q_1)$ and $S_1^\ut(q_1)$,
appearing in the expansions~(\ref{expansV}) and~(\ref{expans1}).

First of all, we parameterize the ourgoing part of the orbit $\gamma$
(inside the neighbourhood $\U$) as a trajectory.
According to hypothesis~\hyp{H2}, we have
\begin{equation}\label{eq:loop}
  \gamma:
  \qquad
  q_1=q_1^0(t),\quad q_2=0,\quad p_1=p_1^0(t),\quad p_2=p_2^0(t),
  \quad t\le t_1,
\end{equation}
where $q_1^0(t)>0$ is an increasing function asymptotic to 0 as $t\to-\infty$.
Reparameterizing $\gamma$ as a function of the
coordinate~$q_1$, we obtain the functions $(S_0^\ut)'$ and $S_1^\ut$
as in~(\ref{loopparam}).
Recall that, to simplify the notation, in~(\ref{eq:defS0S1})
we have rewritten those functions as $S'_0$ and $S_1$ respectively.

We also define the function
\beq\label{defbeta}
  \beta(q_1):=\frac{\det B_0(q_1)}{b_{220}(q_1)}\;.
\eeq

\begin{proposition}\label{teorestrictions}
\ \btm
\item[\rm(a)]
The inner dynamics along $\gamma$ in the neighborhood $\U$ is given by
the differential equation $\ \dot q_1=\beta(q_1)S'_0(q_1)$.
\item[\rm(b)]
The functions in~(\ref{eq:defS0S1}) are given by
\\$\ds S'_0(q_1)=\sqrt{\frac{-2V_0(q_1)}{\beta(q_1)}}$,
\ \ $\ds S_1(q_1)=-\frac{b_{120}(q_1)}{b_{220}(q_1)}\,S'_0(q_1)$.
\item[\rm(c)]
The following equality is satisfied:
$S'_1(q_1)\,\beta(q_1)\,S'_0(q_1)=-V_1(q_1)$.
\etm
\end{proposition}

\proof
We consider the Hamiltonian equations~(\ref{hameq}),
restricted to the orbit $\gamma$,
as well as the fact that $\gamma$ is contained
in the zero energy level of the Hamiltonian:
\bea
  \label{hameqloop1}
  &&\dot q_1=b_{110}(q_1)p_1+b_{120}(q_1)p_2,
  \\
  \label{hameqloop2}
  &&0=b_{120}(q_1)p_1+b_{220}(q_1)p_2,
  \\
  \nonumber
  &&\dot p_1=-\tfrac12\scprod{B'_0(q_1)p}p-V'_0(q_1),
  \\
  \label{hameqloop4}
  &&\dot p_2=-V_1(q_1),
  \\
  \label{hamloop}
  &&H(q_1,0,p_1,p_2)
    =\tfrac12\p{b_{110}(q_1)p_1^{\,2}+2b_{120}(q_1)p_1p_2+b_{220}(q_1)p_2^{\,2}}
     +V_0(q_1)
    =0.
\eea
According to~(\ref{loopparam}),
we can replace $p_1=S'_0(q_1)$ and $p_2=S_1(q_1)$.

As a direct consequence of~(\ref{hameqloop2}),
we obtain the second equality of~(b).
Replacing it in~(\ref{hameqloop1}) and recalling
the definition of $\beta(q_1)$ in~(\ref{defbeta}), we obtain~(a):
\[
  \dot q_1
  =\p{b_{110}(q_1)-\frac{b_{120}(q_1)^2}{b_{220}(q_1)}}S'_0(q_1)
  =\beta(q_1)S'_0(q_1).
\]
Since $q_1^0(t)$ is increasing and $\beta(q_1)>0$,
we deduce that $S'_0(q_1)>0$ for $q_1>0$.
In the same way, we see that~(\ref{hamloop}) can be written in the form
\[
  \tfrac12\beta(q_1)S'_0(q_1)^2+V_0(q_1)=0,
\]
which gives the first equality of~(b).
Finally, replacing $\dot p_2=S'_1(q_1)\,\dot q_1$ in~(\ref{hameqloop4}),
we obtain~(c).
\qed

\bremark
In~(b), we can write both $S'_0(q_1)$ and $S_1(q_1)$
in terms of $V_0(q_1)$ and $B_0(q_1)$.
Inserting this in~(c), we obtain an equality
allowing us to obtain an explicit expression for $V_1(q_1)$.
In other words, the functions  $V_0(q_1)$ and $V_1(q_1)$
in~(\ref{expansV}) cannot be independent,
due to the existence of the orbit $\gamma$ that satisfies $q_2=0$.
\eremark

\bremark
To obtain~(\ref{hameqloop4}), we used that $B_1(q_1)\equiv0$
in~(\ref{expansB}), according to hypothesis~\hyp{H3}.
If this hypothesis is not assumed, after straightforward computations
there appear more terms in
the right hand side of~(\ref{hameqloop4}) and, consequently, in the formula of
item~(c) of this proposition.
\eremark

\subsection{The Riccati equation}\label{secriccati}

In this section, we show that the function $T^\ut(q_1)$,
defined in~(\ref{defT}) or~(\ref{expans1}) from the generating function
$S^\ut(q)$ of the unstable invariant manifold $\W^\ut$ of the point $O$,
is a solution of a first-order differential equation of Riccati type
in the variable $q_1$.
This could be done in the same way for the function $\wt T^\st(\tl q_1)$
associated to the stable manifold $\wt\W^\st$
of $\wt O$, obtaining an analogous Riccati equation in $\tl q_1$.
In this way, the solutions of such Riccati equations provide the slopes of
the two invariant manifolds in a transverse direction to the orbit $\gamma$.
We are assuming in our hypothesis~\hyp{H4} that
both solutions can be extended up to some $q^*_1$ and $\tl q^*_1$,
related by the change $\chi$ in~(\ref{qchange}).
This changes provides the relation between the expansions~(\ref{expans1b})
and~(\ref{expans2}), allowing us to obtain the value of $\wh T^\st(q_1)$
and compare it with $T^\ut(q_1)$ in view of the transversality
condition~(\ref{transverse}).
Thus, in principle both equations should be solved in order to decide whether
the invariant manifolds are transverse along $\gamma$.

Nevertheless, in many cases some kind of reversibility relations are fulfilled,
with an involution relating the two invariant manifolds, as well as their
generating functions (see Section~\ref{secreversible}),
and it will be enough to find the solution of only one of the
Riccati equations.
Since the examples considered in this paper
satisfy some reversibility, in this section we only deal
with the function $T^\ut(q_1)$ associated to the unstable manifold of $O$.

In order to formulate the Riccati equation for $T^\ut(q_1)$,
we define the functions
\bea
  \label{defdelta}
  \delta(q_1)
  &:=
  &b_{120}(q_1)S'_1(q_1),
\\
  \nonumber
  \alpha(q_1)
  &:=
  &Y(q_1)-b_{110}(q_1)S'_1(q_1)^2
\\
  \label{defalpha}
  &&-\tfrac12\pq{b_{112}(q_1)S'_0(q_1)^2
                 +2b_{122}(q_1)S'_0(q_1)S_1(q_1)
                 +b_{222}(q_1)S_1(q_1)^2},
\eea
and recall that $\beta(q_1)$ was defined in~(\ref{defbeta}).

\begin{theorem}\label{teoriccati}
The function $T^{\ut}(q_{1})$ in the expansion~(\ref{expans1}) of the
generating function of the unstable manifold $\W^{\ut}$,
is a solution of the Riccati equation:
\begin{equation}\label{riccati}
  \beta(q_{1})S'_0(q_{1})\,(T^{\ut})'+2\delta(q_{1})\,T^{\ut}
  + b_{220}(q_{1})\,(T^{\ut})^{2}
  =\alpha(q_{1})
\end{equation}
\end{theorem}

\proof
We use that the generating function $S^\ut(q)$ is a solution of the
\emph{Hamilton--Jacobi equation}.
Restricting the Hamiltonian~(\ref{ham}) to the unstable manifold $\W^{\ut}$,
we have an expansion
\begin{displaymath}
H(q,\nabla S^{\ut}(q)) = \Hc_{0}^{\ut}(q_{1}) + \Hc_{1}^{\ut}(q_{1})q_{2}
+\tfrac12\Hc_{2}^{\ut}(q_{1})q_{2}^{\,2}+\Ord(q_{2}^{\,3})= 0,
\end{displaymath}
and we can replace
\[
  \nabla S^\ut(q)
  =\vect{S'_0(q_1)}{S_1(q_1)}
   +\vect{S'_1(q_1)}{T^\ut(q_1)}q_2
   +\vect{\frac12(T^\ut)'(q_1)}{\frac12S^\ut_3(q_1)}q_2^{\,2}
   +\Ord(q_2^{\,3}),
\]
where $S^\ut_3(q_1)$ corresponds to the coefficient of third order in $q_2$,
in the expansion~(\ref{expans1}).
Using also~(\ref{expansV}--\ref{expansB}), we obtain
\begin{eqnarray*}
\Hc_{0}^{\ut}(q_{1})&=&\tfrac{1}{2}\left[b_{110}(q_{1})S'_{0}(q_{1})^{2}
     + 2b_{120}(q_{1})S'_0(q_{1})S_1(q_{1})
     + b_{220}(q_{1})S_1(q_1)^2\right]\\
&&+V_{0}(q_{1}),\\[4pt]
\Hc_{1}^{\ut}(q_{1})&=& b_{110}(q_{1})S'_0(q_{1})
      S'_1(q_{1})+ b_{120}(q_{1})\left[S'_0(q_{1})
        T^{\ut}(q_{1}) + S_1(q_{1})S'_1(q_{1})\right]\\
&&+ b_{220}(q_{1}) S_1(q_{1})T^{\ut}(q_{1})
  + V_{1}(q_{1}),\\[4pt]
\Hc_{2}^{\ut}(q_{1})&=& b_{110}(q_{1})
       \left[S'_0(q_{1})(T^{\ut})'(q_{1})
             + S'_1(q_{1})^{2}\right]\\
&&+2b_{120}(q_{1})
       \left[\tfrac12S'_0(q_{1})S^{\ut}_3(q_{1})
             + S'_1(q_{1})T^{\ut}(q_{1})
             + \tfrac12S_1(q_{1})(T^{\ut})'(q_{1})\right]\\
&&+b_{220}(q_{1})
       \left[S_1(q_{1}) S^{\ut}_3(q_{1})+
             T^{\ut}(q_{1})^2\right]\\
&&+\tfrac12\left[
             b_{112}(q_{1})S'_0(q_{1})^{2} +
            2b_{122}(q_{1})S'_0(q_{1})S_1(q_{1}) +
             b_{222}(q_{1})S_1(q_1)^2\right]\\
&&-Y(q_{1}).
\end{eqnarray*}
Since $H(q,\nabla S^{\ut}(q))=0$, it follows that $\Hc_{j}(q_{1})=0$, for
$j=0,1,2$.
In particular, the expression of $\Hc^\ut_{2}(q_{1})$ above yields,
after some arrangements,
\begin{eqnarray}
&&\left[b_{110}(q_{1})S'_0(q_{1}) + b_{120}(q_{1})
S_1(q_{1})\right](T^{\ut})'
+ 2b_{120}(q_{1})S'_1(q_{1})T^{\ut}+
b_{220}(q_{1})(T^{\ut})^{2}
\nonumber\\
&&\quad =Y(q_{1})-b_{110}(q_{1})S'_1(q_{1})^{2}
\nonumber\\
&&\quad\phantom{=}
  -\tfrac{1}{2}\left[b_{112}(q_{1})S'_0(q_{1})^{2}
                     +2b_{122}(q_{1})S'_0(q_{1})S_1(q_{1})
                     +b_{222}(q_{1})S_1(q_1)^2\right],
\label{eq:riccati-1}
\end{eqnarray}
where the terms multiplied by $S^\ut_3(q_1)$ vanish,
due to the equality~(\ref{hameqloop2}) in
the proof of Proposition~\ref{teorestrictions}.
Now, taking into account~(\ref{hameqloop1}) and
Proposition~\ref{teorestrictions}(a), one has
\begin{displaymath}
b_{110}(q_{1}) S'_0(q_{1})+b_{120}(q_{1})S_1(q_{1})=\beta(q_{1})S'_0(q_{1})
\end{displaymath}
which, together with the definitions in~(\ref{defbeta}) of
$\delta(q_{1})$ and $\alpha(q_{1})$, gives rise to the Riccati
equation~(\ref{riccati}) after substitution in~(\ref{eq:riccati-1}).
\qed

\bremark
Recall that hypothesis~\hyp{H3} is assumed throughout the computations.
Otherwise, one reaches the same form~(\ref{riccati}) of the Riccati
equation, but its coefficients $\delta(q_1)$ and $\alpha(q_1)$
hold additional terms coming from the linear part $B_1(q_1)$
of the expansion~(\ref{expansB}) in $q_2$, of the matrix $B(q)$.
\eremark

The proof of Theorem~\ref{teoriccati} almost completes
the \emph{proof of Theorem~\ref{teoriccati00}}.
To finish it, it is enough to formulate the analogous Riccati equation
for $\wt T^\st(\tl q_1)$, and take into account the considerations
of Section~\ref{secgenerating-loop} on the
transversality condition~(\ref{transverse}), and on the relation between
the functions $\wt T^\st(\tl q_1)$ and $\wh T^\ut(q_1)$
(the latter one taking part in~(\ref{transverse})).
Concerning this relation, recall that the function $\wh T^s(q_1)$
can be determined from
$\wt S_0^\st(\tl q_1)$, $\wt S_1^\st(\tl q_1)$ and $\wt T^\st(\tl q_1)$.
Now, according to the first remark after Proposition~\ref{teorestrictions}
(which applies also to the stable manifold), we see that $\wh T^s(q_1)$
can be determined from $V_0(q_1)$, $B_0(q_1)$ and $\wt T^\st(\tl q_1)$,
which is the last assertion of Theorem~\ref{teoriccati00}.

Coming again to the Riccati equation~(\ref{riccati}),
it is clear that it has a singularity at $q_1=0$ since $S'_0(0)=0$.
Thus, in principle the existence and uniqueness of solution
might not take place.
However, we are going to establish the right \emph{initial condition}
for the equation, and show that the solution is unique.

\begin{lemma}\label{teoinitialcond}
The function $T^\ut(q_1)$ solving~(\ref{riccati})
satisfies the initial condition
\beq\label{initialcond}
  \ds T^\ut(0)
  =\frac{-\delta(0)+\sqrt\Delta}{b_{220}(0)}
  =\frac{\alpha(0)}{\delta(0)+\sqrt\Delta}\;,
  \qquad
  \Delta=\delta(0)^2+b_{220}(0)\alpha(0).
\eeq
\end{lemma}

\proof
Notice that a solution of~(\ref{riccati}) defined at $q_1=0$
satisfies the equality
\[
  2\delta(0)\,T^\ut(0)+b_{220}(0)\,T^\ut(0)^2=\alpha(0),
\]
and hence the initial condition is almost determined
by the differential equation:
$T^\ut(0)=(-\delta(0)\pm\sqrt\Delta)/b_{220}(0)$,
where only the sign `$\pm$' has to be determined.

Taking into account that $S'_0(0)=S_1(0)=0$,
we deduce from Proposition~\ref{teorestrictions}(b) that
\[
  S'_1(0)=-\frac{b_{120}(0)}{b_{220}(0)}\,S''_0(0).
\]
We know from Lemma~\ref{lemma:teohyppoint}(b)
that the matrix $\E^\ut$ in~(\ref{matrixE}) is positive definite.
Then, we have $S''(0)>0$
and
\bean
  &&0<S''_0(0)T^\ut(0)-S'_1(0)^2
  =S''_0(0)\pq{T^\ut(0)+\frac{b_{120}(0)}{b_{220}(0)}\,S'_1(0)}
\\
  &&\phantom{0<S''_0(0)T^\ut(0)-S'_1(0)^2}
  =S''_0(0)\pq{T^\ut(0)+\frac{\delta(0)}{b_{220}(0)}},
\eean
and putting these formulas together, it turns out that we have to choose
the positive sign in~(\ref{initialcond}).
\qed

\bremark
Using similar arguments, it is also easy to show that $\Delta>0$.
However, this is not necessary because the results
of Section~\ref{secgenerating} ensure the existence of the generating
function $S^\ut(q)$ around the origin and, subsequently,
a real value for $T^\ut(0)$.
\eremark

\paragr{Uniqueness of solution}
Despite the singularity of~(\ref{riccati}) at $q_1=0$,
we show in the next proposition that there exists a unique solution
satisfying the initial condition~(\ref{initialcond}),
which gives the function $T^\ut(q_1)$
associated to the unstable manifold $\W^\ut$.

Using the change of variable $q_1=q_1^0(t)$, provided by
the orbit $\gamma$ in the neighborhood $\U$,
we obtain another useful expression for the Riccati equation~(\ref{riccati}),
in terms of the time variable $t$.
For any given function $f(q_1)$, we write $\bar f(t)=f(q_1^0(t))$.
Then, we see from Proposition~\ref{teorestrictions}(a) that
the Riccati equation~(\ref{riccati}) becomes
\[
  \dot{\bar T}^\ut
  +2\bar\delta(t)\,\bar T^\ut
  +\bar b_{220}(t)\,(\bar T^\ut)^2
  =\bar\alpha(t).
\]

\begin{proposition}\label{teouni}
The Riccati equation~(\ref{riccati}) has a unique solution
that satisfies the initial condition~(\ref{initialcond}).
\end{proposition}

\proof
Clearly, the unstable manifold $\W^\ut$ exists,
and its generating function provides
a solution $T^\ut$ of the Riccati equation~(\ref{riccati}) with
the initial condition~(\ref{initialcond}).
Let us show the uniqueness of solution.
If $T(q_1)$ is another solution, the difference $U(q_1)=T(q_1)-T^\ut(q_1)$
is a solution of the associated Bernoulli equation, with $U(0)=0$.
In terms of $t$, the Bernoulli equation for $\bar U(t)$ becomes
\beq\label{bernoulli}
  \dot{\bar U}+\bar\psi(t)\,\bar U
    +\bar b_{220}(t)\,\bar U^2=0,
  \qquad
  \lim_{t\to-\infty}\bar U(t)=0,
\eeq
where we define
\[
  \bar\psi(t)=2(\bar\delta(t)+\bar b_{220}(t)\bar T^\ut(t)).
\]
Denoting $\ds\psi_0=\lim_{t\to-\infty}\bar\psi(t)$,
we deduce from~(\ref{initialcond}) that
\[
  \psi_0=2(\delta(0)+b_{220}(0)T^\ut(0))=2\sqrt\Delta>0.
\]

We prove the uniqueness of the solution $\bar U(t)\equiv0$
for~(\ref{bernoulli}) in a very simple way.
The idea is that the linearization of~(\ref{bernoulli})
tends, as $t\to-\infty$, to the equation $\dot{\bar U}+\psi_0\bar U=0$,
with $\psi_0>0$ (i.e.~the origin is unstable as $t\to-\infty$).
In fact, similar results for higher dimensions were established
in \cite[ch.~13]{CoddingtonL55}.
Denoting $\bar u(t)=\frac12\bar U(t)^2\ge0$, we have:
\[
  \dot{\bar u}(t)
  =\bar U(t)\dot{\bar U}(t)
  =[-\psi_0-(\bar\psi(t)-\psi_0)-\bar b_{220}(t)\bar U(t)]\cdot2\bar u(t).
\]
There exists $t_1<0$ such that
$\abs{\bar\psi(t)-\psi_0}<\psi_0/2$ and
$\abs{\bar b_{220}(t)\bar U(t)}<\psi_0/2$ for any $t<t_1$.
Then,
\[
  \dot{\bar u}(t)\le[-\psi_0+\psi_0/2+\psi_0/2]\cdot2\bar u(t)=0,
  \qquad
  t<t_1,
\]
and we see that it is not possible to have $\ds\lim_{t\to-\infty}\bar u(t)=0$,
unless $\bar u(t)\equiv0$.
\qed

\bremark
One could expect that each solution of the Riccati equation~(\ref{riccati})
generates a 2-dimensional manifold consisting of a 1-parametric family
of trajectories, containing $\gamma$.
Only one of such manifolds, namely the one associated to the initial
condition~(\ref{initialcond}), is the unstable manifold~$\W^\ut$.
\eremark

To end this section, we point out that one could use the
\emph{variational equations} around the orbit $\gamma$
as an alternative method in order to describe the invariant manifolds
$\W^{\ut,\st}$. Such a method is followed (for a more particular case)
in \cite{GelfreichS95}, \cite{RudnevT06}.
Since the variational equations are equivalent to a second-order linear
differential equation, the well-known relation between such linear equations
and Riccati equations via a change of variables provides a
relation between the approach using the variational equations and our
approach using the Hamilton--Jacobi equation, which leads
to a Riccati equation.

The advantatge of using the Riccati equation is that a qualitative
analysis of its phase portrait can be carried out with the help of
simple methods of dynamical systems.
Such a qualitative approach is useful when the solutions of
the variational equations cannot be obtained explicitly,
or they have complicated expressions.
In Sections~\ref{secdevaney} and~\ref{secpend1},
we illustrate with some examples the use of the Riccati equation.
First, in the example of Section~\ref{secdevaney}
we show that this method is simpler than solving explicitly
the corresponding second-order linear equation.
In the example of Section~\ref{secpend1}, the linear equation can be solved
explicitly in some particular cases, but in general it is not integrable.

\subsection{Reversibility relations}\label{secreversible}

In this section, we assume that the Hamiltonian equations~(\ref{hameq})
satisfy a reversibility condition, which relates the two invariant manifolds
$\W^{\ut,\st}$ of the hyperbolic point $O$.
We are going to show that
this reversibility implies a relation between the generating functions
$S^{\ut,\st}(q)$ introduced in~(\ref{expans1}),
and hence between the associated slopes $T^{\ut,\st}(q_1)$.
In fact, we should consider the invariant manifolds $\W^\ut$
and $\wt\W^\st$ of the points $O$ and $\wt O$ respectively,
in order to study whether they are transverse along
a common piece of the orbit $\gamma$.
Nevertheless, if the Hamiltonian satisfies some symmetry or some periodicity
(see Sections~\ref{secdevaney} and~\ref{secperiodic})
then it is not hard to relate the manifolds $\W^\st$ and $\wt\W^\st$.
Consequently, we formulate the results of this section for the invariant
manifolds $\W^{\ut,\st}$, in the neighborhood $\U$ of the first point $O$.

It will be enough, for the examples considered,
to consider the following type of reversibility.
We say that $H$ is \emph{$\Rc$-reversible} if
the Hamiltonian equations~(\ref{hameq}) are reversible
with respect to the linear involution
\beq\label{revers}
  \Rc:(q,p)\mapsto(Rq,-Rp),
\eeq
with a given matrix
\beq\label{matrixrevers}
  R=\symmatrix{r_1}{}{r_2},
  \qquad
  r_1,r_2=\pm1.
\eeq

\begin{lemma}\label{teoreverscond}
The Hamiltonian $H$ in~(\ref{ham}) is $\Rc$-reversible if
the functions $B(q)$ and $V(q)$ satisfy the following identities:
\[
  B(Rq)=R\,B(q)\,R,
  \qquad
  V(Rq)=V(q).
\]
\end{lemma}

\proof
It is well-known that the reversibility with respect to $\Rc$ is equivalent to
the identity $X_H\circ\Rc=-\Rc\,X_H$.
Using that $\Rc^\top J\Rc=-J$ (i.e.~$\Rc$ is `antisymplectic'),
the previous identity becomes $\nabla(H\circ\Rc)=\nabla H$,
which can be written as $H\circ\Rc=H+\const$.
But the point $O$ (the origin of the coordinates $x=(q,p)$)
is a fixed point for $\Rc$, hence the functions $H\circ\Rc$ and $H$
must coincide at this point and the constant vanishes.

Now, we have the equality $H(Rq,-Rp)=H(q,p)$. To finish the proof,
it is enough to write $H$ in terms of the functions $B(q)$ and $V(q)$
as in~(\ref{ham}).
\qed

\bremark
This condition for the $\Rc$-reversibility of the Hamiltonian~(\ref{ham})
is always satisfied if we choose
$r_1=r_2=1$ in~(\ref{matrixrevers}), which can be called
the \emph{trivial} reversibility. This can be enough in some examples,
but in other cases we will be interested in the other
types of reversibility.
\eremark

Clearly, the reversibility gives a relation between the invariant manifolds
of the hyperbolic point $O$: we have $\W^\st=\Rc\,\W^\ut$.
In the next proposition, we deduce from this fact a relation between
the generating functions associated to the invariant manifolds.

\begin{proposition}\label{prop:reversibility}
If the Hamiltonian $H$ in~(\ref{ham}) is $\Rc$-reversible,
then the generating functions defined in~(\ref{expans1})
satisfy the identity
\[
  S^\st(q)=-S^\ut(Rq)
\]
and, for the coefficients of their expansions in $q_2$,
\[
  S^\st_0(q_1)=-S^\ut_0(r_1q_1),
  \qquad
  S^\st_1(q_1)=-r_2S^\ut_1(r_1q_1),
  \qquad
  T^\st(q_1)  =-T^\ut(r_1q_1).
\]
\end{proposition}

\proof
Recall that the manifold $\W^\ut$ is given by the equation
$p=\nabla S^\ut(q)$. Applying the reversibility~(\ref{revers}),
we obtain for $\W^\st$ the equation $-Rp=\nabla S^\ut(Rq)$,
which must coincide with $p=\nabla S^\st(q)$.
Then, we have the equality
$\nabla S^\st=-R\,\nabla S^\ut\circ R=-\nabla(S^\ut\circ R)$,
which implies that $S^\st=-S^\ut\circ R+\const$.
Since in~(\ref{expans1}) we set $S_0^\ut(0)=S_0^\st(0)=0$,
it turns out that the constant vanishes.

Expanding in $q_2$ and taking into account the form of the matrix $R$
in~(\ref{matrixrevers}), we obtain as a simple consequence
the equalities involving the functions $S^{\ut,\st}_0(q_1)$,
$S^{\ut,\st}_1(q_1)$ and $T^{\ut,\st}(q_1)$.
\qed

\subsection{Example: an integrable system on $Q=\Sc^2$}\label{secdevaney}

We consider in this section the classical Neumann problem
on the two $2$-sphere $\Sc^{2}$
(for an account see~\cite{Moser80}, and the same problem is
tackled in~\cite{Devaney78}, but working in the projective space $\Pj^{2}$).
This is an example of integrable Hamiltonian system with hyperbolic
equilibrium points whose invariant manifolds intersect transversely
along heteroclinic orbits (or homoclinic orbits if one works in $\Pj^2$).
We are going to obtain this result of transversality
applying our results on the generating functions of the invariant manifolds.

For 2~d.o.f., consider a particle moving on a 2-sphere $\Sc^{2}$ in $\R^{3}$
under the action of the following potential on $\Sc^{2}$,
\begin{equation}\label{eq:V(x)}
V(x)=-\frac{1}{2}\scprod{Ax}{Ax}=-\frac{1}{2}\abs{Ax}^{2},
\end{equation}
where $A$ is a diagonal matrix, $A=\diag[\lambda_{1}, \lambda_{2},
\lambda_{3}]$. Further, it is assumed that $\lambda_{3}=0$, and that
$0<\lambda_{1} < \lambda_{2}$. As usual, $T\Sc^{2}$ denotes the
tangent bundle,
\begin{displaymath}
T\Sc^{2}=\left\{(x,\dot{x})\in\R^{3}\times\R^{3} : \abs{x}^{2}=1,\;
\mathrm{ and }\; \scprod{x}{\dot{x}}=0\right\}.
\end{displaymath}
Taking into account that $\nabla V(x) = -A^{\top}A x= -A^{2}x$, the
movement of particle on the sphere is described by the second-order equation in
$\R^{3}$,
\begin{equation}\label{eq:flux-e}
   \ddot{x} + \abs{\dot{x}}^{2}x= A^{2}x - \abs{Ax}^{2} x,
\end{equation}
whose orbits lie on $\Sc^{2}$.
It is not hard to check that
this system has 2~hyperbolic equilibrium points
at $(0,0,-1)$ and $(0,0,1)$, with 4~heteroclinic orbits connecting them.
In this section, our heteroclinic orbit $\gamma$
will be the one that goes along the semi-circle $x_2=0$, $x_1>0$.

\paragr{The Hamiltonian}
To tackle the equation~(\ref{eq:flux-e}), we introduce the two (local) charts:
$\varphi:\R^{2}\longrightarrow U :=\Sc^{2}\setminus\{(0,0,1)\}$,
given by
\begin{equation}\label{eq:psi}
 q=(q_{1}, q_{2})\mapsto \varphi(q_{1}, q_{2})=
 \left(\frac{4 q_{1}}{4 + q_{1}^{\,2} + q_{2}^{\,2}},
   \frac{4 q_{2}}{4 + q_{1}^{\,2} + q_{2}^{\,2}},
  \frac{q_{1}^{\,2} + q_{2}^{\,2}-4}{4 + q_{1}^{\,2} + q_{2}^{\,2}}\right),
\end{equation}
and $\tl{\varphi}:\R^{2}\longrightarrow \wt{U} :=\Sc^{2}\setminus\{(0,0,-1)\}$,
given by
\[
 \tl{q}=(\tl{q}_{1},\tl{q}_{2})\mapsto
 \tl{\varphi}(\tl{q}_{1},\tl{q}_{2})=
 \left(\frac{4 \tl{q}_{1}}{4+\tl{q}_{1}^{\,2}+\tl{q}_{2}^{\,2}},
   \frac{4\tl{q}_{2}}{4+\tl{q}_{1}^{\,2}+\tl{q}_{2}^{\,2}},
  -\frac{\tl{q}_{1}^{\,2}+\tl{q}_{2}^{\,2}-4}{4+\tl{q}_{1}^{\,2}+
        \tl{q}_{2}^{\,2}}\right);
\]
so $(U,\varphi)$ and $(\wt{U},\tl{\varphi})$ are the usual \emph{stereographic}
coordinates in $\Sc^{2}$. These coordinates induce the corresponding two charts
$(\U,\varphi^*)$ and $(\wt{\U},\tl{\varphi}^*)$
in the cotangent bundle $T^*\Sc^{2}$,
\bean
   &&(q,p)\in\R^{2}\times\R^{2}\mapsto \varphi^*(q,p)=
     (\varphi(q),\Df\varphi(q)^{-\top}p)\in\U=U\times\R^{2};
\\
   &&(\tl{q},\tl{p})\in\R^{2}\times\R^{2}\mapsto
   \tl{\varphi}^*(\tl{q},\tl{p})=
     (\tl{\varphi}(\tl{q}),\Df\tl{\varphi}(\tl{q})^{-\top}\tl{p})
  \in\wt{\U}=\wt{U}\times\R^{2}.
\eean
In particular, for the origin we have that
$(q,p)=(0,0)=:O \in\R^{2}\times\R^{2}$ maps to $\Pc:=\varphi^*(0,0)
=(0,0,-1,0,0,0)$ in the chart $(\U,\varphi^*)$ and
$(\tl{q},\tl{p})=(0,0)=:\wt{O}\in\R^{2}\times\R^{2}$ maps to
$\wt{\Pc}:=\tl{\varphi}^*(0,0)=
(0,0,1,0,0,0)$, in the chart $(\wt{\U},\tl{\varphi}^*)$.

On the other hand, we can give explicitly the change of coordinates or
overlapping map between the two charts
$(U,\varphi)$ and $(\wt{U},\tl{\varphi})$.
This map $\chi:\Sigma\longrightarrow\wt\Sigma$,
with $\Sigma:= \varphi^{-1}(U\cap \wt{U})=\R^{2}\setminus\{(0,0)\}$ and
$\wt{\Sigma}:=\tl{\varphi}^{-1}(U\cap\wt{U})=\R^{2}\setminus\{(0,0)\}$,
writes
\beq\label{chidev}
   (q_{1},q_{2})\mapsto (\tl{q}_{1},\tl{q}_{2})=\chi(q_{1},q_{2}):=
   \tl{\varphi}^{-1}\circ\varphi (q_{1},q_{2})=
   \left(\frac{4 q_{1}}{q_{1}^{\,2}+q_{2}^{\,2}},
   \frac{4 q_{2}}{q_{1}^{\,2} + q_{2}^{\,2}}\right).
\eeq
Now, if we denote
$q=(q_{1},q_{2})\in\Sigma$, $\wt{q}=
(\tl{q}_{1},\tl{q}_{2})\in \wt{\Sigma}$; $p=(p_{1},p_{2})\in
T^*_{q}\Sigma$, $\tl{p}=(\wt{p}_{1},\wt{p}_{2})\in
T^*_{\tl{q}}\wt{\Sigma}$, then, the extension of $\chi$ to the cotangent
bundles, $\chi^*$, can be written as
\begin{equation}\label{eq:overlap-cot}
    (q,p)\in T^*\Sigma\mapsto (\tl{q},\tl{p})=
         \chi^*(q,p)=(\chi(q),\Df\chi(q)^{-\top}p)\in
         T^*\wt{\Sigma}.
\end{equation}

The system~(\ref{eq:flux-e}) is Lagrangian and, when expressed in the
(intrinsic) coordinates~(\ref{eq:psi}), its Lagrangian function is
\begin{equation}\label{eq:L}
L(q,\dot{q})=\frac{1}{2}\scprod{\dot{q}}{G(q)\dot{q}} - \wh{V}(q),
\end{equation}
where $G(q)=\Df\varphi(q)^{\top}\Df\varphi(q)$ is the
metric tensor of $\Sc^{2}$ in the coordinates chosen, and $\wh{V}(q)=
V\circ\varphi(q)$ corresponds to the restriction of the
potential~(\ref{eq:V(x)}) on the $2$-sphere. Explicitly,
\begin{equation}\label{eq:V-Devaney}
G(q)=
\frac{16}{\p{4+q_{1}^{\,2}+q_{2}^{\,2}}^{2}}\left(\begin{array}{cc} 1 & 0\\
                                                    0 & 1\end{array}\right),
\qquad
\wh{V}(q)=
-\frac{8\left(\lambda_{1}^{\,2}q_{1}^{\,2}+\lambda_{2}^{\,2}q_{2}^{\,2}\right)}
{\p{4 + q_{1}^{\,2}+q_{2}^{\,2}}^{2}}.
\end{equation}
Furthermore, the system~(\ref{eq:flux-e}) can be brought into Hamiltonian form,
taking
\begin{displaymath}
p= (p_{1},p_{2})= \pd{L}{\dot{q}}(q,\dot{q})= G(q)\dot{q},
\end{displaymath}
as the actions conjugated to the coordinates $q=(q_{1},q_{2})$. Therefore,
the corresponding Hamiltonian function writes
\begin{equation}\label{eq:H-Dev}
   H(q,p)= \frac{1}{2}\scprod{B(q)p}{p} + \wh{V}(q)=
\frac{(4+q_{1}^{\,2}+q_{2}^{\,2})(p_{1}^{\,2}+p_{2}^{\,2})}{32}-
\frac{8(\lambda_{1}^{\,2}q_{1}^{\,2}+\lambda_{2}^{\,2}q_{2}^{\,2})}
{\p{4+q_{1}^{\,2}+q_{2}^{\,2}}^{2}},
\end{equation}
where
\begin{equation}\label{eq:B-Devaney}
B(q)=G(q)^{-1}=\frac{\p{4+q_{1}^{\,2}+q_{2}^{\,2}}^{2}}{16}
\biggl(\begin{array}{cc} 1 & 0 \\
                         0 & 1
       \end{array}\biggr).
\end{equation}
We see from the quadratic part of~(\ref{eq:H-Dev}),
\begin{displaymath}
   H_{2}(q,p)= \frac{p_{1}^{\,2}+p_{2}^{\,2}}{2} -
   \frac{\lambda_{1}^{\,2}q_{1}^{\,2}+\lambda_{2}^{\,2}q_{2}^{\,2}}{2},
\end{displaymath}
that the origin $(q,p)=(0,0)$, corresponding to the point
$\Pc=(0,0,-1,0,0,0)$ in the chart $(\U,\varphi^*)$ of $T^*\Sc^{2}$,
is a hyperbolic equilibrium point of the Hamiltonian flow
associated to~(\ref{eq:flux-e}),
with Lyapunov exponents $\pm\lambda_1$, $\pm\lambda_2$.
The same applies to $\wt{\Pc}=(0,0,1,0,0,0)$, which corresponds to
the point $(\tl{q},\tl{p})=(0,0)$ in the chart $(\wt{\U},\tl{\varphi}^*)$ of
$T^*\Sc^{2}$ (see the remark below).
We denote $\W^{\ut,\st}$ the (local) unstable and stable manifolds
of the point $\Pc$ and $\wt{\W}^{\ut,\st}$
the (local) unstable and stable manifolds of the point~$\wt{\Pc}$.

\bremark
The Lagrangian~(\ref{eq:L}) takes the same form in either
local coordinate system, $(U,\varphi)$ or $(\wt{U},\tl{\varphi})$. Thus, the
associated Hamiltonian writes the same irrespectively of which coordinates
$(q,p)$ or $(\tl{q},\tl{p})$ are taken in the phase space
$T^*\Sc^{2}$. Particularly, this implies that:
(i)~the point $\wt{\Pc}$,
represented by $(\wt{q},\wt{p})=(0,0)$ in the chart
$(\wt{\U},\tl{\varphi}^*)$ is also a hyperbolic equilibrium point, and
(ii)~the (local)
unstable and stable manifolds $\W^{\ut,\st}$ of the point $\Pc$ in the
coordinates $(q,p)$, and $\wt{\W}^{\ut,\st}$ of the point $\wt{\Pc}$ in
the coordinates $(\tl{q},\tl{p})$ are given, in terms of generating functions,
by $p = \nabla S^{\ut,\st}(q)$ and $\tl{p}= \nabla S^{\ut,\st}(\tl{q})$
respectively (i.e.~we have $S^{\ut,\st}=\wt S^{\ut,\st}$ as functions).
\eremark

Next, we expand the components of the matrix~(\ref{eq:B-Devaney}) in
the Hamiltonian~(\ref{eq:H-Dev}) in powers of $q_{2}$,
\begin{displaymath}
b_{i j}(q)=b_{i j 0}(q_{1}) +
b_{i j 1}(q_{1})q_{2}+\tfrac{1}{2} b_{i j 2}(q_{1})q_{2}^{\,2}+\cdots,
\qquad
i,j=1,2.
\end{displaymath}
Notice that there are no terms of degree $1$ in $q_{2}$,
i.e.~$b_{ij1}(q_{1})=0$ for $i,j=1,2$, which complies hypothesis~\hyp{H3}.
For the remaining coefficients, we have:
\begin{displaymath}
\begin{array}{lll}
\displaystyle b_{110}(q_{1})= \dfrac{\p{4+q_{1}^{\,2}}^{2}}{16},\quad &
b_{120}(q_{1})=b_{210}(q_{1})= 0,\quad &
b_{220}=\dfrac{\p{4+q_{1}^{\,2}}^{2}}{16},\\ [0.15in]
b_{112}(q_{1})= \dfrac{4+q_{1}^{\,2}}{4},\quad &
b_{122}(q_{1})=b_{212}(q_{1})= 0,\quad &
b_{222}(q_{1})= \dfrac{4+q_{1}^{\,2}}{4}.
\end{array}
\end{displaymath}
A similar expansion can be done for the restricted
potential~(\ref{eq:V-Devaney}),
\begin{displaymath}
\wh{V}(q)=V_{0}(q_{1})+V_{1}(q_{1})q_{2}-
\tfrac{1}{2}Y(q_{1})q_{2}^{\,2}+\cdots,
\end{displaymath}
with
\begin{displaymath}
V_{0}(q_{1})=-\frac{8\lambda_{1}^{\,2}q_{1}^{\,2}}{\p{4+q_{1}^{\,2}}^{2}},
\qquad
V_1(q_1)=0,
\qquad
Y(q_{1}) = \frac{16}{\p{4+q_{1}^{\,2}}^{2}}
\left(\lambda_{2}^{\,2}-\frac{2\lambda_{1}^{\,2}q_{1}^{\,2}}
{4+q_{1}^{\,2}}\right).
\end{displaymath}

\paragr{The outgoing part of the heteroclinic orbit}
Let us consider the local coordinates $(U,\varphi)$.
It can be checked out that the Hamiltonian~(\ref{eq:H-Dev}) has as a solution
a trajectory of the form~(\ref{eq:loop})
(the outgoing part of $\gamma$) with
$q^{0}(t)= (q_{1}^{0}(t),q_{2}^{0}(t))$, $p^{0}(t) = (p_{1}^{0}(t),p_{2}^{0}(t))
= B(q^{0}(t))\dot{q}^{0}(t)$, where $B(q)$ is the
matrix~(\ref{eq:B-Devaney}). Writting them down explicitly,
\beq\label{eq:loop-Dt}
  \gamma:
  \qquad
  \begin{array}{ll}
    q_{1}= q_{1}^{0}(t)= 2\ee^{\lambda_{1}t},
    &q_{2}= q_{2}^{0}(t)= 0,
  \\[2pt]
    p_{1}= p_{1}^{0}(t)
    =\dfrac{16\lambda_{1}\ee^{\lambda_{1}t}}{\p{4 +4\ee^{2\lambda_{1}t}}^{2}},
    &p_{2}= p_{2}^{0}(t) = 0,
  \end{array}
\eeq
which can be also parameterized by $q_{1}$:
\begin{equation}\label{eq:loop-Dq}
\gamma:\qquad
  q_{1}\ge 0,\qquad
  q_{2}= 0,\qquad
  p_{1}= \frac{16\lambda_{1}q_{1}}{\p{4 + q_{1}^{\,2}}^{2}},\qquad
  p_{2}= 0.
\end{equation}
The corresponding trajectory on the phase space $T^*\Sc^{2}$ connects the
point $\Pc$ with the point $\wt{\Pc}$
(as $t$ goes from $-\infty$ to $\infty$, or $q_1$ goes from 0 to $\infty$).
Furthermore, since $q_{1}$ increases
along the orbit and $q_{2}=0$, the outgoing part of $\gamma$ satisfies
hypothesis~\hyp{H2}; then, according
to~(\ref{loopparam}--\ref{eq:defS0S1}) one has, for this example,
\begin{displaymath}
   p_{1}= S'_0(q_{1}) = (S^\ut_0)'(q_{1}) =\frac{16\lambda_{1}q_{1}}{\p{4 +
     q_{1}^{\,2}}^{2}},\qquad\qquad
    p_{2}= S_{1}(q_{1})= S^\ut_1(q_{1})=0,
\end{displaymath}
and the functions defined in~(\ref{defbeta})
and~(\ref{defdelta}--\ref{defalpha}) are found to be
\begin{displaymath}
\beta(q_{1})= \frac{\p{4+q_{1}^{\,2}}^{2}}{16},\qquad
\delta(q_{1})= 0,\qquad \alpha(q_{1})=
\frac{16}{\p{4+q_{1}^{\,2}}^{2}}\left(\lambda_{2}^{\,2}-
\frac{4\lambda_{1}^{\,2}q_{1}^{\,2}}{4+q_{1}^{\,2}}\right).
\end{displaymath}

\bremark
Note that, in agreement with Proposition~\ref{teorestrictions}(a),
the outgoing part of $\gamma$ in~(\ref{eq:loop-Dt})
performs the inner dynamics
$\dot{q}_{1}=\beta({q}_{1})S'_0(q_{1}) = \lambda_{1}q_{1}$,
with the initial condition  $q_{1}(0)= 2$.
\eremark

Therefore, the Riccati equation~(\ref{riccati}) and its
initial condition~(\ref{initialcond}), for the example at hand,
are given by
\begin{equation}\label{eq:Riccati-Dq}
 \lambda_{1}q_{1}(T^{\ut})' + \frac{1}{16}\p{4+q_{1}^{\,2}}^{2}(T^{\ut})^{2}=
\frac{16}{\p{4+q_{1}^{\,2}}^{2}}\left(\lambda_{2}^{\,2}-
  \frac{4\lambda_{1}^{\,2}q_{1}^{\,2}}{4+q_{1}^{\,2}}\right),
\qquad
T^{\ut}(0)=\lambda_{2}.
\end{equation}

\begin{lemma}\label{prop:RS}
The solution of~(\ref{eq:Riccati-Dq}), with the given initial condition,
is defined for all $q_{1}>0$.
\end{lemma}

\begin{figure}[b!]
\subfigure[\label{fig:df-a}]{
\includegraphics[width=0.5\textwidth]{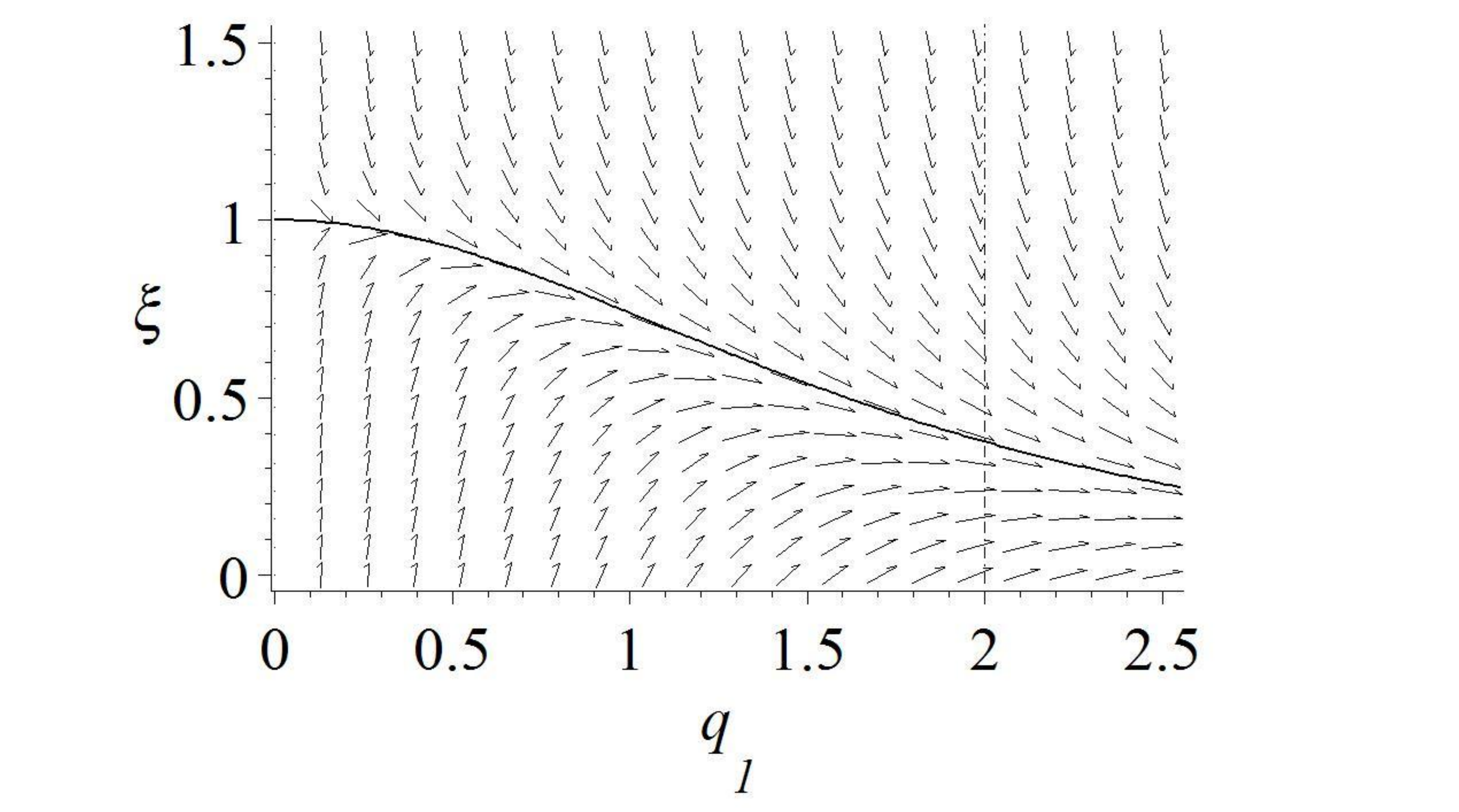}}
\hspace{-0.04\textwidth}
\subfigure[\label{fig:tu-b}]{
\includegraphics[width=0.5\textwidth]{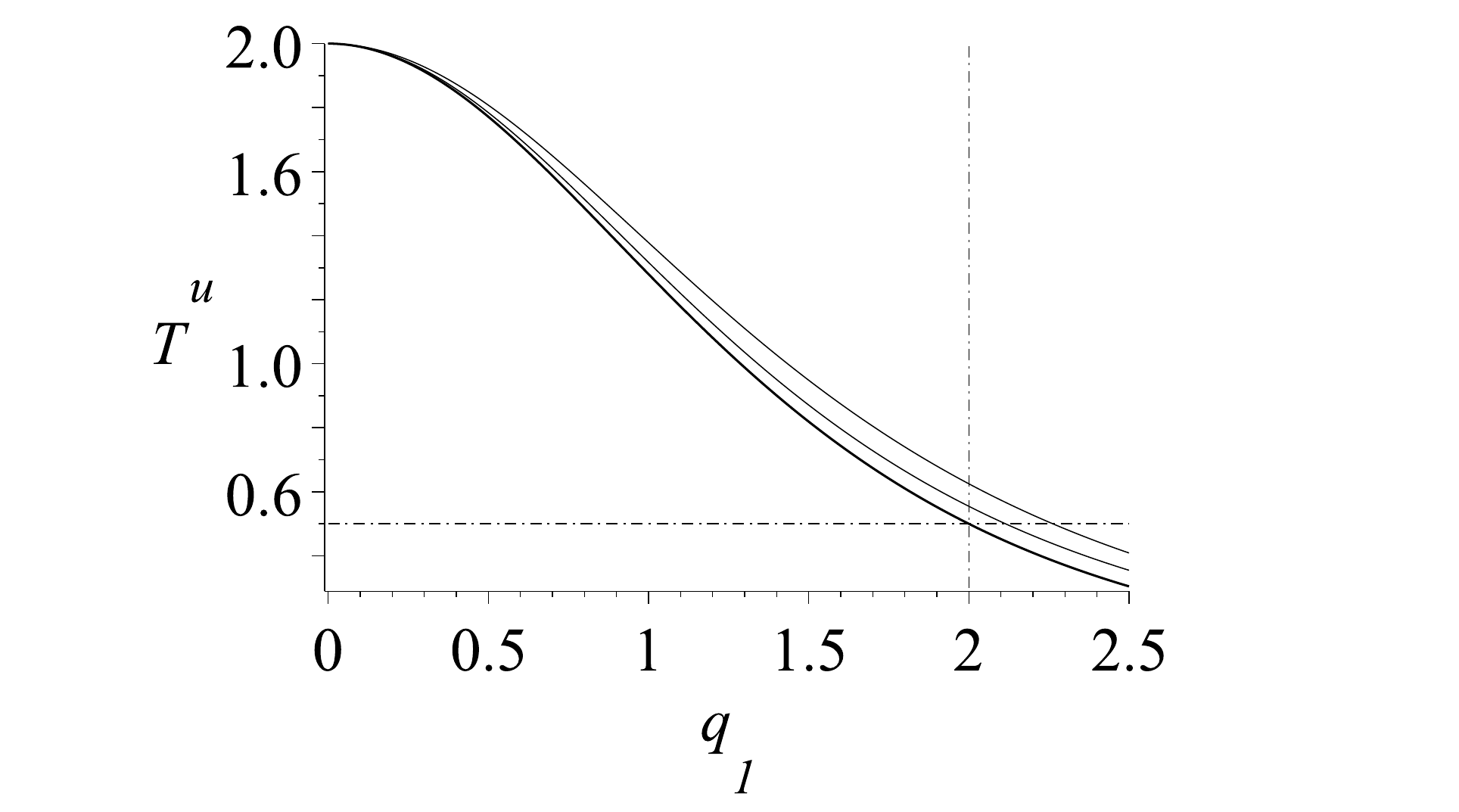}}
\caption{\small
  (a) \emph{Direction field corresponding to the Riccati
  equation~(\ref{eq:xi}), where the solution with
  $\xi(0)=\mu :=\lambda_{2}-\lambda_{1}$, for $\lambda_{2}=2$, $\lambda_{1}=1$,
  is also plotted.}
  (b) \emph{Three solutions of the Riccati
  equation~(\ref{eq:xi}) for $\lambda_{2}=2$: the lowermost one
  corresponds to $\lambda_{1}=2$ and the other two, in ascending order, to
  $\lambda_{1}= 1.75$ and $\lambda_{1}=1$, respectively. Note that the
  corresponding to $\lambda_{1}=\lambda_{2}=2$, holds the point
  $(2,\lambda_{1}/4)=(2,1/2)$.}
}
\end{figure}

\proof
It is checked out immediately that
\begin{equation}\label{eq:T0}
T_{0}^{\ut}(q_{1}):=\dfrac{16\lambda_{1}}{\p{4+q_{1}^{\,2}}^{2}}
\end{equation}
is the solution of the Riccati equation~(\ref{eq:Riccati-Dq}) for
$\lambda_{2}=\lambda_{1}$, with $T^\ut_0(0)=\lambda_{1}$. Let us denote by
$\xi(q_{1}):=T^{\ut}(q_{1})-T_{0}^{\ut}(q_{1})$ and $\mu:=\lambda_{2}-
\lambda_{1}$ ($\mu>0$, since $\lambda_{2}>\lambda_{1}$); therefore, $T^{\ut}(q)$
is a solution of~(\ref{eq:Riccati-Dq}) iff
$\xi(q_{1})$ is a solution of
\begin{equation}\label{eq:xi}
  \lambda_{1}q_{1}\xi' + 2\lambda_{1}\xi
    + \frac{1}{16}\p{4+q_{1}^{\,2}}^{2}\xi^{2}
  = \frac{16}{\p{4 + q_{1}^{\,2}}^{2}}(2\mu\lambda_{1}+\mu^{2}),
  \qquad
  \xi(0)=\mu.
\end{equation}
By Proposition~\ref{teouni}, there exists $\varepsilon > 0$,
such that $\xi(q_{1})$ is defined for $0<q_{1}<\varepsilon$.
The idea is to extend
this local solution. Consider $0<\eta<\varepsilon$ such that $\xi(\eta) >
\mu/2$. From~(\ref{eq:xi}) it follows that $\xi' > 0$ if $q_{1}>\eta$
and $\xi = 0$, i.e.~the direction field points upwards at $\xi =0$
(see Figure~\ref{fig:df-a}); so $\xi(q_{1}) > 0$ for $q_{1}\ge\eta$.
On the other hand,
\begin{eqnarray*}
   \xi' &=& \frac{16/\lambda_{1}}{q_{1}\p{4+q_{1}^{\,2}}^{2}}
   (2\mu\lambda_{1}+\mu^{2})-\frac{2}{q_{1}}\xi -
   \frac{\p{4+q_{1}^{\,2}}^{2}}{16\lambda_{1}q_{1}}\xi^{2}\\
&\le& \frac{16/\lambda_{1}}{\eta\p{4+\eta^{2}}^{2}}(2\mu\lambda_{1}+\mu^{2}) -
     \frac{2}{q_{1}}\xi - \frac{\p{4+\eta^{2}}^{2}}{16\lambda_{1}q_{1}}\xi^{2}
     \le c_{1}
\end{eqnarray*}
for $q_{1} > \eta$, with $c_{1}>0$. This means that the slope is bounded by a
positive quantity and then, this local solution, $\xi(q_{1})$ with
$\xi(0)=\mu$, can be continued, and is defined for all $q_{1}>0$.
\qed

\paragr{Using the reversibility}
One sees from Lemma~\ref{teoreverscond} that the Hamiltonian~(\ref{eq:H-Dev})
satisfies the reversibility defined in Section~\ref{secreversible},
with respect to the linear involution $\Rc:
(q,p)\mapsto (Rq, -Rp)$, with $R=\diag[r_{1},r_{2}]$, and for any
$r_{1},r_{2}=\pm 1$.
In particular,
if we take $r_{1}=r_{2}=1$, then Proposition~\ref{prop:reversibility} yields
\begin{equation}\label{eq:sym-st}
   S_{0}^{\st}(q_{1})= -S_{0}^{\ut}(q_{1})=
   -\frac{2\lambda_{1}q_{1}^{\,2}}{4+q_{1}^{\,2}},
\quad
   S_{1}^{\st}(q_{1})= -S_{1}^{\ut}(q_{1})=0,
\quad
   T^{\st}(q_{1})=-T^{\ut}(q_{1}).
\end{equation}

Let us now consider the other local chart $(\wt{\U},\tl{\varphi}^*)$. The
(local) stable manifolds $\wt{\W}^{\st}$ of the origin
$\wt{O}:=(\wt{q},\wt{p})=(0,0)$
(corresponding, in this coordinate system, to the point $\tl{P}=(0,0,1,0,0,0)$
on the phase space $T^*\Sc^{2}$), can be put as a
graph of a function through the same generating function $S^{\st}$, for the
Hamiltonian takes exactly the same expression in both local charts,
as we pointed out before.
Thus, the manifold $\wt{\W}^{\st}$ is given by
$\wt{p}=\nabla S^{\st}(\wt{q})$.
Therefore, we know from~(\ref{qchange2}) that
the change of coordinates~(\ref{eq:overlap-cot}) yields the identity
$\wh{S}^{\st}= S^{\st}\circ\chi + c$, with $c$ constant,
and this gives
$\wt{\W}^{\st}$ in the coordinates $(q,p)$ of the chart
$(\U,\varphi)$.

Expanding $\wh{S}^{\st}$ around the outgoing part of
$\gamma$ given by~(\ref{eq:loop-Dq}), regarding the
reversibility~(\ref{eq:sym-st}) and the components $\chi=(\chi_1,\chi_2)$
given in~(\ref{chidev}), one has
\begin{eqnarray*}
   \wh{S}^{\st}(q)
   &=& \wh{S}^{\st}_0(q_{1}) + \wh{S}^{\st}_1(q_{1})q_{2} +
   \frac{1}{2}\wh{T}^{\st}(q_{1})q_{2}^{\,2} +\cdots\\ [0.1in]
   &=& c - S^{\ut}_0(\chi_{1}(q)) - S^{\ut}_{1}(\chi_{1}(q))\chi_{2}(q) -
       \frac{1}{2} T^{\st} (\chi_{1}(q))\chi_{2}(q)^{2}+\cdots\\ [0.1in]
   &=& c - S_{0}^{\ut}\left(\frac{4q_{1}}{q_{1}^{\,2}+q_{2}^{\,2}}\right)
      -\frac{1}{2} T^{\ut}\left(\frac{4q_{1}}{q_{1}^{\,2}+q_{2}^{\,2}}\right)
      \frac{16q_{2}^{\,2}}{(q_{1}^{\,2}+q_{2}^{\,2})^{2}}+\cdots\\ [0.2in]
   &=& c - S_{0}^{\ut}\p{\frac{4}{q_{1}}} +
      \frac{1}{2}\left[\frac{8}{q_{1}^{\,3}}(S^{\ut}_{0})'\p{\frac{4}{q_{1}}} -
   \frac{16}{q_{1}^{\,4}} {T}^{\ut}\p{\frac{4}{q_{1}}}
 \right]q_{2}^{\,2}+\cdots,
\end{eqnarray*}
and by comparison of coefficients,
\bean
  &&\wh{S}^{\st}_0(q_{1})= c - S^{\ut}_0\p{\frac{4}{q_{1}}}=
  c-\frac{8\lambda_{1}}{q_{1}^{\,2}+4},\qquad
\\
  &&\wh{T}^{\st}(q_{1})=
  \frac{8}{q_{1}^{\,3}}(S^{\ut}_{0})'\p{\frac{4}{q_{1}}} -
  \frac{16}{q_{1}^{\,4}} T^{\ut}\p{\frac{4}{q_{1}}}=
  \frac{32\lambda_{1}}{(q_{1}^{\,2}+4)^{2}}-
  \frac{16}{q_{1}^{\,4}}T^{\ut}\p{\frac{4}{q_{1}}}
\eean
(according to~(\ref{eq:defS0S1}), the constant $c$ can be fixed if,
for example, we set $\wh{S}_{0}^{\st}(2)=S^{\ut}_0(2)=\lambda_{1}$;
this exacts $c=2\lambda_{1}$).
Then, the difference between $T^{\ut}$
and $\wh{T}^{\st}$ at $q^*_{1}=2$ is
\begin{displaymath}
T^{\ut}(2) - \wh{T}^{\st}(2) = 2T^{\ut}(2)-\frac{\lambda_{1}}{2}
\end{displaymath}
and therefore,
\beq\label{eq:transvdevaney}
  T^{\ut}(2)\ne \frac{\lambda_{1}}4
\eeq
is a necessary and sufficient condition
for transversality between $\W^{\ut}$ and $\wt{\W}^{\st}$ along the
heteroclinic orbit $\gamma$. This is stated in the proposition below.
Notice that, in the chart~(\ref{eq:psi}), for $q=(q^*_1,0)=(2,0)$
we obtain the point $(1,0,0)\in\Sc^2$.

\begin{proposition}
If in the system~(\ref{eq:flux-e}) defined on the $2$-sphere $\Sc^{2}$,
with $A$ the diagonal matrix $A=\diag[\lambda_{1},\lambda_{2},\lambda_{3}]$,
we assume
$\lambda_{2} > \lambda_{1}> \lambda_{3}=0$; then, for its associated Hamiltonian
flow in the phase space $T^*\Sc^{2}$, the stable invariant manifold
$\wt{\W}^{\st}$ of the point $\wt{\Pc}=(0,0,1,0,0,0)$ and the unstable invariant
manifold $\W^{\ut}$ of the point $\Pc=(0,0,-1,0,0,0)$ intersect transversely
along the orbit $\gamma$ (see~(\ref{eq:loop-Dq})).
\end{proposition}

\proof
By the considerations in the previous paragraph, it suffices to check
the transversality condition~(\ref{eq:transvdevaney}),
which is equivalent to $\xi(2)\ne 0$, where
$\xi(q_{1}):= T^{\ut}(q_{1})-T^\ut_{0}(q_{1})$, with $T_{0}^{\ut}(q_{1})$ given
by~(\ref{eq:T0}); but in the proof of Lemma~\ref{prop:RS} it is stated
that $\xi >0$ for $q_{1}> 0$. Thus, in particular, $\xi(2)\ne 0$ is satisfied.
This proves the proposition.
\qed

\paragr{Integrability of the Riccati equation}
Alternatively, the Riccati equation~(\ref{eq:Riccati-Dq}) could have been
solved explicitly in order to check the transversality
condition~(\ref{eq:transvdevaney}). First, a standard change of type
\begin{equation}\label{eq:C-Riccati-to-linear}
   T^{\ut}=\frac{16\lambda_{1}q_{1}}{\p{4+q_{1}^{\,2}}^{2}}\cdot\frac{y'}{y},
\end{equation}
where $T^{\ut}=T^{\ut}(q_{1})$, and $y = y(q_{1})$ is the new unknown
function, transforms~(\ref{eq:Riccati-Dq}) into the second-order linear
equation
\begin{equation}\label{eq:EDO-lineal}
y'' -\frac{3 q_{1}^{\,2}-4}{q_{1}(4+q_{1}^{\,2})}\,y'
-\frac{4\lambda_{2}^{\,2}+(\lambda_{2}^{\,2} - 4\lambda_{1}^{\,2})q_{1}^{\,2}}
{\lambda_{1}^{\,2}q_{1}^{\,2}(4+q_{1}^{\,2})}\,y=0.
\end{equation}

\bremark
This is a particular case of the well-known fact that any Riccati equation
$z'' + a(x) z + b(x)z^{2} = c(x)$
transforms into a second-order linear equation,
$y'' + \left(a(x) - \dfrac{b'(x)}{b(x)}\right)y' - b(x) c(x) y = 0$,
through the change $z = \dfrac{y'}{b(x) y}$\,.
\eremark

Now, the \emph{Kovacic's algorithm} can be applied to investigate
the (Liouvillian)
integrability of second-order differential equations with coefficients in the
class $\C(x)$ of \emph{rational functions} (for a complete description of the
process and applications, we point the reader to~\cite{Kovacic86},
\cite{DuvalL92}, \cite{AcostahMW11} and references therein).
This algorithm is implemented in
some computer algebra systems such as Maple; and when applied to
equation~(\ref{eq:EDO-lineal}), it produces the following two fundamental
solutions:
\bean
&&y_{1}(q_{1})= \left( -4(\lambda_{{2}}+\lambda_{{1}})+(\lambda_{{1}}-
\lambda_{{2}})q_{1}^{\,2} \right) q_{1}^{\lambda_{2}/\lambda_{1}},
\\
&&y_{2}(q_{1})=\frac{ 4(\lambda_{{2}}-\lambda_{{1}})+(\lambda_{{1}}+
\lambda_{2}) q_{1}^{\,2}}{\lambda_{{2}}+\lambda_{1}}\:
q_{1}^{-\lambda_{2}/\lambda_{1}};
\eean
hence, the linear equation~(\ref{eq:EDO-lineal}) is integrable and so it is
the Riccati equation~(\ref{eq:Riccati-Dq}). Actually, transforming back
$y_{1}(q_{1})$ using~(\ref{eq:C-Riccati-to-linear}), one obtains a solution
of~(\ref{eq:Riccati-Dq}),
\[
T^{\ut}(q_{1})= \frac{16\lambda_{1}q_{1}}{\p{4+q_{1}^{\,2}}^{2}}
   \cdot\frac{y'_{1}(q_{1})}{y_{1}(q_{1})}= \frac{-\left(
      16\, \lambda_{2}^{\,2}+16\,\lambda_{1}\lambda_{2}\right)
\left(q_{1}^{\,2}+4 \right) +32\, \lambda_{1}^{\,2}q_{1}^{\,2}}
{ \left( q_{1}^{\,2}+4 \right) ^{2} \left( \left( \lambda_{1}-
\lambda_{2} \right) q_{1}^{\,2}-4\,\lambda_{1}-4\,\lambda_{2} \right)},
\]
with $T^{\ut}(0)=\lambda_{2}$; and since
$\lambda_{2}>\lambda_{1}>0$, it is well defined for all $q_{1}>0$.
Besides, the transversality condition~(\ref{eq:transvdevaney})
can also be checked straightforward from this explicit
form of the solution of the Riccati
equation~(\ref{eq:Riccati-Dq}).
Indeed, for $q_{1}=2$ it gives
\begin{displaymath}
 T^{\ut}(2)= \frac{1}{4}\left(\lambda_{2}+\lambda_{1}-
   \frac{\lambda_{1}^{\,2}}{\lambda_{2}}\right)> \frac{\lambda_{1}}{4}
\end{displaymath}
for $\lambda_{2}>\lambda_{1}>0$.
In Figure~\ref{fig:tu-b}, the function $T^{\ut}$ is plotted
for different values of the parameters.

\section{The case of a 2-torus: $Q=\T^2$}\label{sectorus}

\subsection{The generating functions in the periodic case}\label{secperiodic}

In what comes, that the configuration space will be the 2-torus,
i.e.~$Q=\T^{2}=(\R/2\pi\Z)^{2}$. In view of the form of the
Hamiltonian~(\ref{ham}), this is equivalent to suppose that the quadratic form
$B(q)$ and the potential $V(q)$ are both \emph{$2\pi$-periodic} functions
with respect to $q=(q_{1},q_{2})$. In this way, the hyperbolic point
$O = (0,0,0,0)$ can be identified with $\wt{O}=(2\pi,0,0,0)$,
and therefore a biasymptotic orbit $\gamma$ connecting both ends
is a \emph{homoclinic orbit} or \emph{loop}.

Let us introduce coordinates in a neighborhood of $O$ and in a neighborhood
of $\wt{O}$, denoted $(q,p)=(q_{1},q_{2},p_{1},p_{2})$ and
$(\tl{q},\tl{p})=(\tl{q}_{1},\tl{q}_{2},\tl{p}_{1},\tl{p}_{2})$,
having the origin at $O$ and $\wt O$ respectively.
It is clear that the change~(\ref{qchange}) relating them is given by
\begin{equation}\label{eq:change-in-T2}
  (\tl q_1,\tl q_2)=\chi(q_1,q_2)=(q_1-2\pi,q_2),
  \qquad
  (\tl p_1,\tl p_2)=(p_1,p_2).
\end{equation}
since $\Df\chi(q)^{-\top}=\Id$.

\bremark
Due to the periodicity of the Hamiltonian, the corresponding Hamiltonian
equations write identically in either coordinates $(q,p)$ or $(\tl q,\tl p)$.
The invariant manifolds $\W^{\ut,\st}$ and $\wt\W^{\ut,\st}$ are related by
the periodicity and, therefore, their generating functions
$S^{\ut,\st}$ and $\wt{S}^{\ut,\st}$ are the same functions
(note however, that the functions $S^{\ut,\st}$ depend on $q$
whilst $\wt{S}^{\ut,\st}$ depend on $\tl{q}$).
\eremark

We assume the hypotheses~\hyp{H1--H4} stated in Section~\ref{secsetup},
with some minor modifications in order to adapt them to the case $Q=\T^2$
being studied here:
\btm
\item[\hyp{H2\bprime}] the orbit $\gamma$ satisfies $q_2=0$,
  with $q_1$ increasing along the orbit from $0$ to $2\pi$;
\item[\hyp{H4\bprime}]
  the generating function
  $S^{\ut}(q_{1},q_{2})$ of the invariant manifold $\W^{\ut}$ of
  $O$ is defined in a $\R^{2}$-neighborhood of the segment $0\le q_{1}\le\pi$
  of type $U = \po{-a,\pi+a}\times\po{-\delta,\delta}$,
  for some $a>0$ and
  $\delta > 0$;  and conversely, the generating function
  $\wt{S}^{\st}(\tl{q}_{1},\tl{q}_{2})$ of the invariant manifold
  $\wt{W}^{\st}$ of $\wt{O}$ is defined in a $\R^{2}$-neighborhood
  of the segment $-\pi\le \tl{q}_{1}\le 0$ of type
  $\wt{U}=\po{-\pi-a,a}\times\po{-\delta,\delta}$.
\etm

Applying the change~(\ref{eq:change-in-T2}), we can write the manifold
$\wt W^\st$ in the coordinates $(q,p)$. As in~(\ref{eq:generating-wh}),
it has a generating function, given by
$\wh S^\st=\wt S^\st\circ\chi+{\const}=S^\st\circ\chi+\const$,
that is,
\beq\label{eq:S-hat}
  \wh S^\st(q_1,q_2)=S^\st(q_1-2\pi,q_2)+\sigma.
\eeq
In other words, the generating function $\wh{S}^{\st}(q)$
is the function $\wt{S}^{\st}(\tl{q})$
expressed in the coordinates $q$
and shifted by an additive constant $\sigma$, which will be chosen below.

By hypothesis~\hyp{H4\bprime}, we see that
$S^\ut(q)$ and $\wh S^\st(q)$ can be compared
in their common domain
$U\cap\chi^{-1}(\wt U)=\po{\pi-a,\pi+a}\times\po{-\delta,\delta}$,
which is a $\R^2$-neighborhood of $q^*_1=\pi$.
As we see below, our choice of this point to check
the transversality allows us to take advantage
of the reversibility properties of the Hamiltonian.

\paragr{The generating functions along the loop}
As in the general setting described in Section~\ref{secgenerating-loop},
the loop $\gamma$ is described by the coefficients $S_0$ and $S_1$,
of orders~0 and~1, in the expansion in $q_2$ of the generating function of
any of the invariant manifolds containing it.
More precisely, the loop $\gamma$ is given by the
equations~(\ref{loopparam}--\ref{eq:defS0S1}),
which correspond to the coefficients in
the expansions~(\ref{expans1}) and~(\ref{expans2}).
Such coefficients are, initially, functions of $q_1$
defined in a finite interval.
Nevertheless, using that they must coincide for the unstable and stable
manifolds, together with the $2\pi$-periodicity in~$q_1$ of the Hamiltonian,
we are going to show that such coefficients are
$4\pi$-periodic functions of $q_1$.

Combining~(\ref{loopparam}--\ref{eq:defS0S1}) with the
identity~(\ref{eq:S-hat}), we have
\bea
  \nonumber
  \gamma:
  &&p_1=S'_0(q_1):=(S^\ut_0)'(q_1)=(\wh S^\st_0)'(q_1)=(S^\st_0)'(q_1-2\pi),
\\
  \label{eq:DS0S1}
  &&p_2=S_1(q_1):=S^\ut_1(q_1)=\wh S^\st_1(q_1)=S^\st_1(q_1-2\pi).
\eea
As a parameterization of $\gamma$, this is valid for $0<q_1<2\pi$,
but the functions in the equalities can be extended up
to a neighborhood of $0\le q_1\le2\pi$.

In fact, it is easy to see from the above equations~(\ref{eq:DS0S1})
and the reversibility relations of
Proposition~\ref{prop:reversibility} with $r_{1}=r_{2}=1$ |this is the
trivial $\Rc$-reversibility, see the remark after Lemma~\ref{teoreverscond}|,
that $S'_0$ and $S_1$ can be extended as $2\pi$-antiperiodic functions:
\[
  \begin{array}{ll}
    S'_{0}(q_{1})=(S^\st_0)'(q_1-2\pi)=-(S^\ut_0)'(q_1-2\pi)=-S'_0(q_{1}-2\pi),
  \\[4pt]
    S_{1}(q_{1})=S^\st_1(q_1-2\pi)=-S^\ut_1(q_1-2\pi)=-S_{1}(q_{1}-2\pi).
  \end{array}
\]
Thus, the functions $S'_0(q_1)$ and $S_1(q_1)$ can be defined
for all $q_{1}\in\R$, and are $2\pi$-antiperiodic, and hence they are
$4\pi$-periodic and have zero average.
It turns out that the coefficient $S_{0}(q_1)$ (without derivative)
is also a $4\pi$-periodic function. Indeed, one can easily check
that a primitive of any $2\pi$-antiperiodic function is $4\pi$-periodic.

\begin{figure}[b!]
\subfigure[\label{fig:V0}]{
\includegraphics[width=0.4\textwidth]{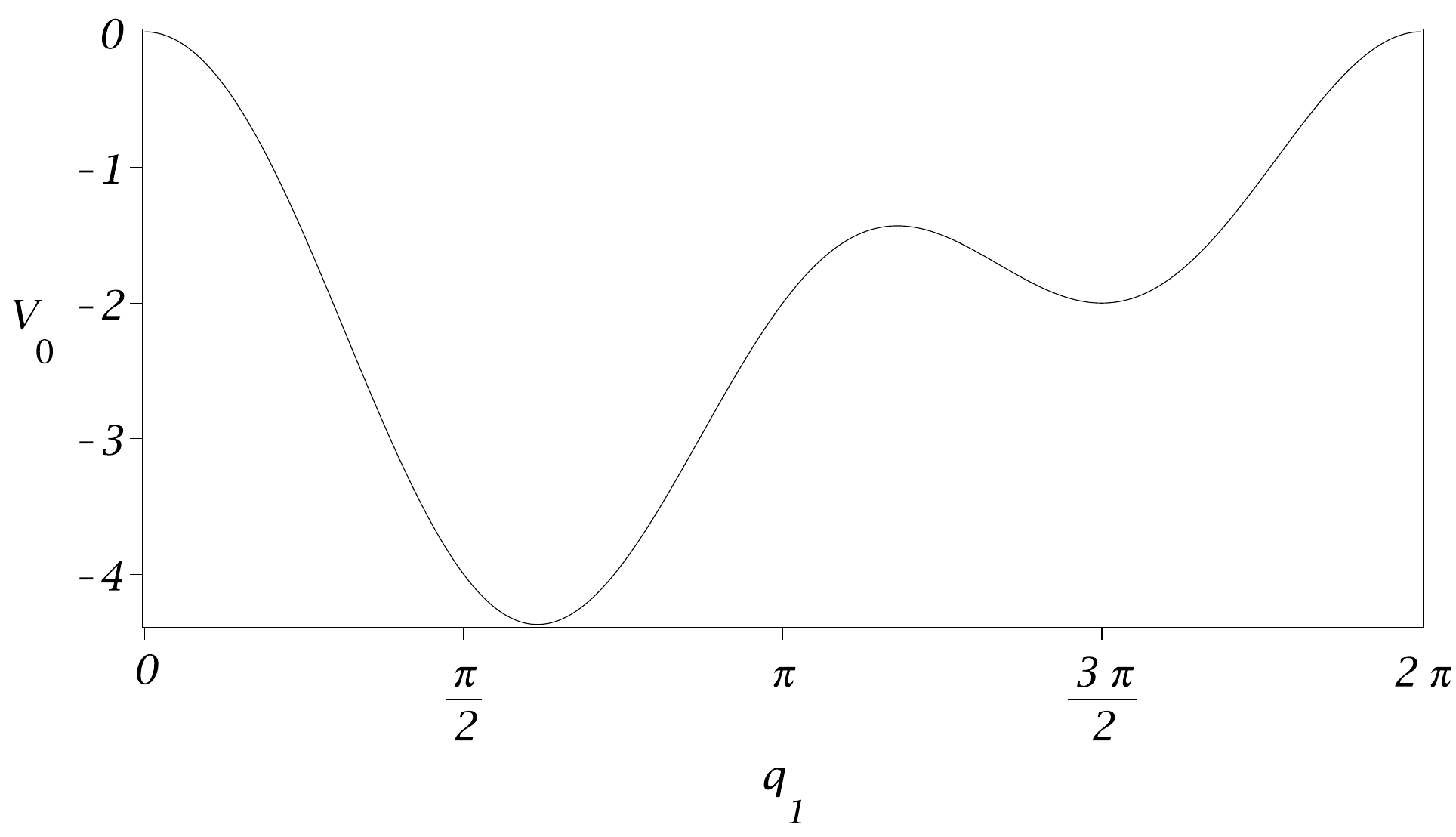}}
\hspace{-0.02\textwidth}
\subfigure[\label{fig:DS0}]{
\includegraphics[width=0.55\textwidth]{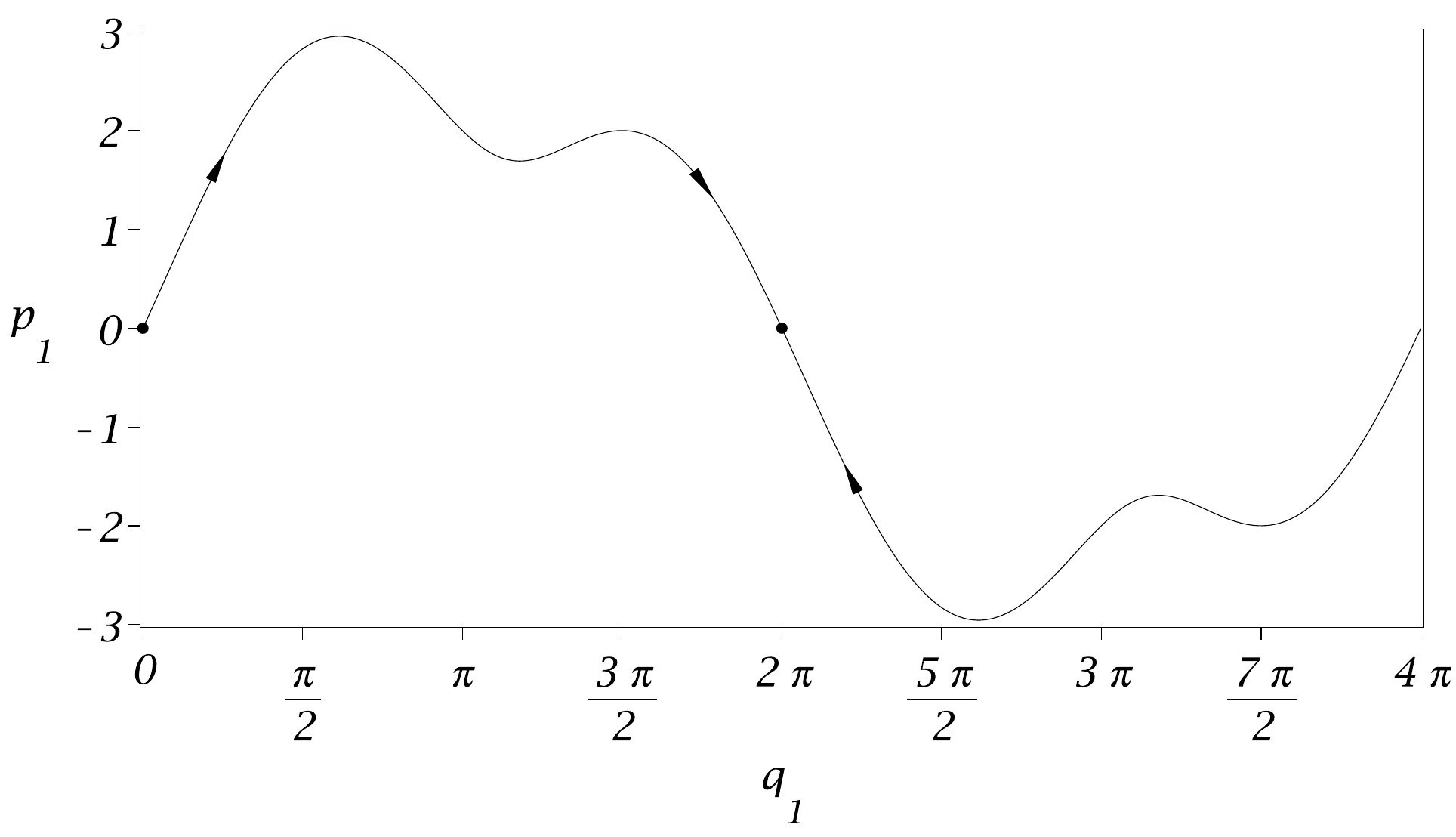}}
\caption{\small
  (a) \emph{Graph of the $2\pi$-periodic potential
        $V_0(q_1)=\cos q_1-\sin q_1+\cos2q_1-\frac12\sin2q_1-2$.}
  (b) \emph{Projection onto the plane $(q_1,p_1)$ of the loops, given by
        the graph $p_1=S'_0(q_1)$ and assuming $\beta(q_1)=1$
        (see Proposition~\ref{teorestrictions}(b));
        the loop $\gamma$ has $q_1$ increasing from $0$ to $2\pi$,
        and a second loop obtained from the periodicity has $q_1$ decreasing
        from $4\pi$ to $2\pi$.}
}
\end{figure}

We point out that this period $4\pi$ is due to the fact that
the functions $S'_0(q_1)$ and $S_1(q_1)$ parameterize two loops of~$O$.
Indeed, the loop $\gamma$ is obtained for $0<q_1<2\pi$ and satisfies $p_1>0$,
and another loop is obtained for $2\pi<q_1<4\pi$ and satisfies $p_1<0$
(see an example in Figures~\ref{fig:V0} and~\ref{fig:DS0}).

Now it is clear, by integrating from 0 to $2\pi$
the first equation in~(\ref{eq:DS0S1}), that we can write
\[
  S_0(q_1):=S^\ut_0(q_1)=\wh S^\st_0(q_1)=S^\st_0(q_1-2\pi)+\sigma,
\]
where $\sigma$ is the constant in~(\ref{eq:S-hat}), given by
\beq\label{sigma}
  \sigma
  =\int_0^{2\pi}p_1\,\df q_1
  =S_0(2\pi)-S_0(0)
  =S^\ut_0(2\pi)
  =-S^\st_0(-2\pi),
\eeq
and we have taken into account our choice
$S^{\ut,\st}_0(0)=0$ in~(\ref{expans1}).

\paragr{The generating functions of the invariant manifolds}
Let us return to the invariant manifolds $\W^\ut$ and~$\wt\W^\st$,
described in the coordinates $(q,p)$ by their generating functions
$S^\ut(q_1,q_2)$ and $\wh S^\st(q_1,q_2)$.
We know from the previous paragraph the relations concerning the coefficients
of orders~0 and~1 in their expansions in $q_2$,
and now we are interested in the coefficients of order~2,
which are used in order to study the transversality.

Expanding in $q_2$ the identity~(\ref{eq:S-hat}), the coefficients
of order~2 satisfy the equality
\[
   \wh{T}^{\st}(q_{1})= T^{\st}(q_{1}-2\pi).
\]
By hypothesis~\hyp{H4\bprime}, the functions $T^\ut$ and $\wh T^{\st}$
are defined for $q_1\in\po{-a,\pi+a}$ and
$q_1\in\po{\pi-a,2\pi+a}$ respectively.
Then, they can be compared for $q_1\in\po{\pi-a,\pi+a}$,
which makes a difference with the fact that the coefficients
of orders~0 and~1 (giving only the loop in~(\ref{eq:DS0S1}))
are functions defined for all $q_{1}\in\R$.

\paragr{Reversibility relations and transversality}
Finally, we assume the reversibility of the Hamiltonian
in the sense stated in Section~\ref{secreversible},
i.e.~with respect to an involution $\Rc$ of the form~(\ref{revers})
with the matrix $R=\diag[r_1,r_2]$ given in~(\ref{matrixrevers}).
Since we want to compare $\W^{\ut}$ with $\wt{\W}^{\st}$ along $\gamma$,
we shall fix $r_{1}=-1$, for this will be the case of interest:
the loop is mapped into itself by the reversibility, $\Rc\gamma=\gamma$,
since it reverses the sign of $q_{1}$, but keeps the sign of $p_{1}$.
In other words, for $r_1=-1$ the restriction of $\Rc$
to the plane $(q_1,p_1)$ is a symmetry with respect to $q_1=0$
(or equivalently $q_1=\pi$).

\begin{proposition}\label{prop:rev-T2}
Under the hypothesis of this section,
if the Hamiltonian~(\ref{ham}) is $\Rc$-reversible with $r_{1}=-1$
(and $r_2=\pm1$), then the following equality holds:
\begin{equation}\label{eq:Shat-Su}
  \wh{S}^{\st}(2\pi-q_{1}, r_{2}q_{2})= -S^{\ut}(q_{1},q_{2})+\sigma,
  \qquad
  (q_1,q_2)\in\po{\pi-a,\pi+a}\times\po{-\delta,\delta},
\end{equation}
where $\sigma$ is the constant in~(\ref{eq:S-hat}).
Besides, its coefficients satisfy
\beq\label{eq:ST-rev}
  \begin{array}{lcll}
   S_{0}(2\pi-q_{1}) &= &-S_{0}(q_{1})+\sigma,&\qquad q_1\in\R,
  \\[4pt]
   S_{1}(2\pi-q_{1}) &= &-r_{2}S_{1}(q_{1}),&\qquad q_1\in\R,
  \\[4pt]
   \wh{T}^{\st}(2\pi-q_{1})&= &-T^{\ut}(q_{1}),
   &\qquad q_1\in\po{\pi-a,\pi+a}.
  \end{array}
\eeq
\end{proposition}

\proof
We use the reversibility relations given
by Proposition~\ref{prop:reversibility}, with $r_1=-1$.
The first equality~(\ref{eq:Shat-Su}) follows from~(\ref{eq:S-hat}) and the
first relation in Proposition~\ref{prop:reversibility}:
\[
  \wh S^\st(2\pi-q_1,r_2q_2)=S^\st(-q_1,r_2q_2)+\sigma=-S^\ut(q_1,q_2)+\sigma.
\]
Expanding this equality in $q_2$, we get the remaining
three equalities~(\ref{eq:ST-rev}).

However, to ensure the validity for $q_1\notin\po{\pi-a,\pi+a}$ of the first
two equalities~(\ref{eq:ST-rev}), we should deal with the whole loop.
For a given point $(q_1,0,p_1,p_2)\in\gamma$, we have in~(\ref{eq:DS0S1})
the equations $p_1=S'_0(q_1)$ and $p_2=S_1(q_1)$. Since $\Rc$ takes this point
into another point of $\gamma$, we also have
$p_1=S'_0(2\pi-q_1)$ and $-r_2p_2=S_1(2\pi-q_1)$.
Then, the functions $S_0$ and $S_1$ satisfy the identities
$S'_0(2\pi-q_1)=S'_0(q_1)$ and $S_1(2\pi-q_1)=-r_2S_1(q_1)$.
The second one corresponds to our statement,
and integrating the first one we obtain $S_0(2\pi-q_1)=-S_0(q_1)+\const$,
where for $q_1=0$ we see that the constant is $\sigma$,
in agreement with~(\ref{sigma}).
\qed

As a direct consequence of this proposition, we see that in this case $Q=\T^2$
the transversality condition~(\ref{transverse}) becomes very simple.
Using the third equality~(\ref{eq:ST-rev}) at $q^*_1=\pi$,
we get that the necessary and sufficient condition for the transversality
of the invariant manifolds $\W^\ut$ and $\wt\W^\st$ along the loop $\gamma$,
can be expressed as follows:
\[
   T^{\ut}(\pi)\ne 0.
\]

\subsection{Example: two identical connected pendula}\label{secpend1}

As an example with the configuration manifold $Q=\T^2$, we consider
the system of two identical pendula connected by an interacting potential.
This model was considered in \cite{GelfreichS95},
where the transversality of the invariant manifolds along a loop
was established by using the variational equations around it
(for other properties of this model, see for instance \cite[\S23D]{Arnold89}).
Here, our aim is to apply the method developed in this paper to
recover the result of transversality in a somewhat more general model.

Initially, in symplectic coordinates
$(\xi,\eta)=(\xi_1,\xi_2,\eta_1,\eta_2)\in\T^2\times\R^2$
we consider the Hamiltonian
\beq\label{identpend}
  H(\xi,\eta)
  =\tfrac12(\eta_1^{\,2}+\eta_2^{\,2})+(\cos\xi_1-1)+(\cos\xi_2-1)
   +f(\xi_1)(1-\cos(\xi_2-\xi_1)),
\eeq
where the interacting term is given by
a $2\pi$-periodic \emph{even} function $f(\xi_1)$, such that $0\le f(0)<1/2$.
Under this inequality, the origin is a hyperbolic equilibrium point
with (different) Lyapunov exponents $\pm1$, $\pm b$,
where we define $b=\sqrt{1-2f(0)}$, $0<b\le 1$.
The case of a constant function $f$, considered in \cite{GelfreichS95},
corresponds to a linear spring connecting the two pendula.
We point out that the term $1-\cos(\xi_2-\xi_1)$ allows us to keep
the periodicity in $\xi_1$, but could be replaced
for instance by $(\xi_2-\xi_1)^2/2$.

It is clear that, if we assume $f\equiv0$ in~(\ref{identpend}),
the system is separable and integrable:
it consists of the product of two identical pendula
with coordinates $(\xi_j,\eta_j)$, $j=1,2$.
We have, for each pendulum, the (positive) separatrix
\beq\label{separatrix1}
  \xi_j=4\arctan\ee^t,
  \qquad
  \eta_j=2\sin\frac{\xi_j}2=\frac2{\cosh t}.
\eeq
Combining the two separatrices, with a free initial condition
for one of them, say the first one,
we have the following 1-parameter family of loops indexed by $s\in\R$
(see Figure~\ref{fig:pend1}),
\beq\label{loops2pend}
  \begin{array}{lll}
    \bar\gamma_s:\qquad
    &\bar\xi_1(t,s)=4\arctan\ee^{s+t},
    &\bar\xi_2(t)=4\arctan\ee^t,
  \\[2pt]
    &\bar\eta_1(t,s)=\dfrac2{\cosh(s+t)},
    &\bar\eta_2(t)=\dfrac2{\cosh t}.
  \end{array}
\eeq

If we consider $f\not\equiv0$ in~(\ref{identpend}),
in general the family of loops is destroyed, but the concrete loop
$\gamma=\bar\gamma_0$
remains unchanged, as well as the dynamics along it.
Our aim is to provide a sufficient condition on the function $f$,
ensuring that the invariant manifolds intersect transversely along this loop.

In order to apply the results of Section~\ref{secriccati0}, we perform
the symplectic change of coordinates
\beq\label{qchange1}
  \vect{\xi_1}{\xi_2}=\varphi(q)=\vect{q_1}{q_1+q_2},
  \qquad
  \vect{\eta_1}{\eta_2}=\Df\varphi(q)^{-\top}p=\vect{p_1-p_2}{p_2},
\eeq
where it is clear that $\xi=\varphi(q)$ is a well-defined change on $\T^2$.
In the new coordinates $(q,p)$,
the expression of the Hamiltonian $H$ takes the form~(\ref{ham}), with
\[
  B(q)=\symmatrix1{-1}2,
  \quad
  V(q)=(\cos q_1-1)+(\cos(q_1+q_2)-1)+f(q_1)(1-\cos q_2).
\]
It is easy to see from Lemma~\ref{teoreverscond} that $H$ is $\Rc$-reversible,
with $r_1=r_2=-1$ (recall that $f$ is an even function).
Now, the loop $\gamma$ lies on the line $q_2=0$:
\[
  q_1^0(t)=4\arctan\ee^t,
  \quad
  q_2^0(t)=0,
  \quad
  p_1^0(t)=\frac4{\cosh t},
  \quad
  p_2^0(t)=\frac2{\cosh t}.
\]
We know from~(\ref{loopparam})
that the loop is also given by
\[
  q_2=0,
  \quad
  p_1=S_0'(q_1)=4\sin\frac{q_1}2,
  \quad
  p_2=S_1(q_1)=2\sin\frac{q_1}2.
\]
Notice that these are $4\pi$-periodic functions, in agreement with the
considerations of Section~\ref{secperiodic}.

\bremark
As an alternative to~(\ref{qchange1}), we might have chosen as
$\xi=\varphi(q)$ a rotation taking the line $\xi_1=\xi_2$ into $q_2=0$.
This change would lead to a diagonal matrix $B(q)$ and,
although it is not a one-to-one change on $\T^2$,
some arrangements can be done.
However, we chose the change~(\ref{qchange1}) since it is a particular
case of the change considered for the example of Section~\ref{secpend2}.
\eremark

Now, we can write down the Riccati equation~(\ref{riccati})
allowing us to study the transversality.
Developing $V(q)$ in $q_2$, we get:
\[
  V_0(q_1)=2(\cos q_1-1),
  \quad
  V_1(q_1)=-\sin q_1,
  \quad
  Y(q_1)=\cos q_1-f(q_1).
\]
Then, the Riccati equation for the unstable manifold $\W^\ut$ is
\beq\label{riccatiGS}
  2\sin\frac{q_1}2\,(T^\ut)'-2\cos\frac{q_1}2\,T^\ut+2(T^\ut)^2
  =-\p{\sin^2\frac{q_1}2+f(q_1)},
  \qquad
  T^\ut(0)=\frac{1+b}2,
\eeq
where the initial condition comes from Lemma~\ref{teoinitialcond}.
We give below a sufficient condition on $f$, ensuring that
the solution of~(\ref{riccatiGS}) can be extended
up to $q_1=\pi$, and satisfies the transversality condition $T^\ut(\pi)\ne0$.

We provide another useful form for the Riccati equation~(\ref{riccatiGS})
by carrying out the change of variable $x=-\cos(q_1/2)$,
with a further modification in order to remove the linear term.
Writing $T^\ut(q_1)=(\bar T^\ut(x)-x)/2$ and $f(q_1)=\bar f(x)$,
we get the equation
\beq\label{riccatiGS2}
  (1-x^2)\deriv{\bar T^\ut}x+(\bar T^\ut)^2=-2(1-x^2)+(1-2\bar f(x)),
  \qquad
  \bar T^\ut(-1)=b,
\eeq
and the transversality condition becomes $\bar T^\ut(0)\ne0$.

As in Section~\ref{secdevaney} the Riccati equation~(\ref{riccatiGS}) can be
transformed, through the change $\bar T^\ut(x)=\dfrac{1-x^2}{y(x)}\,\dfrac{\df
  y}{\df x}$, into the second-order linear differential equation
\beq\label{legendre}
  (1-x^2)\,\dderiv yx-2x\deriv yx+\p{2-\frac{1-2\bar f(x)}{1-x^2}}y=0.
\eeq

We are going to establish the transversality of the invariant manifolds
along the loop $\gamma$ by checking the transversality condition,
for a wide class of functions $f$ in~(\ref{identpend}).
We start with the simplest case of a constant function,
already considered in \cite{GelfreichS95}.

\begin{lemma}\label{teotransvGS0}
If we consider $1-2\bar f(x)\equiv b^2$ constant, with $0<b<1$,
then the solution of~(\ref{riccatiGS2}) is
\beq\label{solriccatiGS2}
  \bar T^\ut_b(x)=b-\dfrac{1-x^2}{b-x},
\eeq
which can be extended up to the interval $x\in\pco{-1,b}$,
and satisfies $\bar T^\ut_b(0)=b-1/b<0$.
\end{lemma}

\proof
In this case, the second-order equation~(\ref{legendre}) is
a generalized \emph{Legendre equation},
for which a fundamental system of solutions is given by
\[
  y_1(x)=(b+x)\p{\frac{1-x}{1+x}}^{b/2},
  \qquad
  y_2(x)=(b-x)\p{\frac{1+x}{1-x}}^{b/2}=y_1(-x).
\]
We get, from the second one, the solution of~(\ref{riccatiGS2})
as given in~(\ref{solriccatiGS2}).
This solution obviously satisfies the transversality condition,
provided $0<b<1$ (of course, the separable case $b=1$ has to be excluded).
\qed

\begin{proposition}\label{teotransvGS}
If the function $\bar f$ satisfies $0<\bar f(x)<1/2$ for all $x\in[-1,0]$,
then the solution $\bar T^\ut(x)$ of~(\ref{riccatiGS2})
can be extended up to the interval $x\in[-1,0]$,
and satisfies $\bar T^\ut(0)<0$.
\end{proposition}

\proof
The argument is similar to the one given for the example
of Section~\ref{secdevaney}.
With a given $\bar f(x)$, we write~(\ref{riccatiGS2})
in the normalized form $\ds\deriv{\bar T^\ut}x=F(x,\bar T^\ut)$,
where the function $F$ is defined for
$(x,\bar T^\ut)\in\poc{-1,0}\times\R$.
Let $c$ and $d$ be positive constants such that
$0<c^2\le1-2\bar f(x)\le d^2<1$ for any $x$.
Replacing $1-2\bar f(x)$ by $c$ and $d$ in~(\ref{riccatiGS2}),
we define functions $F_c(x,\bar T^\ut)$ and $F_d(x,\bar T^\ut)$,
corresponding to respective differential equations of the type considered
in Lemma~\ref{teotransvGS0}.
We know that their associated solutions $\bar T_c^\ut$ and $\bar T_d^\ut$
are both defined for $x\in\pco{-1,c}$ since $c<d$.
Combining the initial conditions
$\bar T_c^\ut(-1)=c$, $\bar T^\ut(-1)=b$, $\bar T_d^\ut(-1)=d$,
that satisfy $c<b<d$,
with the fact that the inequality $F_c\le F\le F_d$ is true for $x>-1$,
we see that the solution $\bar T^\ut(x)$ is confined between
$\bar T_c^\ut(x)$ and $\bar T_d^\ut(x)$,
and hence it can also be extended up to the
whole interval $\pco{-1,c}$,
and satisfies $\bar T^\ut(0)<\bar T_d^\ut(0)<0$.
\qed

\bremark
If the function $f$ in~(\ref{identpend}) is a trigonometric polynomial,
then the change of variable between $q_1$ and $x$ leads to
a second-order linear differential equation~(\ref{legendre}),
whose coefficients are rational functions of $x$; hence, as in
Section~\ref{secdevaney} one can apply Kovacic's algorithm to see
whether~(\ref{legendre}) is integrable (in the Liouville sense) or not and, when
it is, obtain explicit solutions. In such a case, the transversality condition
could be checked directly.
Nevertheless, applying Proposition~\ref{teotransvGS} we can check
the transversality in many cases,
even when it is not possible to solve~(\ref{legendre}) explicitly.
As a concrete example, choosing $f(\xi_1)=(2-\cos\xi_1)/8$
in~(\ref{identpend}), we have $\bar f(x)=(3-2x^2)/8$,
and applying Kovacic's algorithm (as implemented in Maple) for this case it
turns out that the equation~(\ref{legendre}) cannot be solved explicitly
but, on the other hand, without solving it we see from
Proposition~\ref{teotransvGS} that, for this example,
the transversality condition is fulfilled.
\eremark

\section{The case of a perturbed Hamiltonian}\label{secperturbed}

\subsection{Mel{\cprime}nikov integrals}\label{secmelnikov}

Now, we study the invariant manifolds in a Hamiltonian
depending on a small perturbation parameter $\eps$,
\beq\label{hameps}
  H_\eps(q,p)=H(q,p)+\eps H_*(q,p).
\eeq
As in Section~\ref{sectorus},
we assume that the configuration manifold is $Q=\T^2$.
For the unperturbed case $\eps=0$, we assume that the Hamiltonian $H_0=H$
is in the situation described in Section~\ref{secperiodic}.
In particular, it has a hyberbolic point $O$ at the origin,
with a loop $\gamma$ contained in $q_2=0$.

Concerning the invariant manifolds in a neighborhood of the loop,
we are going to consider two different situations:
\btm
\item[\hyp{A}] the unperturbed invariant manifolds $\W^{\ut,\st}$ intersect
  transversely along the loop $\gamma$;
\item[\hyp{B}] the unperturbed invariant manifolds $\W^{\ut,\st}$ coincide
  along the loop $\gamma$.
\etm
We stress that loops of both types~\hyp{A} and~\hyp{B} may coexist
in the same Hamiltonian $H$
(and each of them can be moved to $q_2=0$ by a specific symplectic change),
as illustrated by the example of Section~\ref{secpend2}.

Notice that a loop $\gamma$ of type~\hyp{A} is \emph{isolated}
in the sense that, if an orbit close to $\gamma$ is a loop, it is not contained
in a small neighborhood of $\gamma$.

On the other hand, if $\gamma$ is a loop of type~\hyp{B}
it belongs to a \emph{$1$-parameter family of loops},
filling the \emph{2-dimensional homoclinic manifold} or \emph{separatrix}
$\W:=\W^\ut=\W^\st$.

For $\eps\ne0$, in both cases~\hyp{A} and~\hyp{B} the loop $\gamma$
will not survive in general.
Our aim is to study whether, for $\eps$ small enough,
there exists a \emph{perturbed loop $\gamma_\eps$}
at distance $\Ord(\eps)$ from $\gamma$,
and the perturbed manifolds $\W_\eps^\ut$ and $\W_\eps^\st$
intersect transversely along $\gamma_\eps$,
inside the energy level containing them.
We show that, in the case~\hyp{A}, this problem is quite simple
and can be solved by applying the implicit function theorem on a transverse
section $q_1=\const$.
Instead, in the case~\hyp{B} we have a problem of
\emph{splitting of separatrices}, which must be re-scaled in order to apply
the implicit function theorem,
leading to a necessary and sufficient condition which can be expressed
in terms of Mel{\cprime}nikov integrals.
In other words, we have for~\hyp{A} and~\hyp{B} a \emph{regular} and
a \emph{singular} perturbation problem respectively.

It is clear that for $\eps$ small enough we have a perturbed hyperbolic
equilibrium point $O_\eps$ inherited from $O$,
whose invariant manifolds $\W_\eps^{\ut,\st}$
can be expressed in terms of perturbed generating functions:
$p=\nabla S_\eps^\ut(q)$ and $p=\nabla\wh S_\eps^\st(q)$
(see Section~\ref{secgenerating-loop} for the notations).
We point out that the identity~(\ref{eq:S-hat})
also applies here to the perturbed case:
$\wh S^\st_\eps(q_1,q_2)=S^\st_\eps(q_1-2\pi,q_2)+\sigma_\eps$.

Let us consider the expansion in $\eps$ of the generating functions,
\beq\label{perturbedS}
  S_\eps^\ut(q)=S^\ut(q)+\eps S_*^\ut(q)+\Ord(\eps^2),
  \qquad
  \wh S_\eps^\st(q)=\wh S^\st(q)+\eps\wh S_*^\st(q)+\Ord(\eps^2).
\eeq
On the other hand, we know from the ideas introduced
in Section~\ref{secgenerating-loop} that the transversality of the
perturbed manifolds, inside the energy level containing them,
can be studied from the expansion in $q_2$, up to order~2,
of the generating functions~(\ref{perturbedS}).
For the unstable manifold, we write
\[
  S_\eps^\ut(q)
  =S_{0,\eps}^\ut(q_1)+S_{1,\eps}^\ut(q_1)\,q_2
   +\tfrac12T_\eps^\ut(q_1)\,q_2^{\,2}+\Ord(q_2^{\,3}),
\]
with
\bean
  &&S_{j,\eps}^\ut(q_1)=S_j(q_1)+\eps S_{j,*}^\ut(q_1)+\Ord(\eps^2),
  \quad j=0,1,
\\
  &&T_\eps^\ut(q_1)=T^\ut(q_1)+\eps T_*^\ut(q_1)+\Ord(\eps^2),
\eean
and analogously for $\wh S_\eps^\st(q)$.
Notice that, for $\eps=0$, we have written
$S_j(q_1):=S_j^\ut(q_1)=\wh S_j^\st(q_1)$, $j=0,1$,
since these terms coincide due to the existence
of the unperturbed loop $\gamma$
(recall the definitions~(\ref{eq:defS0S1})).

Following a standard terminology (see for instance \cite{DelshamsG00}),
we define the \emph{splitting potential} as the difference between
the generating functions for the unstable and stable manifolds:
\beq\label{splpot}
  \Lc_\eps(q)=\Delta S_\eps(q):=S_\eps^\ut(q)-\wh S_\eps^\st(q)
\eeq
(for any functions $f^\ut$ and $\wh f^\st$,
we are using the notation $\Delta f=f^\ut-\wh f^\st$).
We point out that the name `splitting potential' is normally used
in the context of splitting of separatrices, i.e.~our case~\hyp{B},
but we are considering here a somewhat more general situation.
Expanding $\Lc_\eps$ in $q_2$, we write
\[
  \Lc_\eps(q)=\Delta S_{0,\eps}(q_1)+\Delta S_{1,\eps}(q_1)\,q_2
   +\tfrac12\Delta T_\eps(q_1)\,q_2^{\,2}+\Ord(q_2^{\,3}).
\]
If we consider the expansion $\Lc_\eps=\Lc+\eps\Lc_*+\Ord(\eps^2)$,
we have
\bean
  &&\Lc(q)=\tfrac12\Delta T(q_1)\,q_2^{\,2}+\Ord(q_2^{\,3}),
  \\[2pt]
  &&\Lc_*(q)=\Delta S_{0,*}(q_1)+\Delta S_{1,*}(q_1)\,q_2
    +\tfrac12\Delta T_*(q_1)\,q_2^{\,2}+\Ord(q_2^{\,3}),
\eean
where we used that $\Delta S_0(q_1)=\Delta S_1(q_1)=0$,
by the existence of the unperturbed loop $\gamma$.

Now, we can rewrite the two hypotheses on the loop,
in terms of the unperturbed generating functions.
Recalling the transversality condition~(\ref{transverse}) for the case~\hyp{A},
we have:
\btm
\item[\hyp{A}] $\Delta T(q_1)\ne0$ for any $q_1$;
\item[\hyp{B}] $\Delta S(q)=0$ for any $q=(q_1,q_2)$.
\etm

\begin{theorem}\label{teoloopsAB}
\ \btm
\item[\rm(a)]
If $\gamma$ is a loop of type~\hyp{A},
then for $\eps$ small enough
there exists a perturbed loop $\gamma_\eps$ biasymptotic to $O_\eps$,
and the invariant manifolds $\W_\eps^{\ut,\st}$
intersect transversely along $\gamma_\eps$.
\item[\rm(b)]
If $\gamma$ is a loop of type~\hyp{B} and the conditions
\beq\label{transvcond2}
  \pd{\Lc_*}{q_2}(\pi,0)=\Delta S_{1,*}(\pi)=0,
  \qquad
  \pdd{\Lc_*}{q_2}(\pi,0)=\Delta T_*(\pi)\ne0
\eeq
are fulfilled,
then for $\eps\ne0$ small enough
there exists a perturbed loop $\gamma_\eps$ biasymptotic to $O_\eps$,
and the invariant manifolds $\W_\eps^{\ut,\st}$
intersect transversely along $\gamma_\eps$.
\etm
\end{theorem}

\proof
As in Section~\ref{secgenerating-loop},
the energy level $\Nc_\eps$ containing both perturbed invariant
manifolds is given by the equation $H_\eps-H_\eps(O_\eps)=0$.
Since for $\eps=0$ we have $\pd H{p_1}=\dot q_1\ne0$ on $\gamma$,
we see from the implicit function that,
near $\gamma$ and for $\eps$ small enough, we can parameterize $\Nc_\eps$ as
$p_1=g_\eps(q_1,q_2,p_2)$. In the coordinates $(q_1,q_2,p_2)$ of $\Nc_\eps$,
the intersections between the invariant manifolds
$\W_\eps^\ut$ and $\W_\eps^\st$ are given by the equation
\[
  p_2=\pd{S_\eps^\ut}{q_2}(q_1,q_2)=\pd{\wh S_\eps^\st}{q_2}(q_1,q_2),
  \qquad
  \textrm{i.e.}
  \quad
  \pd{\Lc_\eps}{q_2}(q_1,q_2)=0.
\]
It is enough to consider the intersections in a transverse section,
say $q_1=\pi$, and solve the equation for $q_2$.
Expanding in $q_2$, we see that the intersections are given by the solutions of
\beq\label{inters}
  \pd{\Lc_\eps}{q_2}(\pi,q_2)
  =\Delta S_{1,\eps}(\pi)+\Delta T_\eps(\pi)\,q_2+\Ord(q_2^{\,2})=0.
\eeq
Now, expanding in $\eps$ and recalling that $\Delta S_1=0$,
the equation becomes
\beq\label{inters1}
  \p{\Delta T(\pi)\,q_2+\Ord(q_2^{\,2})}
  +\eps\p{\Delta S_{1,*}(\pi)+\Delta T_*(\pi)\,q_2+\Ord(q_2^{\,2})}
  +\Ord(\eps^2)=0.
\eeq
Therefore, in the case~\hyp{A} there is a solution $q_2=\Ord(\eps)$, and hence
the invariant manifolds intersect transversely along a perturbed loop
$\gamma_\eps$, for $\eps$ small enough.

In the case~\hyp{B}, the fact that $\Lc(q)=\Delta S(q)=0$
for any $q$ implies that
$\ds\pd\Lc{q_2}(\pi,q_2)=\Delta T(\pi)\,q_2+\Ord(q_2^{\,2})=0$,
and hence the expansion~(\ref{inters1}) begins with terms of order~1 in $\eps$.
Then, dividing by $\eps$ the whole equation we obtain:
\[
  \Delta S_{1,*}(\pi)+\Delta T_*(\pi)\,q_2+\Ord(q_2^{\,2})+\Ord(\eps)=0.
\]
Since we assumed in~(\ref{transvcond2}) that
$\Delta S_{1,*}(\pi)=0$ and $\Delta T_*(\pi)\ne0$,
we see that there is a solution $q_2=\Ord(\eps)$, and hence
the invariant manifolds intersect transversely along a perturbed loop
$\gamma_\eps$, for $\eps\ne0$ small enough.
\qed

\bremark
In the case~\hyp{B}, it is easy to see from the existence of
the perturbed loop $\gamma_\eps$ that $\Delta S_{1,*}(q_1)=0$ for any~$q_1$.
Indeed, parameterizing $\gamma_\eps$ by $q_1$, we write $q_2=f_\eps(q_1)$.
It is satisfied that $\ds\pd{\Lc_\eps}{q_2}(q_1,f_\eps(q_1))=0$.
Expanding this equation in $\eps$ and using that for $\eps=0$
we have $\Lc=0$ and $f(q_1)=0$, we obtain:
$\Delta S_{1,*}(q_1)=\ds\pd{\Lc_*}{q_2}(q_1,0)=0$.
\eremark

\bremark
We point out that, for a loop of type~\hyp{B},
it is not enough to impose a condition like~(\ref{transverse}),
and we need as an additional condition that the coefficients
$S_{1,*}^\ut(\pi)$ and $\wh S_{1,*}^\st(\pi)$ coincide.
Otherwise, the existence of
a perturbed loop $\gamma_\eps$ close to $\gamma$
could not be established in general.
\eremark

The remaining part of this section concerns the case~\hyp{B},
and our aim is to provide an explicit condition allowing us
to check~(\ref{transvcond2}) in a concrete example.
For $\eps=0$, the separatrix $\W$ is a graph $p=\nabla S(q)$,
with the generating function
$S(q):=S^\ut(q)=\wh S^\st(q)$,
and the \emph{inner dynamics} on $\W$ is given
by the ordinary differential equation
\beq\label{innerdyn}
  \dot q=B(q)\nabla S(q).
\eeq
For a given $q=(q_1,q_2)$, denoting $\check q(t,q)$
the solution of~(\ref{innerdyn}) that satisfies the initial condition
$\check q(0,q)=q$, we have a loop
\beq\label{loops}
  \check x(t,q)=(\check q(t,q),\check p(t,q)),
  \qquad
  \textrm{with \ $\check p(t,q)=\nabla S(\check q(t,q))$}.
\eeq
Of course, it will be enough to consider the initial condition
in a direction transverse to $q_2=0$:
\beq\label{kappa}
  q=\kappa(s)=(\kappa_1(s),\kappa_2(s)),\ \ s\in\po{-a,a},
  \qquad
  \textrm{with \ $\kappa(0)=(\pi,0)$ \ and \ $\kappa'_2(0)\ne0$.}
\eeq
In this way, we have a 1-parameter family of loops indexed by $s$,
that we denote $\bar\gamma_s$, given by
$\bar q(t,s):=\check q(t,\kappa(s))$.
The particular loop contained in $q_2=0$ is $\bar\gamma_0=\gamma$,
and we write $q^0(t)=(q_1^0(t),0)=\bar q(t,0)$.
Notice that
$(q_1,q_2)=\bar q(t,s)$, $(t,s)\in\R\times\po{-a,a}$,
is a change of parameters for the separatrix $\W$,
in a neighborhood of $\gamma$.

As a possible approach, we could see the functions $T_\eps^\ut(q_1)$ and
$\wh T_\eps^\st(q_1)$ as solutions of Riccati equations as in~(\ref{riccati}),
with initial conditions at $q_1=0$ and $q_1=2\pi$ respectively.
Developing such equations in $\eps$ and taking the terms of order~1,
the functions $T_*^\ut(q_1)$ and $\wh T_*^\st(q_1)$ become
solutions of \emph{linear differential equations},
and hence they can be expressed in terms of integrals.
Comparing such integrals, one would obtain
a condition for the transversality which could be checked.

Nevertheless, similar integrals of Mel{\cprime}nikov type can be obtained,
in a more direct way, by developing in $\eps$ the Hamilton--Jacobi equation
satisfied by the perturbed generating functions~(\ref{perturbedS})
and comparing the unstable and stable ones.
An analogous approach was previously used in \cite{LochakMS03}
(see also \cite{Sauzin01}),
for the case of invariant manifolds of tori associated to a simple resonance,
i.e.~with a 1-d.o.f. hyperbolic part.
We are going to show that a first order approximation for $\eps$ small
of the splitting potential $\Lc_\eps$ defined in~(\ref{splpot})
is given by the \emph{Mel{\cprime}nikov potential}, defined as the integral
\beq\label{melnipot}
  L(q)=-\int_{-\infty}^\infty\pq{H_*(\check x(t,q))-H_*(O)}\,\df t,
\eeq
absolutely convergent due to the fact that
$\check x(t,q)$ tends to $O$, with exponential bounds, as $t\to\pm\infty$.
Essentially, the Mel{\cprime}nikov potential is defined as the integral
of the perturbation $H_*$ on the unperturbed loops.
See also \cite{Treschev94} for another definition of
a Mel{\cprime}nikov integral along the loops contained in a separatrix,
and \cite{DelshamsR97} for the case of a hyperbolic fixed point of
a symplectic map.

\begin{theorem}\label{teomelnikov}
The splitting potential is given at first order in $\eps$
by the Mel{\cprime}nikov potential: for some constant $c$,
\[
  \Lc_\eps(q)=\eps(L(q)+c)+\Ord(\eps^2).
\]
\end{theorem}

\proof
Expanding in $\eps$ the generating functions in~(\ref{splpot}),
it is enough to show that the difference of their first order terms
is given by the Mel{\cprime}nikov potential:
\beq\label{melnipot2}
  S_*^\ut(q)-\wh S_*^\st(q)=L(q)+c.
\eeq
As in the proof of Theorem~\ref{teoriccati},
we use that the perturbed generating functions
satisfy the Hamilton--Jacobi equation, but now we expand it in $\eps$.
We start with the generating function associated to the unstable manifold:
for any $q$, we have
\[
  H_\eps(q,\nabla S^\ut_\eps(q))
  =\Hc(q)+\eps\Hc^\ut_*(q)+\Ord(\eps^2)
  =\eps H_*(O)+\Ord(\eps^2)
\]
(to obtain the constant at the right hand side, we expand $H_\eps(O_\eps)$
in $\eps$ and we use that $H(O)=0$ and $\nabla H(O)=0$).
The terms of orders~0 and~1 are
\bean
  &&\Hc(q)=H(q,\nabla S(q)),
  \\
  &&\Hc^\ut_*(q)
    =\scprod{\pd Hp(q,\nabla S(q))}{\nabla S^\ut_*(q)}+H_*(q,\nabla S(q)).
\eean
Since for any $q$ we have $\Hc^\ut_*(q)=H_*(O)$, we deduce from~(\ref{ham})
that $\nabla S^\ut_*(q)$ satisfies the following
\emph{linear partial differential equation\/}:
\beq\label{linearpde}
  \scprod{B(q)\nabla S(q)}{\nabla S^\ut_*(q)}=-[H_*(q,\nabla S(q))-H_*(O)].
\eeq
In the solution of this equation, the \emph{characteristic curves}
play an essential r\^ole. Since they are given by the
ordinary differential equation~(\ref{innerdyn}),
the characteristic curves are projections of the unperturbed
loops~(\ref{loops}) onto the configuration space: $q=\check q(t,q)$.

What follows is a standard argument for
the Poincar\'e--Mel{\cprime}nikov method.
Using that $S^\ut_*(\check q(t,q))$ tends to $S^\ut_*(0,0)$
as $t\to-\infty$ for any $q$, we see that
\bean
  &&S^\ut_*(q)-S^\ut_*(0,0)
  =\int_{-\infty}^0\deriv{}t\pq{S^\ut_*(\check q(t,q))}\,\df t
  =\int_{-\infty}^0
     \scprod{\nabla S^\ut_*(\check q(t,q))}{\dot{\check q}(t,q)}\,\df t
\\
  &&\phantom{S^\ut_*(q)-S^\ut_*(0)}
  =-\int_{-\infty}^0[H_*(\check x(t,q))-H_*(O)]\,\df t,
\eean
where we have used~(\ref{innerdyn}) together with~(\ref{linearpde}).
Proceeding analogously with the generating function associated to the
stable manifold, we have
\[
  \wh S^\st_*(2\pi,0)-\wh S^\st_*(q)
  =-\int_0^\infty[H_*(\check x(t,q))-H_*(O)]\,\df t,
\]
and adding the two expressions we get the identity~(\ref{melnipot2}),
with $c=S^\ut_*(0,0)-\wh S^\st_*(2\pi,0)$.
\qed

We see from this result that,
in terms of the Mel{\cprime}nikov potential,
the transversality condition~(\ref{transvcond2}) becomes
\beq\label{transvcond3}
  \pd L{q_2}(\pi,0)=\Delta S_{1,*}(\pi)=0,
  \qquad
  \pdd L{q_2}(\pi,0)=\Delta T_*(\pi)\ne0.
\eeq

Now, we show in the following easy lemma that the Mel{\cprime}nikov
potential $L$ is constant along each loop or, in other words,
it is a first integral of the inner dynamics
on the separatrix, given by~(\ref{innerdyn}).

\begin{lemma}\label{teofirstint}
The value of $L(\check q(\tau,q))$ does not depend on $\tau$.
\end{lemma}

\proof
It is enough to carry out the change of variable $t\mapsto t-\tau$
in the integral~(\ref{melnipot}), together with using the identity
$\check x(t,\check q(\tau,q))=\check x(t+\tau,q)$.
\qed

As a consequence, the splitting potential can be studied in terms
of the parameter $s$ introduced in~(\ref{kappa}).
With this in mind, and denoting $\bar x(t,s):=\check x(t,\kappa(s))$,
we define the \emph{reduced Mel{\cprime}nikov potential} as
\beq\label{redumelnipot}
  \wt L(s):=L(\kappa(s))
  =-\int_{-\infty}^\infty\pq{H_*(\bar x(t,s))-H_*(O)}\,\df t,
  \quad
  s\in\po{-a,a}.
\eeq
We show in the next result that the transversality condition
can also be written in terms of $\wt L(s)$
(for an illustration, see the example of Section~\ref{secpend2}).

\begin{proposition}\label{teomelnikov2}
If the reduced Mel{\cprime}nikov potential $\wt L(s)$ has
a \emph{nondegenerate critical point} at $s=0$,
then for $\eps\ne0$ small enough
there exists a perturbed loop $\gamma_\eps$ biasymptotic to $O_\eps$,
and the invariant manifolds $\W_\eps^{\ut,\st}$
intersect transversely along~$\gamma_\eps$.
\end{proposition}

\proof
We first show that
\[
  \wt L'(0)=0
  \qquad\Longleftrightarrow\qquad
  \pd L{q_2}(\pi,0)=0
\]
and, when this is satisfied, we also show that
\[
  \wt L''(0)\ne0
  \qquad\Longleftrightarrow\qquad
  \pdd L{q_2}(\pi,0)\ne0.
\]
Then, the transversality condition~(\ref{transvcond3})
is fulfilled when the reduced potential $\wt L(s)$
has a nondegenerate critical point at $s=0$.

To show the two equivalences, we use the identity
\beq\label{Lident}
  L(\bar q(t,s))=\wt L(s),
\eeq
which comes from Lemma~\ref{teofirstint} and the fact that
$\bar q(0,s)=\kappa(s)$. The first derivative with respect to $s$
of the identity~(\ref{Lident}) provides
\[
   \pd L{q_1}(\bar q(t,s))\,\pd{\bar q_1}{s}(t,s)
  +\pd L{q_2}(\bar q(t,s))\,\pd{\bar q_2}{s}(t,s)=\wt L'(s),
\]
Considering $s=0$ and using that $L$ is constant along $q_2=0$,
we get the equality
$\ds\pd L{q_2}(\bar q(t,0))\,\pd{\bar q_2}{s}(t,0)=\wt L'(0)$.
Since by~(\ref{kappa}) we have $\ds\pd{\bar q_2}{s}(0,0)=\kappa'_2(0)\ne0$,
we obtain the first equivalence.
Now, assuming that $\wt L'(0)=0$ we deduce that $\ds\pd L{q_2}(\bar q(t,0))=0$
at least in a neighborhood of $t=0$.
Thus, we have $\ds\pdd L{q_1}(\pi,0)=\pddm L{q_1}{q_2}(\pi,0)=0$,
and hence the second derivative with respect to~$s$ of~(\ref{Lident}),
at $s=t=0$, provides
\[
  \pdd L{q_2}(\pi,0)\,\kappa'_2(0)^2=\wt L''(0),
\]
which implies the second equivalence.
\qed

\bremark
If one defines the \emph{splitting function}
as $\M_\eps(q)=\ds\pd{\Lc_\eps}{q_2}(q)$,
this provides a measure of the distance between the invariant manifolds
in the $p_2$-direction. At first order, this function can be approximated
in terms of the \emph{Mel{\cprime}nikov function} $M(q)=\ds\pd L{q_2}(q)$,
or the \emph{reduced} one $\wt M(s)=\wt L'(s)$, which has to have
a \emph{simple zero} at $s=0$, as the condition for the
transversality of the invariant manifolds.
\eremark

To end this section, we revisit the above results assuming for
the perturbed Hamiltonian a suitable type of reversibility
(see Section~\ref{secreversible}).

\begin{proposition}\label{teoeven}
If the Hamiltonian $H_\eps$ is \emph{$\Rc$-reversible} with $r_1=r_2=-1$,
then $\Lc_\eps(\pi,q_2)$ and $L(\pi,q_2)$ are even functions in $q_2$.
\end{proposition}

\proof
We see from the considerations of Section~\ref{secperiodic},
in particular Proposition~\ref{prop:rev-T2}, that
the generating functions of the invariant manifolds $\W^{\ut,\st}_\eps$
are related by
$\wh S^\st_\eps(2\pi-q_1,-q_2)=-S^\ut_\eps(q_1,q_2)+\sigma_\eps$.
Then, we obtain
\[
  \Lc_\eps(\pi,q_2)=S^\ut_\eps(\pi,q_2)+S^\ut_\eps(\pi,-q_2)+\sigma_\eps,
\]
an even function.
This applies also to the Mel{\cprime}nikov potential $L(\pi,q_2)$,
though one could also check that this is an even function directly
from its definition~(\ref{melnipot}).
\qed

Clearly, the fact that $\Lc_\eps(\pi,q_2)$ is even in $q_2$
implies that the equation~(\ref{inters}) always has $q_2=0$ as a
solution, implying the existence of a perturbed loop $\gamma_\eps$
(whose $q_2$-coordinate coincides with $\gamma$ at $q_1=\pi$ but,
in general, we will have $\gamma_\eps\ne\gamma$).
Nevertheless, an additional condition has to be imposed in order to
ensure the transversality of the invariant manifolds
along the loop~$\gamma_\eps$.

Additionally, if we choose $q=\kappa(s)$ in~(\ref{kappa}) in such a way that
the family of unperturbed loops satisfies $\bar\gamma_{-s}=\Rc\bar\gamma_s$,
then the reduced Mel{\cprime}nikov potential $\wt L(s)$
is an even function in $s$.
Then, it always has a critical point at $s=0$, and we only
have to check its nondegeneracy, $\wt L''(0)\ne0$,
in order to ensure transversality.

\subsection{Example: two different weakly connected pendula}\label{secpend2}

In symplectic coordinates $(\xi,\eta)\in\T^2\times\R^2$
we consider the Hamiltonian
\[
  H_\eps(\xi,\eta)
  =\tfrac12(\eta_1^{\,2}+\eta_2^{\,2})+(\cos\xi_1-1)+\lambda^2(\cos\xi_2-1)
   +\eps(1-\cos(\xi_2-\xi_1)),
\]
where $\eps$ is a small parameter, and $\lambda\ge1$ is a fixed value
(recall that the case $\lambda=1$ has already been considered
in Section~\ref{secpend1}).

\begin{figure}[b!]
\subfigure[\label{fig:pend1}]{
\includegraphics[width=0.5\textwidth]{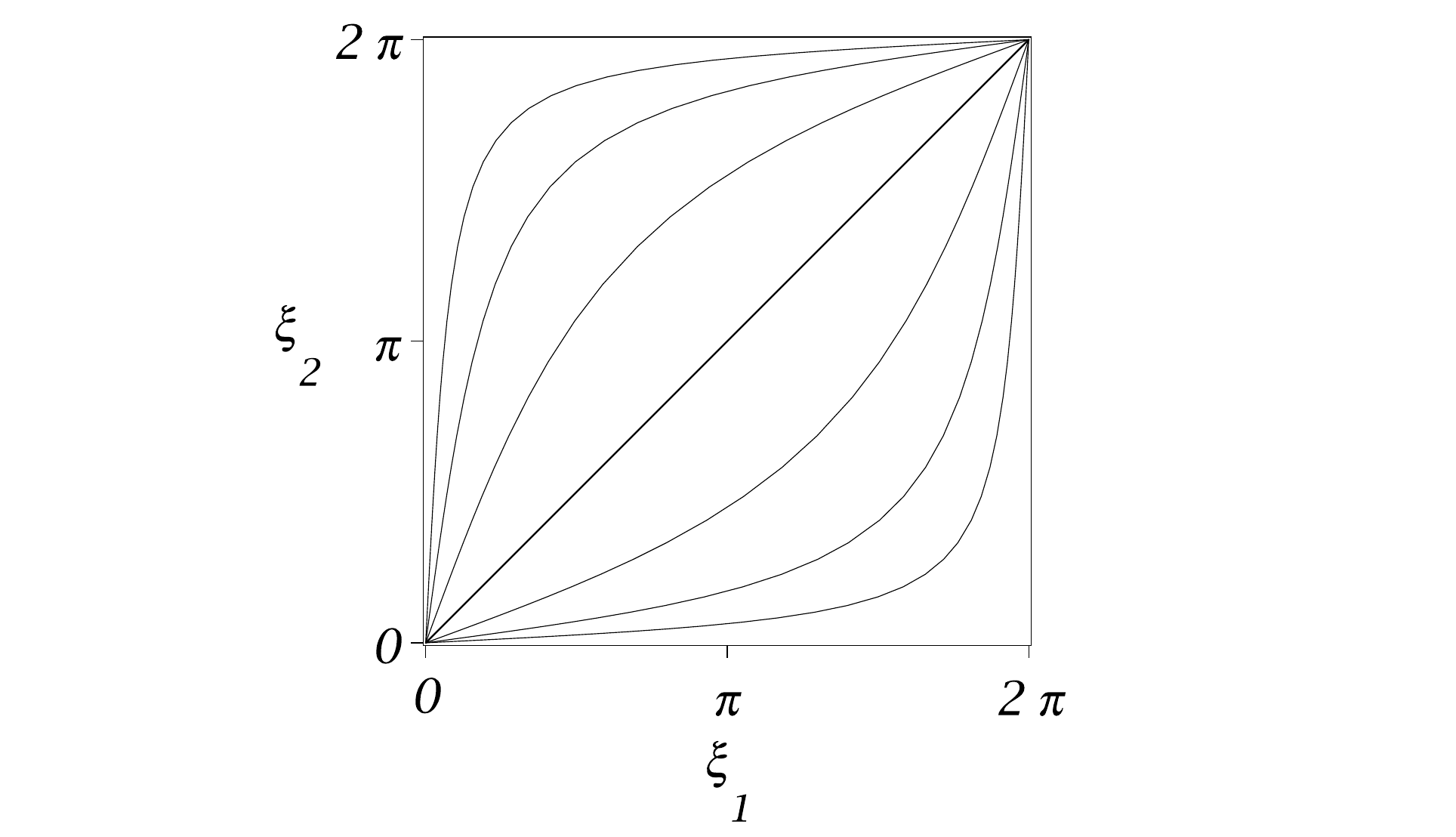}}
\hspace{-0.04\textwidth}
\subfigure[\label{fig:pend2}]{
\includegraphics[width=0.5\textwidth]{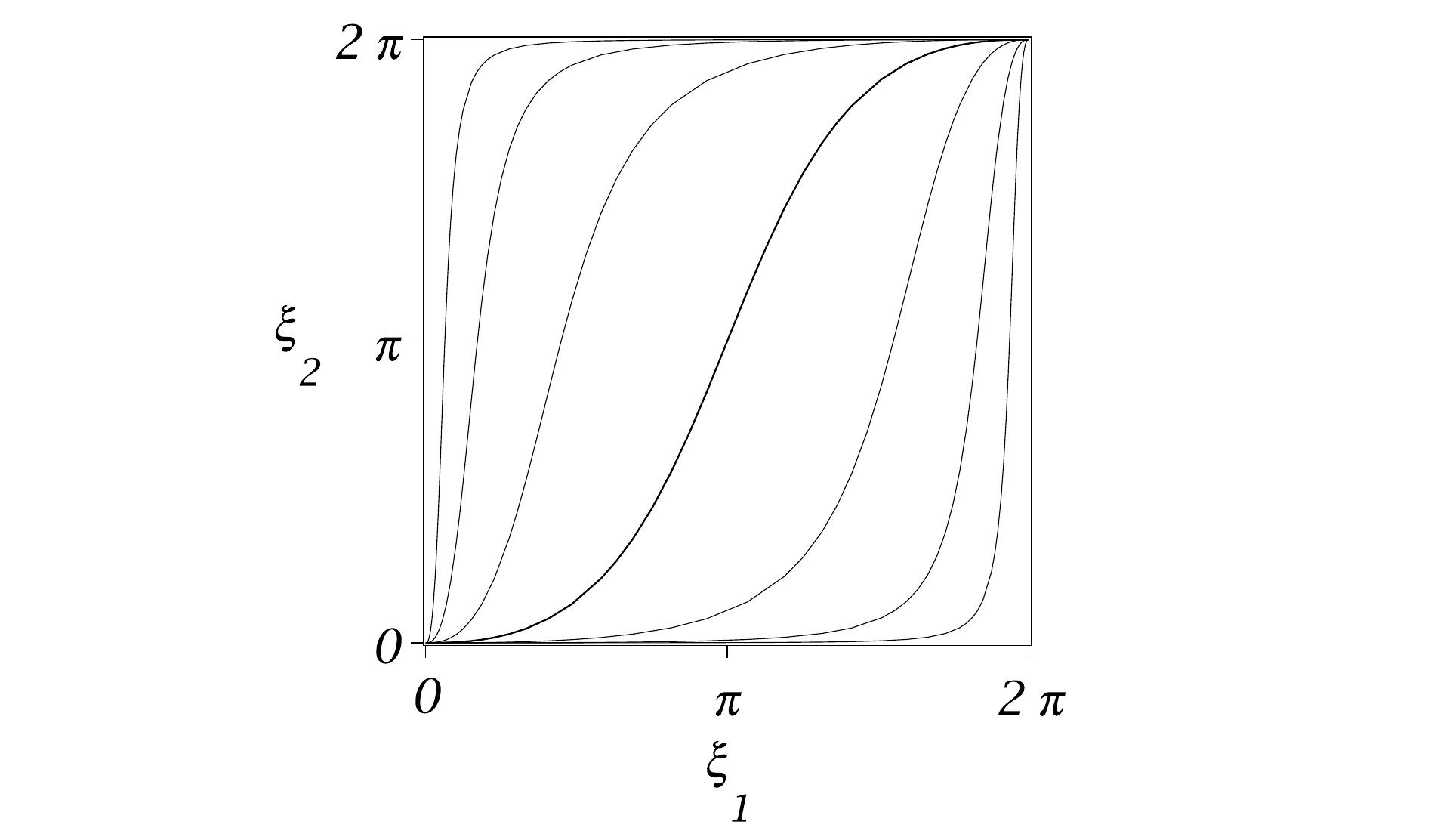}}
\caption{\small
  \emph{The 1-parameter family of loops in the separable case,}
  \ (a) \emph{for $\lambda=1$;}
  \ (b) \emph{for $\lambda=2.5$.}
}
\end{figure}

For $\eps=0$ the system is separable, and consists of
two pendula with Lyapunov exponents $\pm1$ and $\pm\lambda$,
generalizing the separable case of Section~\ref{secpend1}.
On the region $\eta_1,\eta_2\ge0$ (the other three ones are symmetric),
it has a 1-parameter family of loops $\bar\gamma_s$, $s\in\R$,
plus two special loops $\wh\gamma_1$ and $\wh\gamma_2$
with a different topological behavior.
The two special loops are given by the separatrix of first/second pendulum,
with the equilibrium point of the second/first pendulum:
\bean
  &&\wh\gamma_1:
  \qquad
  \xi_1(t)=4\arctan\ee^t,
  \quad
  \eta_1(t)=\dfrac2{\cosh t}\,,
  \quad
  \xi_2(t)=\eta_2(t)=0,
\\
  &&\wh\gamma_2:
  \qquad
  \xi_1(t)=\eta_1(t)=0,
  \quad
  \xi_2(t)=4\arctan\ee^{\lambda t},
  \quad
  \eta_2(t)=\dfrac{2\lambda}{\cosh\lambda t}\,;
\eean
and the 1-parameter family is given by the separatrices of the two pendula,
with a free initial condition in one of them
(see Figure~\ref{fig:pend2}):
\[
  \begin{array}{lll}
    \bar\gamma_s:
    \qquad
    &\bar\xi_1=\bar\xi_1(t,s)=4\arctan\ee^{t-s},
    \quad
    &\bar\xi_2=\bar\xi_2(t)=4\arctan\ee^{\lambda t},
  \\[4pt]
    &\bar\eta_1=\bar\eta_1(t,s)=\dfrac2{\cosh(t-s)}\,,
    \quad
    &\bar\eta_2=\bar\eta_2(t)=\dfrac{2\lambda}{\cosh\lambda t}\,.
  \end{array}
\]

It is not hard to see that
the loops $\wh\gamma_1$ and $\wh\gamma_2$ are of type~\hyp{A}.
Indeed, the unperturbed unstable/stable invariant manifolds
of the loop $\wh\gamma_1$ are given by the separatrix of the first pendulum,
times the local unstable/stable curves of the equilibrium point of
the second pendulum:
recalling~(\ref{separatrix1}), their equations are
$\eta_1=2\sin(\xi_1/2)$, $\eta_2=\pm2\lambda\sin(\xi_2/2)$
(with the signs $+$/$-$ for the unstable/stable manifold respectively),
and it is clear that they intersect transversely.
Then, we see as a consequence of Theorem~\ref{teoloopsAB}(a)
that the equilibrium point $O_\eps=O$ has a perturbed loop $\wh\gamma_{1,\eps}$
whose invariant manifolds intersect transversely along it
for $\eps$ small enough
(from the Hamiltonian equations, one sees that for $\eps\ne0$
there is no orbit contained in $\xi_2=0$,
and hence $\wh\gamma_{1,\eps}\ne\wh\gamma_1$).
The same considerations are valid for the loop $\wh\gamma_2$,
since we did not use the fact that $\lambda\ge1$,
and both loops are completely analogous.

Now, we study the loop $\gamma=\bar\gamma_0$,
corrsponding to $s=0$ in the 1-parameter family.
The existence of this family implies that this is a loop of
type~\hyp{B}, whose separatrix $\W$ is the manifold defined
by the equations
$\eta_1=2\sin(\xi_1/2)$, $\eta_2=2\lambda\sin(\xi_2/2)$, $0<\xi_1,\xi_2<2\pi$.
The manifold $\W$ contains the whole family of loops $\bar\gamma_s$, $s\in\R$.

In order to introduce new coordinates $(q,p)$ such that the loop $\gamma$ is
contained in $q_2=0$, we take into account that the $\xi$-projection of
the loop $\gamma$ is the graph
\[
  \xi_2=h(\xi_1):=4\arctan\p{\tan^\lambda\frac{\xi_1}4}.
\]
This function $h:\T\rightarrow\T$ is of class $\Cc^r$,
where $r$ is the integer part of $\lambda$,
and $h$ is analytic if $\lambda$ is integer.
Notice that we can consider $h$ as an odd function, since the equality
$h(2\pi-\xi_1)=2\pi-h(\xi_1)$ is fulfilled.

Thus, we consider the change of coordinates
\bean
  &&\vect{\xi_1}{\xi_2}=\varphi(q)=\vect{q_1}{h(q_1)+q_2},
\\
  &&\vect{\eta_1}{\eta_2}
  =\Df\varphi(q)^{-\top}p=\mmatrix1{-h'(q_1)}01\vect{p_1}{p_2}.
\eean
(notice that $x=\varphi(q)$ is a well-defined change on $\T^2$).
In the new coordinates $(q,p)$ the Hamiltonian takes
the form~(\ref{hameps}) where, for the unperturbed part $H$,
we have
\bea
  \label{Heps0a}
  &&B(q)=\symmatrix1{-h'(q_1)}{1+h'(q_1)^2},
\\[2pt]
  \label{Heps0b}
  &&V(q)=(\cos q_1-1)+\lambda^2(\cos(h(q_1)+q_2)-1),
\eea
and the perturbation $\eps H_*$ is given by
\beq\label{Heps1}
  H_*(q)=1-\cos(h(q_1)-q_1+q_2).
\eeq
Let us write down, for $\eps=0$, the $q$-coordinates of the
unperturbed loops,
\beq\label{loops2}
  \bar\gamma_s:
  \quad
  \bar q_1=\bar\xi_1
  =4\arctan\ee^{t-s},
  \quad
  \bar q_2=\bar\xi_2-h(\bar\xi_1)
  =4\arctan\ee^{\lambda t}-4\arctan\ee^{\lambda(t-s)},
\eeq
and it is clear that the loop $\gamma=\bar\gamma_0$
is now contained in $q_2=0$.
As in~(\ref{kappa}), the initial conditions in~(\ref{loops2})
are given by $\kappa(s)=(4\arctan\ee^{-s},\pi-4\arctan\ee^{-\lambda s})$,
a direction transverse to $q_2=0$ for $s=0$.

As we easily check, the perturbed Hamiltonian $H_\eps$
in~(\ref{Heps0a}--\ref{Heps1}) is $\Rc$-reversible with $r_1=r_2=-1$
(using that the functions $h$ and $h'$ are odd and even respectively).
Then, we see from Proposition~\ref{teoeven} that there exists
a perturbed loop $\gamma_\eps$ close to $\gamma$.

To study the transversality of the perturbed invariant manifolds along
$\gamma_\eps$, we use the reduced Mel{\cprime}nikov potential:
\bea
  \nonumber
  \wt L(s)
  &=&-\int_{-\infty}^\infty H_*(\bar q)\,\df t
  =-\int_{-\infty}^\infty\pq{1-\cos(h(\bar q_1)-\bar q_1+\bar q_2)}\,\df t
\\
  \label{melnikov2}
  &=&-\int_{-\infty}^\infty\pq{1-\cos(\bar\xi_2-\bar\xi_1)}\,\df t,
\eea
where the expressions
for $\bar \xi_1=\bar\xi_1(t,s)$ and $\bar \xi_2=\bar \xi_2(t)$
have been introduced in~(\ref{loops2pend}).

The following result implies,
according to Proposition~\ref{teomelnikov2},
that the invariant manifolds
intersect transversely along the perturbed loop $\gamma_\eps$
for $\eps\ne0$ small enough,
at least for the interval of $\lambda$ considered.

\begin{figure}[b!]
\begin{center}
  \includegraphics[width=0.5\textwidth]{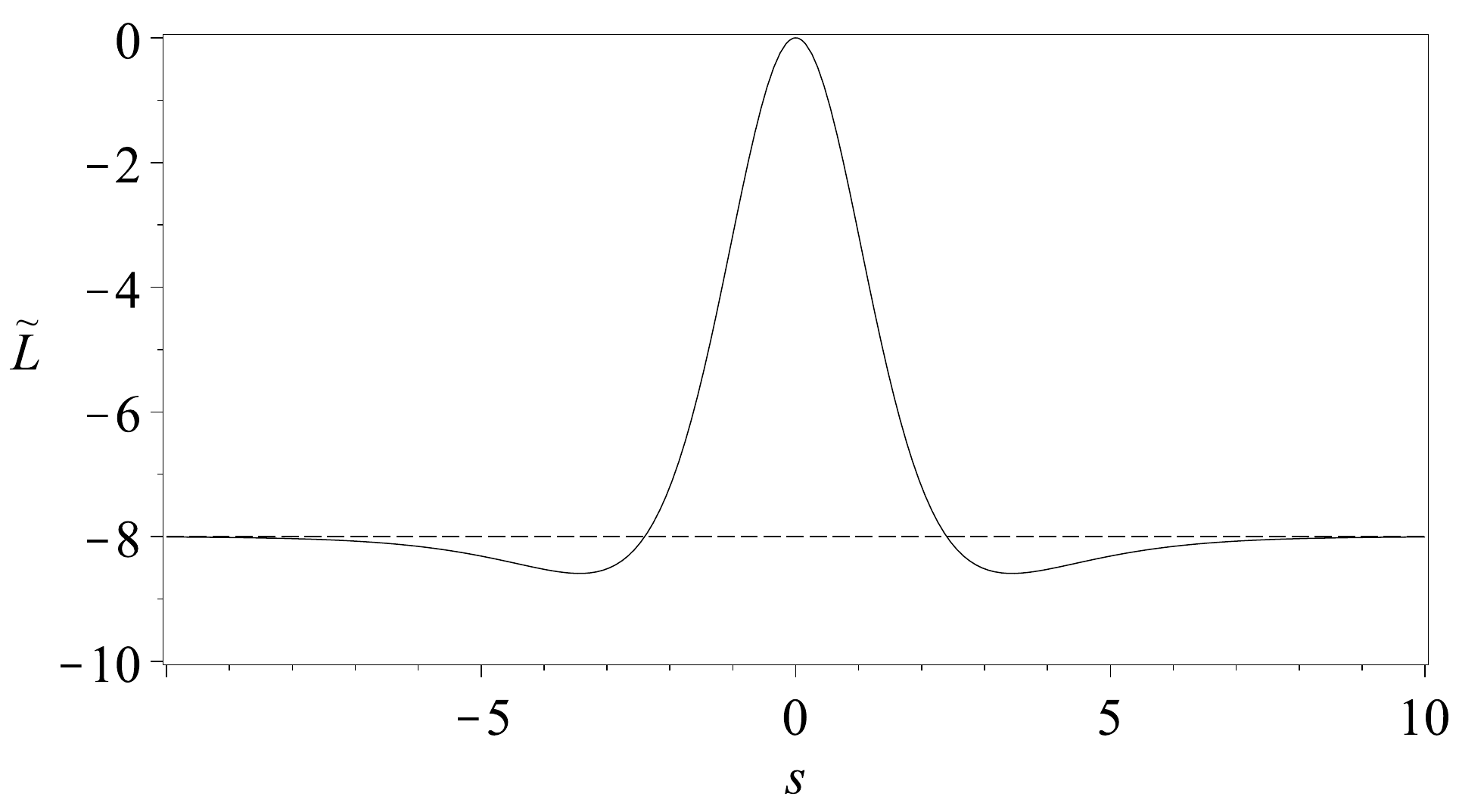}
\end{center}
\caption{\small
  \emph{The reduced Mel{\cprime}nikov potential $\wt L(s)$ for $\lambda=1$.}
}
\label{fig:melni}
\end{figure}

\begin{proposition}\label{teonondeg}
Assuming that $1\le\lambda\le\lambda_0\simeq3.68078$,
the reduced Mel{\cprime}nikov potential~(\ref{melnikov2})
has a nondegenerate critical point at $s=0$.
\end{proposition}

\proof
Differentiating~(\ref{melnikov2}) with respect to $s$, we obtain:
\beq\label{ddL}
  \begin{array}{l}
    \ds\wt L'(s)
    =\int_{-\infty}^\infty\sin(\bar\xi_2-\bar\xi_1)\,\pd{\bar\xi_1}s\,\df t,
  \\
    \ds\wt L''(s)
    =-\int_{-\infty}^\infty
      \pq{\cos(\bar\xi_2-\bar\xi_1)\,\p{\pd{\bar\xi_1}s}^2
          -\sin(\bar\xi_2-\bar\xi_1)\,\pdd{\bar\xi_1}s}\,\df t,
  \end{array}
\eeq
where we have
$\ds\pd{\bar\xi_1}s=-2\sin(\bar\xi_1/2)$ and
$\ds\pdd{\bar\xi_1}s=\sin\bar\xi_1$.

Due to the $\Rc$-reversibility and the fact that
the family of unperturbed loops satisfies $\bar\gamma_{-s}=\Rc\bar\gamma_s$,
it turns out that $\wt L(s)$ is an even function, and hence $\wt L'(0)=0$
(this is a consequence of Proposition~\ref{teoeven},
or may also be checked directly).

Now, let us see that $L''(0)<0$, at least for a large interval of values
of the parameter $\lambda$.
We are going to show that for $s=0$ the function inside the integral
in~(\ref{ddL}),
\[
  f(t)
  =4\cos(\bar\xi_2-\bar\xi_1)\,\sin^2\frac{\bar\xi_1}2
   -\sin(\bar\xi_2-\bar\xi_1)\,\sin\bar\xi_1,
\]
is positive for any $t\in\R$.
Notice that for $s=0$ we have
$\bar\xi_1=4\arctan\ee^t$ and $\bar\xi_2=h(\bar\xi_1)=4\arctan\ee^{\lambda t}$.

For any $t>0$, the difference $\bar\xi_2-\bar\xi_1$ is positive,
and reaches its maximum value $\Xi_\lambda$ when $t$ is such that
$\cosh\lambda t=\lambda\cosh t$.
Using that $\bar\xi_2=h(\bar\xi_1)$ is increasing in $\lambda$ for $t>0$,
then we see that $\Xi_\lambda$ is also increasing from 0 to $\pi$
as $\lambda$ goes from 1 to $\infty$.
Therefore, there exists a value $\lambda_0$ such that $\Xi_{\lambda_0}=\pi/2$,
and numerically we see that $\lambda_0\simeq3.68078$.

If $1\le\lambda\le\lambda_0$, then for $t>0$ we have $\pi<\bar\xi_1<2\pi$ and
$0\le\bar\xi_2-\bar\xi_1\le\pi/2$, and this implies that
$f(t)\ge0$ for any $t>0$.
Similar inequalities show that $f(t)\ge0$ for $t<0$, and we also have $f(0)=4$.
\qed

\bremark
Of course, it is not necessary to have a function $f(t)$ positive for all $t$,
in order to obtain a positive integral.
Numerically, one can see that the result of Proposition~\ref{teonondeg}
is valid for other values of $\lambda$, larger than $\lambda_0$.
\eremark

Finally we point out, for integer values of $\lambda$, the Mel{\cprime}nikov
potential~(\ref{melnikov2}) could be computed explicitly
by writing it as the integral of a rational function.
Applying standard trigonometric and hyperbolic formulas,
and replacing $t\to\frac s2+t$, the integral becomes
\[
  \wt L(s)=-2\int_{-\infty}^\infty
  \pq{\frac{\sinh\p{\frac s2-t}+\sinh\lambda\p{\frac s2+t}}
           {\cosh\p{\frac s2-t}\cdot\cosh\lambda\p{\frac s2+t}}
  }^2\,\df t.
\]
This can be transformed into the integral of a rational function,
if $\lambda$ is integer,
through the change $x=\tanh t$.
This change is analogous to the one carried out in Section~\ref{secpend1},
transforming the Riccati equation~(\ref{riccatiGS}) into
an equation with rational coefficients
(in fact, the change $x=\tanh t$ is the composition of $q_1=4\arctan\ee^t$
and $x=-\cos(q_1/2)$).

For instance, in the simplest case $\lambda=1$ we obtain
\bean
  \wt L(s)
  &=&-8\sinh^2{\tfrac s2}
   \cdot\int_{-1}^1\frac{\df x}{\p{\cosh^2\frac s2-x^2\sinh^2\frac s2}^2}
\\
  &=&-4\tanh{\tfrac s2}\cdot\p{\frac s{\cosh^2\frac s2}+2\tanh\tfrac s2}.
\eean
A numerical inspection of this function (see Figure~\ref{fig:melni})
shows that it has 2~additional critical points,
associated for $\eps\ne0$ small enough to other perturbed loops.

\paragr{Acknowledgments}
We would like to express our sincere thanks to Sergey V.~Bolotin
for his indications about the existence of homoclinic orbits,
and to Primitivo Acosta-Hum\'anez and J.~Tom\'as L\' azaro
for useful discussions and remarks about Riccati equations
and the Kovacic's algorithm.

\small

\def\noopsort#1{}\def\cprime{$'$}

\end{document}